%% file: mes_overview_article_arxiv.tex
\documentclass[a4paper,12pt]{article}
\usepackage{etex}
\usepackage{float} 
\usepackage[off]{auto-pst-pdf}
\usepackage{amsmath,amsfonts,amssymb}
\usepackage{paralist} 
\usepackage{epstopdf}
\usepackage{latexsym}
\usepackage{mathrsfs,amssymb} 
\usepackage[intlimits]{empheq} 
\usepackage[arrow, matrix, curve]{xy}
\usepackage{amsmath}
\usepackage{amssymb}
\usepackage{empheq}
\usepackage{mathrsfs}
\usepackage{lmodern}
\usepackage[noindentafter]{titlesec}
\usepackage[latin1]{inputenc}
\usepackage[T1]{fontenc}
\usepackage{amsthm}
\usepackage{dsfont}
\usepackage{pgf}
\usepackage{graphicx}
\usepackage[nottoc]{tocbibind}
\usepackage[linktocpage=true,colorlinks=false]{hyperref}
\usepackage{xcolor}
\usepackage{pfnote}
\usepackage{fnpos} 
\usepackage{color}
\usepackage{blindtext}
\usepackage{paralist}
\usepackage{makeidx}
\usepackage{a4wide}
\usepackage{pdfpages}
\usepackage[vcentermath]{youngtab}
\usepackage{shuffle}
\usepackage{graphicx}
\usepackage{tikz}
\usepackage{tikz-cd}
\usepackage{stmaryrd}
\usepackage[toc,page]{appendix}

\usetikzlibrary{calc}
\usetikzlibrary{arrows}

\makeindex 
\newlength{\mvlineoffs}

\theoremstyle{definition}
\newtheorem{ex}{\bfseries \upshape Example}[section]
\newtheorem{dfn}[ex]{\bfseries \upshape Definition}

\newtheorem{rem}[ex]{\bfseries \upshape Remark}
\newtheorem{conj}[ex]{\bfseries \upshape Conjecture}
\newtheorem{prop}[ex]{\bfseries \upshape Proposition}
\newtheorem{lem}[ex]{\bfseries \upshape Lemma}
\newtheorem{thm}[ex]{\bfseries \upshape Theorem}
\newtheorem{constr}[ex]{\bfseries \upshape Construction}

\newtheorem{thmx}{Theorem}
\newtheorem{conjx}[thmx]{Conjecture}
\newtheorem{propx}[thmx]{Proposition}
\newtheorem{qu}{Question}
\theoremstyle{plain}

\newtheorem{cor}[ex]{\bfseries \upshape Corollary}

\newenvironment{prf}{\begin{proof}[{\bf Proof}]}{\end{proof}}

\newcommand{\h}{\mathfrak H}
\newcommand{\N}{\ensuremath{\mathds{N}}}	
\newcommand{\Z}{\ensuremath{\mathds{Z}}}	
\newcommand{\Q}{\ensuremath{\mathds{Q}}}	
\newcommand{\R}{\ensuremath{\mathbb{R}}}	
\newcommand{\C}{\ensuremath{\mathds{C}}}
\newcommand{\Ha}{\ensuremath{\mathbb{H}}}

\newcommand{\Li}{\operatorname{Li}}
\newcommand{\Lit}{\widetilde{\operatorname{Li}}}
\newcommand{\dif}{ \operatorname{d}}

\newcommand{\filw}{ \operatorname{Fil}^{\operatorname{W}}}
\newcommand{\fild}{ \operatorname{Fil}^{\operatorname{D}}}
\newcommand{\fille}{ \operatorname{Fil}^{\operatorname{L}} }
\newcommand{\filwle}{ \operatorname{Fil}^{\operatorname{W},\operatorname{L}}}

\newcommand{\thmref}[1] {Theorem \ref{#1}}

\newcommand{\RR}{{\color{red}{R}}}
\newcommand{\UU}{{\color{blue}{U}}}

\newcommand{\sh}{\shuffle}
\newcommand{\mtt}[3] {{\Large \begin{vsmallmatrix} 
  #1\\
 	#2 \\
	#3
\end{vsmallmatrix}}}

\newcommand*{\braceme}[6][]{
\draw[
    shift={(#3:#2)},
    right to reversed-right to reversed,
    shorten >=-.75\pgflinewidth,
    #1
    ] (0,0)
        arc[radius=#2, start angle=#3, end angle=#3+(#4-#3)/2] node[rotate=#3+(#4-#3)/2-90,above=2pt] (#5) {#6};
\draw[
    shift={({#3+(#4-#3)/2}:#2)},
    left to reversed-left to reversed,
    shorten <=-.75\pgflinewidth,
    #1
    ] (0,0)
        arc[radius=#2, start angle=#3+(#4-#3)/2, end angle=#4];
}
\definecolor{lightred}{RGB}{255,114,104}

\DeclareMathOperator{\Sl}{SL}

\DeclareMathOperator{\MD}{\mathcal{MD}}

\DeclareMathOperator{\bMD}{\mathcal{BD}}
\DeclareMathOperator{\bMDG}{\bMD_{gen}}
\DeclareMathOperator{\MDA}{q\mathcal{M}\mathcal{Z}}
\DeclareMathOperator{\MDE}{\mathcal{MD}^\textrm{even}}

\DeclareMathOperator{\MZ}{\mathcal{MZ}}
\DeclareMathOperator{\MZB}{\mathcal{MZB}}

\DeclareMathOperator{\OZ}{\mathsf{Z}}
\DeclareMathOperator{\qMZV}{\mathsf{qMZV}}

\DeclareMathOperator{\ds}{ds}
\DeclareMathOperator{\eds}{\mathsf{eds}}
\DeclareMathOperator{\fds}{\mathsf{fds}}
\DeclareMathOperator{\rds}{\mathsf{rds}}

\DeclareSymbolFont{extraup}{U}{zavm}{m}{n}
\DeclareMathSymbol{\varheart}{\mathalpha}{extraup}{86}
\DeclareMathSymbol{\vardiamond}{\mathalpha}{extraup}{87}
\DeclareRobustCommand{\mb}{\genfrac{[}{]}{0pt}{}}
\DeclareRobustCommand{\mt}{\genfrac{|}{|}{0pt}{}}

\begin{document}
\title{{ \bf Multiple Eisenstein series and\\$q$-analogues of multiple zeta values}}
\author{{\sc Henrik Bachmann}}

\date{\today}
\maketitle
\abstract{This work is an example driven overview article of recent works on the connection of multiple zeta values, modular forms and $q$-analogues of multiple zeta values given by multiple Eisenstein series.}
\tableofcontents

\section*{Introduction}
\setcounter{section}{0}
We study a specific connection of multiple zeta values and modular forms given by multiple Eisenstein series. This work is an example driven overview article and summary of the results obtained in the works \cite{BK},\cite{BT},\cite{Ba2} and \cite{BK2}.

Multiple zeta values are real numbers that are natural generalizations of the Riemann zeta values. These are defined for integers $s_1 \geq 2$ and $s_2,\dots,s_l \geq 1$ by
\[ \zeta( s_1 , \dots , s_l ) :=  \sum_{n_1 >n_2 >\dots>n_l>0 } \frac{1}{n_1^{s_1} \dots n_l^{s_l} } \,. \]
Such real numbers were already studied by Euler in the $l=2$ case in the 18th century. Because of its occurrence in various fields of mathematics and theoretical physics these real numbers had a comeback in the mathematical and physical research community in the late 1990s due to works by several people such as  D. Broadhurst, F. Brown, P. Deligne, H. Furusho, A. Goncharov, M. Hoffman, M. Kaneko, D. Zagier et al.. Denote the $\Q$-vector space of all multiple zeta values of weight $k$ by 
\[ \MZ_k :=\big <\,\zeta(s_1,\dots,s_l) \, \big| \, s_1 + \dots + s_l = k 
\textrm{ and } l>0 \big>_\Q   \]
and write $\MZ$ for the space of all multiple zeta values. One of the main interests is to understand the $\Q$-linear relations between these numbers. The first one is given by $\zeta(2,1) = \zeta(3)$ and there are several different ways to prove this relation (\cite{BB}). Using the representation of multiple zeta values as an ordered sum, their product can be written as a linear combination of multiple zeta values of the same weight, i.e. the space $\MZ$  has the structure of a $\Q$-algebra. For example it is
\begin{align}\label{eq:stuffle1}
 \zeta(2) \cdot \zeta(3) &= \zeta(2,3) + \zeta(3,2) + \zeta(5) \,, \\ \label{eq:stuffle2}
 \zeta(3) \cdot \zeta(2,1) &= \zeta(3,2,1) + \zeta(2,3,1) + \zeta(2,1,3) + \zeta(5,1) + \zeta(2,4) \,. 
\end{align}
This way to express the product, which will be studied in Section $1$ in more detail, is called the stuffle product (also named harmonic product). Besides this, a representation of multiple zeta values as iterated integrals yields another way to express the product of two multiple zeta values, which is called the shuffle product. For the above examples, this is given by 
\begin{align}\label{eq:shuffle1}
\zeta(2) \cdot \zeta(3) &=  \zeta(2,3) + 3 \zeta(3,2) + 6 \zeta(4,1)  \,, \\ \label{eq:shuffle2}
\zeta(3) \cdot \zeta(2,1) &= \zeta(2,1,3)+\zeta(2,2,2)+2\zeta(2,3,1)+2\zeta(3,1,2) +5\zeta(3,2,1)+9\zeta(4,1,1) \,.
\end{align}
Since \eqref{eq:stuffle1} and \eqref{eq:shuffle1} are two different expressions for the product $\zeta(2) \cdot \zeta(3)$ we obtain the linear relation $\zeta(5) = 2 \zeta(3,2) + 6 \zeta(4,1)$. These relations are called the double shuffle relations. Conjecturally all $\Q$-linear relations between multiple zeta values can be proven by using an extended version of these types of relations (\cite{IKZ}). Often relations between multiple zeta values are not proven by using double shuffle relations, since there are easier ways to prove them in some cases. The relation $\zeta(4) = \zeta(2,1,1)$ for example, has an easy proof using the iterated integral expressions for multiple zeta values. A famous result by Euler is, that every even zeta value $\zeta(2k)$ is a rational multiple of $\pi^{2k}$ and in particular we have, for example,
\begin{equation}\label{eq:eulerrelation}
 \zeta(2)^2 = \frac{5}{2} \zeta(4)\,,\quad \zeta(4)^2 =  \frac{7}{6} \zeta(8) \,, \quad \zeta(6)^2 = \frac{715}{691} \zeta(12) \,.
\end{equation}
The relations \eqref{eq:eulerrelation} can also be proven with the double shuffle relations, but for general $k$ there is no explicit proof of Eulers relations using only double shuffle relations so far.

Since the double shuffle relations just give relations in a fixed weight it is conjectured that the space $\MZ$ is a direct sum of the $\MZ_k$, i.e. there are no relations between multiple zeta values with different weight. 

Surprisingly there are several connections of these numbers to modular forms for the full modular group. Recall, modular forms are holomorphic functions in the complex upper half-plane fulfilling certain functional equations. One of the most famous connection is the Broadhurst-Kreimer conjecture. 
\begin{conjx}\label{conj:bk}(Broadhurst-Kreimer conjecture \cite{bk}) The generating series of for the dimension $\dim_\Q\left( \MZ_{k,l} \right)$ of weight $k$ multiple zeta values of length $l$ modulo lower lengths can be written as
\[   \sum_{\substack{k\geq 0 \\ l\geq 0}} \dim_\Q\left( \MZ_{k,l} \right) X^k Y^l = \frac{1 + \mathbb{E}(X) Y}{1 - \mathbb{O}(X) Y + \mathbb{S}(X) Y^2 - \mathbb{S}(X) Y^4} \,, \]
where
\[ \mathbb{E}(X) = \frac{X^2}{1-X^2}  \,,\quad \mathbb{O}(X) = \frac{X^3}{1-X^2} \,, \quad  \mathbb{S}(X) = \frac{X^{12}}{(1-X^4)(1-X^6)} \,.\]
\end{conjx} 
The connection to modular forms arises here, since $\mathbb{S}(X) = \sum_{k\geq 0} \dim S_k(\Sl_2(\Z)) X^k$ is the generating function of the dimensions of cusp forms for the full modular group. In the formula of the Broadhurst-Kreimer conjecture one can see, that cusp forms give rise to relations between double zeta values, i.e. multiple zeta values in the length $l=2$ case. For example in weight $12$, the first weight where non-trivial cusp forms exist, there is the following famous relation
\begin{equation}\label{eq:exotic}
 \frac{5197}{691} \zeta(12) =  168 \zeta(5,7)+150 \zeta(7,5) + 28 \zeta(9,3) \,.
\end{equation}
Even though we are not focused on this conjecture, the concept of obtaining relations of multiple zeta values by cusp forms also appears in our context of multiple Eisenstein series and $q$-analogues of multiple zeta values. It is known that every modular form for the full modular group can be written as a polynomial in classical Eisenstein series. These are for even $k>0$ given by 
\begin{align*}
 G_k(\tau) =\frac{1}{2} \sum_{\substack{(m,n) \in \Z^2 \\ (m,n) \neq (0,0)}} \frac{1}{(m\tau + n)^k }= \zeta(k) + \frac{(-2\pi i)^k}{(k-1)!} \sum_{n = 1}^\infty \sigma_{k-1}(n) q^{n}\,,
\end{align*}
where $\tau \in \Ha$ is an element in the upper half-plane, $q=\exp(2\pi i \tau)$ and $\sigma_k(n) = \sum_{d | n} d^k$ denotes the classical divisor-sum. 
In \cite{GKZ} the authors introduced a direct connection of modular forms to double zeta values following ideas of Don Zagier introduced in \cite{dz2}. They defined double Eisenstein series $G_{s_1,s_2} \in \C[\![q]\!]$ which are a length two generalization of classical Eisenstein series and which are given by a double sum over ordered lattice points. These functions have a Fourier expansion given by sums of products of multiple zeta values and certain $q$-series with the double zeta value $\zeta(s_1,s_2)$ as their constant term. In \cite{Ba} the author treated the multiple cases and calculated the Fourier expansion of multiple Eisenstein series $G_{s_1,\dots,s_l} \in \C[\![q]\!]$. The result of \cite{Ba} was that the Fourier expansion of multiple Eisenstein series is again a $\MZ$-linear combination of multiple zeta values and the $q$-series $g_{t_1,\dots,t_m} \in \C[\![q]\!]$ defined by  $g_{t_1,\dots,t_m}(\tau) := (-2\pi i)^{t_1+\dots+t_m} [t_1,\dots,t_m]$ with $q=e^{2\pi i \tau}$ and
\[ [t_1,\dots,t_m] := \sum_{\substack{u_1 > \dots > u_m> 0 \\ v_1, \dots , v_m >0}} \frac{v_1^{t_1-1} \dots v_m^{t_m-1}}{(t_1-1)!\dots(t_m-1)!}   \cdot q^{u_1 v_1 + \dots + u_m v_m} \,. \]  

\begin{thmx}(\cite{Ba})
For $s_1,\dots,s_l \geq 2$ the $G_{s_1,\dots,s_l}$ can be written as a $\MZ$-linear combination of the above functions $g_{t_1,\dots,t_m}$. 
\end{thmx} 

For example:
\begin{align*}
G_{3,2,2}(\tau) =& \zeta(3,2,2) + \left( \frac{54}{5} \zeta(2,3) + \frac{51}{5} \zeta(3,2) \right) g_2(\tau) + \frac{16}{3} \zeta  (2,2)g_3(\tau) \\
&+3 \zeta(3) g_{2,2}(\tau) + 4 \zeta(2) g_{3,2}(\tau) + g_{3,2,2}(\tau)\,.
\end{align*}

The starting point of the thesis \cite{Ba4} was the fact that there are more multiple zeta values than multiple Eisenstein series, since $\zeta(s_1,\dots,s_l)$ exists for all $s_1 \geq 2, s_2, \dots, s_l \geq 1$ and the $G_{s_1,\dots,s_l}$ just exists when all $s_j \geq 2$. The main objective was to answer the following question

\begin{qu} \label{qu1}
\emph{
 What is a "good" definition of a "regularized" multiple Eisenstein series, such that for each multiple zeta value  $\zeta(s_1,\dots,s_l)$ with $s_1>1$,$s_2,\dots,s_l \geq 1$ there is a $q$-series
\[ G^{reg}_{s_1,\dots,s_l} =  \zeta(s_1,\dots,s_l)+ \sum_{n>0} a_n q^n \in \C[\![q]\!] \]
with this multiple zeta value as the constant term in its Fourier expansion and which equals the multiple Eisenstein series in the cases $s_1,\dots,s_l \geq 2$?}
\end{qu}
By "good" we mean that these regularized multiple Eisenstein series should have the same, or at least as close as possible, algebraic structure similar to multiple zeta values. Our answer to this question was approached in several steps which will be described in the following i)-iii). First (i)  the algebraic structure of the functions $g$ was studied. During this investigation it turned out, that these objects, or more precisely the $q$-series $[s_1,\dots,s_l]$ are very interesting objects in their own rights. It turned out that in order to understand their algebraic structure it was necessary to study a more general class of $q$-series, called bi-brackets in (ii). The results on bi-brackets and brackets then were used, together with a beautiful connection of the multiple Eisenstein series to the coproduct structure of formal iterated integrals, to answer the above question in (iii).  \newline

{\bf i)} To answer Question \ref{qu1} the algebraic structure of the functions $g$ or more precisely the algebraic structure of the $q$-series $[s_1,\dots,s_l]$ was studied in \cite{BK}. It turned out that these $q$-series, whose coefficients are given by weighted sums over partitions of $n$, are, independently to their appearance in the Fourier expansion of multiple Eisenstein series, very interesting objects.  We will denote the $\Q$-vector space spanned by all these brackets and $1$ by $\MD$. Since we also include the rational numbers, the normalized Eisenstein series $\widetilde{G}_k(\tau) := (-2\pi i)^{-k}G_k(\tau)$ are contained in $\MD$. For example we have
\[ \widetilde{G}_2 = -\frac{1}{24} + [2] \,,\quad \widetilde{G}_4 = \frac{1}{1440} + [4] \,, \quad \widetilde{G}_6 = -\frac{1}{60480} + [6] \,.\]
The algebraic structure of the space $\MD$ was studied in \cite{BK} and one of the main result was the following 
\begin{thmx}(\cite{BK}) \label{thm:bk} The $\Q$-vector space spanned by all brackets equipped with the usual multiplication of formal $q$-series is a $\Q$-algebra, with the algebra of modular forms with rational coefficients as a subalgebra. 
\end{thmx} 
In fact, the product fulfills a quasi-shuffle product and the notion of quasi-shuffle products will be made precise in Section \ref{sec:brackets}. Roughly speaking, this means that the product of two brackets can be expressed as a linear combination of brackets similar to the stuffle product \eqref{eq:stuffle1},\eqref{eq:stuffle2}  of multiple zeta values. For example we will see that 
\begin{align*}
 [2] \cdot [3] &= [3,2] + [2,3] + [5] - \frac{1}{12} [3] \,,\\
 [3] \cdot [2,1] &= [3,2,1]+[2,3,1]+[2,1,3]+[5,1]+[2,4]+\frac{1}{12}[2,2]-\frac{1}{2}[2,3]-\frac{1}{12}[3,1] \,,
\end{align*}
i.e. up to the lower weight term  $- \frac{1}{12} [3]$ and $\frac{1}{12}[2,2]-\frac{1}{2}[2,3]-\frac{1}{12}[3,1]$ this looks exactly like \eqref{eq:stuffle1},\eqref{eq:stuffle2}. One might ask if there is also something which corresponds to the shuffle product \eqref{eq:shuffle1} of multiple zeta values. It turned out that for the lowest length case, this has to do with the differential operator $\dif = q \frac{d}{dq}$. In \cite{BK} it was shown that 
\begin{equation}\label{eq:sh23}
 [2] \cdot [3] = [2,3] + 3[3,2] + 6[4,1] - 3[4] + \dif[3] \,,
\end{equation}
which, again up to the term $- 3[4] + \dif[3]$,  looks exactly like the shuffle product \eqref{eq:shuffle1} of multiple zeta values. In particular it follows that $\dif[3]$ is again in the space $\MD$ and in general it was shown that
\begin{thmx} (\cite{BK})\label{thm:derivintro}
 The operator $\dif = q \frac{d}{dq}$ is a derivation on $\MD$.
\end{thmx}

{\bf ii)} Equation \eqref{eq:sh23} above was the motivation to study a larger class of $q$-series, which will be called bi-brackets. While the quasi-shuffle product of brackets also exists in higher length, the second expression for the product, corresponding to the shuffle product, does not appear in higher length if one just allows derivatives as "error terms". The bi-brackets can be seen as a generalization of the derivative of brackets. For $s_1,\dots,s_l \geq 1$, $r_1,\dots,r_l \geq 0$ we define these bi-brackets by
\begin{align*}
\mb{s_1, \dots , s_l}{r_1,\dots,r_l} :=\sum_{\substack{u_1 > \dots > u_l> 0 \\ v_1, \dots , v_l >0}} \frac{u_1^{r_1}}{r_1!} \dots \frac{u_l^{r_l}}{r_l!} \cdot \frac{v_1^{s_1-1} \dots v_l^{s_l-1}}{(s_1-1)!\dots(s_l-1)!}   \cdot q^{u_1 v_1 + \dots + u_l v_l} \in \Q[\![q]\!]\,.
\end{align*}
In the case $r_1=\dots=r_l=0$ these are just ordinary brackets. The products of these seemingly larger class of $q$-series have two representations similar to the stuffle and shuffle product of multiple zeta values in arbitrary length. For our example, the analogue of the shuffle product \eqref{eq:shuffle2} for brackets can now be expressed as 
\begin{align*}
[3] \cdot [2,1] \,&= [2,1,3]+[2,2,2]+2[2,3,1]+2[3,1,2]+5[3,2,1]+9[4,1,1]\\
&+\mb{2,3}{0,1} + 2 \mb{3,2}{0,1} + 3 \mb{4,1}{1,0} - [2,3] - 2 [3,2] -6[4,1]\,. 
\end{align*}
We will see in Section \ref{sec:dsh} that these double shuffle structure can be described, using the so called partition relation, in a nice combinatorial way. This gives a large family of linear relations between bi-brackets. In fact numerical calculations show, that there are so many relations, that we have the following surprising conjecture
\begin{conjx}\label{conj:mdbd}
Every bi-bracket can be written in terms of brackets, i.e. $\MD = \bMD$. 
\end{conjx}
Using the algebraic structure of the space of bi-brackets we now review the definition of shuffle brackets $[s_1,\dots,s_l]^\sh$ and stuffle $[s_1,\dots,s_l]^\ast$ version of the ordinary brackets as certain linear combination of bi-brackets as introduced in \cite{Ba2}. These objects fulfill the same shuffle and stuffle products as multiple zeta values respectively. Both constructions use the theory of quasi-shuffle algebras first developed by Hoffman in \cite{H} and later generalized in \cite{hi}. We summarize the results in the following Theorem. 
 
\begin{thmx}(\cite{Ba2})
\begin{enumerate}[i)]
\item The space $\bMD$ spanned by all bi-brackets $ \mb{s_1, \dots , s_l}{r_1,\dots,r_l}$  forms a $\Q$-algebra with the space of (quasi-)modular forms and the space  $\MD$ of brackets as subalgebras. There are two ways to express the product of two bi-brackets which correspond to the stuffle and shuffle product of multiple zeta values. 
\item There are two subalgebras $\MD^\sh \subset \bMD$ and $\MD^\ast \subset \MD$ spanned by elements $[s_1,\dots,s_l]^\sh$ and $[s_1,\dots,s_l]^\ast$ which fulfill the shuffle and stuffle products, respectively, and which are in the length one case given by the bracket $[s_1]$.
\end{enumerate}
\end{thmx}
For example, similarly to the relation between multiple zeta values above we have 
\[ [2,3]^\ast + [3,2]^\ast + [5] = [2] \cdot [3] = [2,3]^\sh + 3 [3,2]^\sh + 6 [4,1]^\sh \,. \] \newline

{\bf iii)} A particular reason for studying the $[s_1,\dots,s_l]^\sh$ is due to their use in the regularization of multiple Eisenstein series, i.e. they are needed in the answer of the original Question \ref{qu1}. This was implicitly done in \cite{BT} by proving an explicit connection of the Fourier expansion of multiple Eisenstein series to the coproduct on formal iterated integrals introduced by Goncharov in \cite{G}. This connection was already known to the authors of \cite{GKZ} in the length two case. Without knowing this connection it was then rediscovered independently by the authors of \cite{BT} during a research stay of the second author at the DFG Research training Group 1670 at the University of Hamburg in 2014. The result of this research stay was the work \cite{BT}, in which the authors used this connection to give a definition of the shuffle regularized multiple Eisenstein series. Later, the present author combined the result of \cite{BT} and the algebraic structure of bi-brackets to give a more explicit definition of shuffle regularized multiple Eisenstein series using bi-brackets in \cite{Ba2}. 

Formal iterated integrals are symbols $I(a_0; a_1, \dots, a_n; a_{n+1})$ with $a_j \in \{ 0, 1 \}$ that fulfill identities like real iterated integrals. We will write $I(3,2)$ for $I(1;00101;0)$ and we will see that the elements of the form $I(s_1,\dots,s_l)$, obtained in the same way as $I(3,2)$, form a basis of the space of formal iterated integrals in which we are interested. The space of these integrals has a Hopf algebra structure with the multiplication given by the shuffle product and the coproduct $\Delta$ given by an explicit formula which we will review in Section \ref{sec:formaliteratedintegrals}. For example it is
\[ \Delta( I(3,2) ) = 1 \otimes I(3,2)  + 3 I(2) \otimes I(3) + 2  I(3) \otimes I(2) +  I(3,2) \otimes 1  \,.  \]
Compare this with the Fourier expansion of the double Eisenstein series $G_{3,2}$
\[ G_{3,2}(\tau) = \zeta(3,2) + 3 g_2(\tau) \zeta(3) + 2 g_3(\tau) \zeta(2) + g_{3,2}(\tau) \,. \]
Since $\Delta(I(s_1,\dots,s_l))$ exists for all $s_1,\dots,s_l \geq 1$ this comparison suggested a definition of shuffle regularized multiple Eisenstein series $G^\sh_{s_1,\dots,s_l}$ by sending the first component of the coproduct of $I(s_1,\dots,s_l)$ to a $(-2\pi i)$-multiple of the shuffle bracket and the second component to shuffle regularized multiple zeta values. In $\cite{BT}$ it was proven that this construction gives back the original multiple Eisenstein series in the cases $s_1,\dots,s_l \geq 2$. Together with the results on the shuffle brackets in $\cite{Ba2}$ we obtain the following 
\begin{thmx}(\cite{BT},\cite{Ba2}) \label{thm:4}
For all $s_1,\ldots,s_l\ge1$ there exist shuffle regularized multiple Eisenstein series $G^{\sh}_{s_1,\ldots,s_l} \in \C[\![q]\!]$ with the following properties:
\begin{enumerate}[i)]
\item They are holomorphic functions on the upper half-plane (by setting $q=\exp(2\pi i \tau)$) having a Fourier expansion with the shuffle regularized multiple zeta values as the constant term. 
\item They fulfill the shuffle product.
\item They can be written as a linear combination of multiple zeta values, powers of $(-2\pi i)$ and shuffle brackets $[\dots]^\sh \in \bMD$. 
\item For integers $s_1,\ldots,s_l\ge2$ they equal the multiple Eisenstein series
\[ G^{\sh}_{s_1,\ldots,s_l}(\tau)=G_{s_1,\ldots,s_l}(\tau) \]
and therefore they fulfill the stuffle product in these cases. 
\end{enumerate}
\end{thmx} 

We now study the $\Q$-algebra spanned by the $G^\sh$ and its relation to multiple zeta values. Theorem \ref{thm:4} iv) gives a subset of the double shuffle relations between the $G^\sh$, since the stuffle product is just fulfilled for the case $s_1,\dots,s_l \geq 2$. A natural question is, if they also fulfill the stuffle product when some indices $s_j$ are equal to $1$. For some cases this was proven in \cite{Ba2}. For example it was shown, that
\begin{equation} \label{eq:stufflew5}
 G^\sh_2 \cdot G^\sh_{2,1} = G^\sh_{2,1,2} + 2 G^\sh_{2,2,1} + G^\sh_{2,3} +G^\sh_{4,1}\,.
\end{equation}
The method to prove this was to introduce stuffle regularized multiple Eisenstein series $G^\ast_{s_1,\dots,s_l}$, which fulfill by construction the stuffle product and which equal the classical multiple Eisenstein series in the $s_1,\dots,s_l \geq 2$ cases. Since both $G^\ast$ and $G^\sh$ can be written in terms of multiple zeta values and bi-brackets it was possible to compare these two regularization. It was shown that all $G^\sh$ appearing in \eqref{eq:stufflew5} equal the $G^\ast$ ones, from which this equation followed. In contrast to the shuffle regularized multiple Eisenstein series the stuffle regularized ones could not be defined for all $s_1,\dots,s_l \geq 1$, but we have the following results: 

\begin{thmx}(\cite{Ba2}) For all $s_1,\dots,s_l \geq 1$ and $M \geq 1$ there exists $G_{s_1,\dots,s_l}^{\ast,M} \in \C[\![q]\!]$ with the following properties
\begin{enumerate}[i)]
\item They are holomorphic functions on the upper half-plane (by setting $q=\exp(2\pi i \tau)$) having a Fourier expansion with the stuffle regularized multiple zeta values as the constant term. 
\item They fulfill the stuffle product. 
\item In the case where the limit $G^\ast_{s_1,\dots,s_l} :=  \lim_{M\to\infty} G_{s_1,\dots,s_l}^{\ast,M}$ exists, the functions $G^\ast_{s_1,\dots,s_l}$ are a linear combination of multiple zeta values, powers of $(-2\pi i)$ and bi-brackets.
\item For $s_1,\dots,s_l \geq 2$ the $G^\ast_{s_1,\dots,s_l}$ exist and equal the classical multiple Eisenstein series
\[ G_{s_1,\dots,s_l}(\tau) = G^\ast_{s_1,\dots,s_l}(\tau) \,.\] 
\end{enumerate}
\end{thmx} 
It is still an open question which extended double shuffle relations of multiple zeta values are also fulfilled for the $G^\sh$. Or equivalently, under what circumstances the product of two $G^\sh$ can be expressed using the stuffle product formula. Clearly there are some double shuffle relations which can't be fulfilled by multiple Eisenstein series. For example not all of the Euler relations \eqref{eq:eulerrelation} are fulfilled since $G_2^2$ is not a multiple of $G_4$ as $G_2$ is not modular and $G_6^2$ is not a multiple of $G_{12}$ as there are cusp forms in weight $12$. In Section \ref{sec:stufflereg} we will explain this failure in terms of the double shuffle relations which are fulfilled by multiple Eisenstein series.\newline

After the discussion above, we believe that Question \ref{qu1} got a satisfying answer given by the regularized multiple Eisenstein series $G^\sh$ and $G^\ast$. To go back from multiple Eisenstein series to multiple zeta values one can consider the projection to the constant term. But there is another direct connection of brackets, and therefore also of the subalgebra of modular forms, to multiple zeta values. The brackets can be seen as a $q$-analogue of multiple zeta values. A $q$-analogue of multiple zeta values is said to be a $q$-series which gives back multiple zeta values in the case $q \rightarrow 1$. Define for $k\in \N$ the map $Z_k: \Q[\![q]\!] \rightarrow \R \cup \{ \infty \}$ by 
\[ Z_k(f) = \lim_{q \to 1} (1-q)^{k}f(q) \,.\]
\begin{propx}(\cite[Prop. 6.4]{BK})\label{prop:zkintro}
For $s_1\geq 2$ and $s_2,\dots,s_l\geq 1$ the map $Z_k$ sends a bracket to the corresponding multiple zeta value, i.e. 
\[ Z_k\left( [s_1, \dots , s_l ] \right) = \left\{
\begin{array}{cl} \zeta(s_1,\dots,s_l)\,, &  s_1+\dots+s_l = k ,  \\ 0\,, &s_1+\dots+s_l < k
\,. \end{array}   \right. \]
\end{propx}
Since every relation of multiple zeta values in a given weight $k$ is, by Proposition \ref{prop:zkintro}, in the kernel of the map $Z_k$, this kernel was studied in \cite{BK} with the following result
\begin{thmx}(\cite[Thm. 1.13]{BK})\label{thm:qana}
\begin{enumerate}[i)]
\item For any $f\in \MD$ which can be written as a linear combination of brackets with weight $\leq k-2$ we have $\dif f \in \ker Z_k$.
\item Any cusp form for $\Sl_2(\Z)$ of weight $k$ is in the kernel of $Z_k$.
\end{enumerate}
\end{thmx} 
We give an example for Theorem \ref{thm:qana}  ii): Using the theory of brackets (Corollary \ref{cor:delta}) we can prove for the cusp form $\Delta = q\prod_{n>0} \left(1-q^n \right)^{24} \in S_{12}(\Sl_2(\Z))$  the representation 
\begin{align}\label{eq:deltal21}
 -\frac{1}{2^6\cdot 5 \cdot 691}  \Delta  &=  168 [5,7]+150 [7,5]+28 [9,3] \notag \\
&+\frac{1}{1408} [2] - \frac{83}{14400}[4] +\frac{187}{6048} [6] - \frac{7}{120} [8] - \frac{5197}{691} [12] \,.
\end{align}
Letting $Z_{12}$ act on both sides of \eqref{eq:deltal21} one obtains a new proof for the relation \eqref{eq:exotic}, i.e., 
\[ \frac{5197}{691} \zeta(12) =  168 \zeta(5,7)+150 \zeta(7,5) + 28 \zeta(9,3) \,. \] 
Another reason for studying the enlargement of the brackets given by the bi-brackets is the following: In weight $4$ one has the following relation of multiple zeta values $\zeta(4) = \zeta(2,1,1)$, i.e. it is $[4] - [2,1,1] \in \ker Z_4$. But this element can't be written as a linear combination of cusp forms, lower weight brackets or derivatives. But one can show, by using the double shuffle relations of bi-brackets, that 
\begin{equation} \label{eq:bibracketker}
 [4] - [2,1,1] = \frac{1}{2}\left( \dif[1] + \dif[2] \right) - \frac{1}{3} [2] - [3] + \mb{2,1}{1,0}
\end{equation}
 and $ \mb{2,1}{1,0} \in \ker Z_4$. To describe the kernel of the map $Z_k$ was in fact our first motivation to study the bi-brackets. Equation \eqref{eq:bibracketker} is also an example for the above mentioned Conjecture \ref{conj:mdbd}, since it shows that the bi-bracket $\mb{2,1}{1,0}$ can be written as brackets and therefore is an element in $\MD$.  \newline  \newline

\section*{Outlook and related work} 
In the following paragraphs a.)-g.) we want to mention some related works and give an outlook to open questions. \newline 

{\bf a.)} There are still a lot of open questions concerning multiple Eisenstein series as well as the space of (bi-)brackets. 
After the above mentioned works \cite{BK},\cite{Ba2} and \cite{BT} we now have a good definition of regularized multiple Eisenstein series given by the $G^\sh$. For the structure of the space spanned by these series there are still several open questions. 
\begin{enumerate}[i)]
\item What exactly is the failure of the stuffle product for the $G^{\shuffle}$ and when does it hold?
\item For which indices $s_1,\dots,s_l \in \N$ do we have $G^{\shuffle}_{s_1,\ldots,s_r}(\tau) = G^{\ast}_{s_1,\ldots,s_r}(\tau)$? Is there an explicit connection between these two regularizations similar to the regularized multiple zeta values given by the map $\rho$ in \cite{IKZ}?
\item What is the dimension of the space of (shuffle) regularised multiple Eisenstein series? Is there an explicit basis similar to the Hoffman basis of multiple zeta values (Which is given by all multiple zeta values $\zeta(s_1,\dots,s_l)$ with $s_j \in \{2,3\}$)?
\item Which linear combinations of multiple Eisenstein series are modular forms for $\Sl_2(\Z)$? Is there an explicit way to describe the modular defect? 
\item Is the space of multiple Eisenstein series closed under the derivative $\dif = q \frac{d}{dq}$ ? Meanwhile this question was also already addressed in \cite{Ba5}.
\item What is the kernel of the projection to the constant term? Does it consist of more than derivatives and cusp forms?
\item Is there a general theory behind the connection of the Fourier expansion of multiple Eisenstein series and the Goncharov coproduct? Can we equip the space of multiple Eisenstein series with a coproduct structure in an useful way?
\end{enumerate}
Especially the last questions seems to be interesting since the connection to the coproduct of formal iterated integrals is quite mysterious and it seems that there might be a geometric interpretation for this connection. \newline

{\bf b.)} Several $q$-analogues of multiple zeta values were studied in recent years. The first works on this area are  \cite{db}, \cite{Zh0}, \cite{Sch} and \cite{YOZ}. Possible double shuffle structures are discussed for example in \cite{yt}, \cite{JMK}, \cite{S} and \cite{Zh}, where the last one gives also a nice overview of various different $q$-analogue models. Often these $q$-analogues have a product structure similar to the stuffle product of multiple zeta values. To obtain something which corresponds to the shuffle product one usually needs to modify the space and add extra elements (like derivatives) or consider index sets $(s_1,\dots,s_r)$ with $s_j \in \Z$ or $s_j \geq 0$. The picture is similar for bi-brackets, where we consider double indices $ \mb{s_1, \dots , s_l}{r_1,\dots,r_l}$ to obtain an analogue for both products in a very natural way. This gives a lot of linear relations similar to the double shuffle relations. Numerical experiments suggest, that every bi-bracket can be written as a linear combination of brackets and therefore (conjecturally) every relation of bi-brackets gives rise to relations between multiple zeta values by applying the map $Z_k$. \newline  

{\bf c.)} In the case of multiple zeta values one way to give upper bounds for the dimension is to study the double shuffle space (\cite{IKZ}, \cite{io}). Similarly, one can study the partition shuffle space 
\[
\mathbb{PS}(k-l,l) = \big\{ f \in \Q[X_1,..,X_l,Y_1,..,Y_l] \,\big| \, 
\deg f = k-l, \,\, f \big|_P -f = f \big|_{\operatorname{Sh}_j} =0 \,\, \forall j \big\},
\]
for bi-brackets, where $|_P$ is the involution given by the partition relation (see Section \ref{sec:bibrackgen}, \eqref{eq:partition}) and  $|_{\operatorname{Sh}_j}$ is given by the sum of all shuffles of type $j$ similar to the one in \cite{io}. Counting the number of these polynomials it is possible to give upper bounds for the dimensions of the space of bi-brackets. This approach therefore enabled us to prove the conjecture $\MD=\bMD$ up to weight $7$ in a current work in progress (\cite{BK3}). Therefore, considering the space $\mathbb{PS}(k-l,l)$ in more detail might be crucial to understand the structure of bi-brackets. \newline

{\bf d.)} In this work we were interested in modular forms for the full modular group and therefore studied the level $1$ case. In \cite{KT} the authors studied double Eisenstein series and double zeta values of level $2$. They also derive the Fourier expansion of these series which involves similar calculation as in the level $1$ case. One result is, that they derive the dimension of the space of double Eisenstein series and give also an upper bound for the dimension of double zeta values of level 2, which involves the dimension of the spaces of cusp forms of level 2. Beside the work on Level 2 double Eisenstein series there are also work for level $N$ double Eisenstein series of H. Yuan and J. Zhao in \cite{YZ}. Later these authors also considered a level $N$ version of the brackets in $\cite{YZ2}$.  \newline

{\bf e.)} At the end of \cite{KT} the authors give a proof for an upper bound of the dimension of double zeta values in even weight. We want to recall this result, since the presented results in the present work might be able to use these ideas for higher lengths. Consider the space spanned by all normalized double Eisenstein series $(-2\pi i)^{-r-s} G_{r,s}(\tau)$ in even weight $k=r+s$. Denote by $\pi_i$ the projection of this space to the imaginary part. Using the Fourier expansion of double Eisenstein series the authors can write down the matrix representation of $\pi_i$ explicitly. Together with well known results on period polynomials they obtain 
\[ \dim_\Q \langle \zeta(r,k-r) \mid 2 \leq r \leq k-1 \rangle_\Q \leq \frac{k}{2} - 1 - \dim S_k \,.\] 
Due to the Broadhurst-Kreimer conjecture \ref{conj:bk} it is  conjectured that this is actually an equality. The key fact here is, that it is possible to write down an explicit basis of the imaginary part and the matrix representation of $\pi_i$. 
To also obtain upper bounds for the dimensions of multiple zeta values in higher lengths, one might try to use the exact same method as in the length two case. 
The imaginary part of the (again normalized with the factor $(-2\pi i)^{-k}$) multiple Eisenstein series is more complicated since it involves the functions $g$ in different length, where it is known that they are not linearly independent anymore. But the algebraic structure of the $g$ or more precisely of the brackets $[..]$ are subject of the current work. It is quite possible that the results on the brackets enable one to study the projection of the imaginary part of multiple Eisenstein series to obtain upper bounds for the Broadhurst-Kreimer conjecture. \newline

{\bf f.)} The multiple Eisenstein series and the bi-brackets itself also have connections to counting problems in enumerative geometry:
\begin{enumerate}[i)]
\item In \cite{ar} and \cite{R} the author studies $q$-series $A_k(a) \in \Q[\![q]\!]$ which arises in counting certain types of hyperelliptic curves. One of the results is, that the $A_k(q)$ are contained in the ring of quasi-modular forms. The connection to the brackets is given by the fact that $A_k(q) = [\underbrace{2,\dots,2}_k]$. The results of \cite{ar} can also be obtained by using an explicit calculation of the Fourier expansion of  $G_{2,\dots,2}$ which will be done in an upcoming work \cite{Ba3}. 
\item In \cite{O} and \cite{QY} the authors connect certain $q$-analogues of multiple zeta values to Hilbert schemes of points on surfaces. These $q$-analogues are just particular linear combinations of brackets as explained in \cite{BK2} and Section \ref{sec:otherqana}. 
\item The coefficients of bi-brackets also occur naturally when counting flat surfaces \cite{Zo}, i.e. certain covers of the torus. 
\end{enumerate}

{\bf g.)} There also exists different "multiple"-versions of classical Eisenstein series. One of them is treated in \cite{BTs}, where the authors discuss the series defined by
\begin{align*}
\mathfrak{G}_{2p_1,\ldots,2p_r}(\tau)& =\sum_{m\in \mathbb{Z}} \sum_{n_1 \in \mathbb{Z}\atop (m,n_1)\not=(0,0)}\cdots \sum_{n_r\in \mathbb{Z}\atop (m,n_r)\not=(0,0)}\prod_{j=1}^{r}\frac{1}{(m+n_j\tau)^{2p_j}}
\end{align*}
for $r\in \mathbb{N}_{\geq 2}$ and $p_1,\ldots,p_r\in \mathbb{N}$ and prove (Theorem 2) that for $r\in \mathbb{N}_{\geq 2}$ and $p_1,\ldots,p_r\in \mathbb{N}$, 
\begin{align*}
& \tau^{2(p_1+\cdots+p_r)}\mathfrak{G}_{2p_1,\ldots,2p_r}(\tau)\in \mathbb{Q}\left[ \tau^2,\, \pi^2,\,G_{2}(\tau),\,G_4(\tau),\,G_6(\tau)\right].
\end{align*}
The methods used to prove these statements are similar to the methods used in the calculation of the Fourier expansion of multiple Eisenstein series. But besides this there does not seem to be a direct connection to the multiple Eisenstein series presented here. 
 
\section*{Acknowledgment} This paper has served as the introductory part of my cumulative thesis written at the University of Hamburg. First of all I would like to thank my supervisor Ulf K\"uhn for his continuous, encouraging and patient support during the last years. Besides this I also want to thank several people for supporting me during my PhD project by whether giving me suggestion and ideas, letting me give  talks on conferences and seminars, proof reading papers or having general discussions on this topic with me. A big "thank you" goes therefore to Olivier Bouillot, Kathrin Bringmann, David Broadhurst, Kurusch Ebrahimi-Fard, Herbert Gangl, Jos\'{e} I. Burgos Gil, Masanobu Kaneko, Dominique Manchon, Nils Matthes, Martin M\"oller, Koji Tasaka, Don Zagier, Jianqiang Zhao  and Wadim Zudilin. Finally I would like to thank the referee for various helpful comments and remarks.

\section{Multiple Eisenstein series}
In this section we are going to introduce multiple zeta values and present the multiple Eisenstein series and their Fourier expansion. Especially the construction of the Fourier expansion of multiple Eisenstein series in Section \ref{sec:mes} was rewritten for this survey. It will be a shortened version of the construction given in \cite{Ba} using results by Bouillot obtained in \cite{Bo}. This section is not part of the works $\cite{BK},\cite{BK2},\cite{BT}$ and $\cite{BK2}$.
Before we discuss multiple Eisenstein series, we give a short review of multiple zeta values and their algebraic structure given by the stuffle and shuffle product. In order to describe these two products we will use quasi-shuffle algebras, introduced by Hofmann in \cite{H}, which will also be needed later when we deal with the generating series of multiple divisor-sums (brackets) and their generalizations given by the bi-brackets. 

\subsection{Multiple zeta values and quasi-shuffle algebras}
Multiple zeta values are natural generalizations of the Riemann zeta values that are defined\footnote{Some authors use the opposite convention $0 < n_1 <\dots < n_l$ in the definition of multiple zeta values. This is in particular the case for the work \cite{BT}, where this opposite convention is used for multiple zeta values and multiple Eisenstein series.} for integers $s_1 > 1$ and $s_i \geq 1$ for $i>1$ by
\[ \zeta( s_1 , \dots , s_l ) :=  \sum_{n_1 >n_2 >\dots>n_l>0 } \frac{1}{n_1^{s_1} \dots n_l^{s_l} } \,. \]
We denote the  $\Q$-vector space of all multiple zeta values of weight $k$ by
\[ \MZ_k :=\big <\,\zeta(s_1,\dots,s_l) \, \big| \, s_1 + \dots + s_l = k 
\textrm{ and } l>0 \big>_\Q.  \]
It is well known that  the product of two multiple zeta values can be written as a linear combination of multiple zeta values of the same weight by using the stuffle or shuffle relations (See for example \cite{IKZ}, \cite{Zu2}). Thus they generate a $\Q$-algebra $\MZ$. There are several connections of these numbers to modular forms for the full modular group.
In the smallest length the stuffle product reads
\begin{align*}
\zeta(s_1) \cdot \zeta(s_2) &= \sum_{n_1>0} \frac{1}{n_1^{s_1}} \sum_{n_2>0} \frac{1}{n_2^{s_2}} \\
&= \sum_{n_1>n_2>0} \frac{1}{n_1^{s_1}n_2^{s_2}} + \sum_{n_2>n_1>0} \frac{1}{n_1^{s_1}n_2^{s_2}} + \sum_{n_1=n_2>0} \frac{1}{n_1^{s_1+s_2}} \\
&=\zeta(s_1,s_2) + \zeta(s_2,s_1) + \zeta(s_1+s_2) \,. 
\end{align*}
For length $1$ times length $2$ the same argument gives
\begin{align*}
\zeta(s_1) \cdot \zeta(s_2,s_3) &= \zeta(s_1,s_2,s_3) +\zeta(s_2,s_1,s_3)+\zeta(s_2,s_3,s_1)\\
&+\zeta(s_1+s_2,s_3)+\zeta(s_2,s_1+s_3) \,. 
\end{align*}
The second expression for the product, the shuffle product, comes from the iterated integral expression of multiple zeta values. For example it is
\[ \zeta(2,3) = \int_{{\small 1 > t_1 > \dots > t_5 > 0}} \underbrace{\frac{dt_1}{t_1} \cdot \frac{dt_2}{1-t_2}}_{2} \cdot \underbrace{ \frac{dt_3}{t_3}\cdot  \frac{dt_4}{t_4} \cdot \frac{dt_5}{1-t_5}}_{3} \,. \]
Multiplying two of these integrals one obtains again a linear combination of multiple zeta values as for example
\[ \zeta(2) \cdot \zeta(3) =  \zeta(2,3) + 3 \zeta(3,2) + 6 \zeta(4,1) \,.\]
More generally the smallest length case is given by 
\begin{equation} \label{eq:shufflelen1}
 \zeta(s_1)\cdot \zeta(s_2) = \sum_{\substack{a+b=s_1+s_2 \\ a > 1}} \left( \binom{a-1}{s_1-1}+\binom{a-1}{s_2-1} \right) \zeta(a,b) \,.
\end{equation} 
To describe these two product structures precisely we will use the language of quasi-shuffle algebras as introduced in \cite{H} and \cite{hi}. 

\begin{dfn} \label{def:quasishuffle}
Let $A$ (the alphabet) be a countable set of letters, $\Q A$ the $\Q$-vector space generated by these letters and $\Q\langle A \rangle$ the noncommutative polynomial algebra over $\Q$ generated by words with letters in $A$. For a commutative and associative product $\diamond$ on $\Q A$, $a,b \in A$ and $w,v \in \Q\langle A\rangle$ we define on $\Q\langle A\rangle$ recursively a product by $1 \odot w=w \odot 1=w$ and
\begin{equation} \label{qshp} aw \odot bv := a(w \odot bv) + b(aw \odot v) + (a \diamond b)(w \odot v) \,. \end{equation}
 By a result of Hoffman (\cite[Thm. 2.1]{hi}) $(\Q\langle A \rangle,\odot)$ is a commutative $\Q$-algebra which is called a \emph{quasi-shuffle algebra}.
\end{dfn}

To describe the stuffle and the shuffle product for multiple zeta values we need to deal with two different alphabets $A_{xy}$ and $A_z$. The first alphabet is given by $A_{xy} := \{ x,y \}$  and we set $\h = \Q\langle A_{xy} \rangle$ and $\h^1 = 1 \cdot \Q +  \h y$, with $1$ being the empty word. It is easy to see that $\h^1$ is generated by the elements $z_j = x^{j-1} y$ with $j \in \N$, i.e. $\h^1 = \Q\langle A_{z} \rangle$ with the second alphabet $A_z := \{ z_1,z_2, \dots \}$. Additionally, we define $\h^0 = 1 \Q + x \h y$. 

\begin{enumerate}[i)]
\item On $\h^1$ we have the following quasi-shuffle product with respect to the alphabet $A_z$, called the \emph{stuffle product}. We denote it by $\ast$ and define it as the quasi-shuffle product with $z_j \diamond z_i = z_{j+i}$. For $a,b \in \N$ and  $w,v \in \h^1$ we therefore have:
\begin{align*}
\begin{split}
 z_a w \ast z_b v &= z_a(w \ast z_b v) + z_b(z_a w \ast v)  + z_{a+b}( w \ast v) \,.
\end{split}
\end{align*}
By $(\h^1, \ast)$ we denote the corresponding $\Q$-algebra. 
\item On the alphabet $A_{xy}$ we define the \emph{shuffle product} as the quasi-shuffle product with $\diamond \equiv 0$, and by $(\h^1,\sh)$ we denote the corresponding $\Q$-algebra. 
\end{enumerate}
It is easy to check that $\h^0$ is closed under both products $\ast$ and $\sh$ and therefore we have also the two algebras $(\h^0,\ast)$ and $(\h^0,\sh)$.

By the definition of multiple zeta values as an ordered sum and by the iterated integral expression one obtains algebra homomorphisms $Z: (\h^0, \ast) \rightarrow \MZ$ and $Z: (\h^0, \shuffle) \rightarrow \MZ$ by sending $w=z_{s_1} \dots z_{s_l} $ to $\zeta(w) = \zeta(s_1,\dots,s_l)$, since the words in $\h^0$ correspond exactly to the indices for which the multiple zeta values are defined. It is a well known fact, that these algebra homomorphisms can be extended to $\h^1$:
\begin{prop}(\cite[Prop. 1]{IKZ})\label{prop:mzvreg}
There exist algebra homomorphisms
\begin{align*}
Z^\ast : (\h^1, \ast) \longrightarrow \MZ  \quad \text{ and } \quad 
Z^\shuffle : (\h^1, \shuffle) \longrightarrow \MZ \,, 
\end{align*}
which are uniquely determined by $Z^\ast(w) = Z^\shuffle(w) = \zeta(w)$ for $w\in \h^0$ and by their images on the word $z_1$, which we set $0$ here, i.e. \mbox{$Z^\ast(z_1) = Z^\shuffle(z_1) = 0$}. 
\end{prop} \qed

\subsection{Multiple Eisenstein series and the calculation of their Fourier expansion}\label{sec:mes}
The Riemann zeta values appear as the constant term in the Fourier expansion of classical Eisenstein series. These series are defined for $\tau \in \Ha$ by 
\begin{equation}\label{eq:defgklattice}
G_k(\tau) = \frac{1}{2} \sum_{\substack{(m,n) \in \Z^2 \\ (m,n) \neq (0,0)}} \frac{1}{(m\tau + n)^k } .
\end{equation}
where $k>2$ is the called the weight. 
 Splitting the summation into the parts $m=0$ and $m \in \Z \backslash {0}$ we obtain for even $k$
\[ G_k(\tau) = \frac{1}{2} \sum_{n \neq 0} \frac{1}{n^k} + \sum_{m = 1}^\infty \left( \sum_{n \in \Z} \frac{1}{(m\tau + n)^k} \right) \,. \]
To calculate the Fourier expansion of the sum on the right, one uses the  well known Lipschitz summation formula ($q=e^{2\pi i \tau}$)
\begin{equation}\label{lipschitz-sum}
 \sum_{d \in \Z} \frac{1}{(\tau + d)^k} = \frac{(-2 \pi i)^k}{(k-1)!} \sum_{m=1}^\infty m^{k-1} q^m \,,
\end{equation}
which is valid for $k>1$. With \eqref{lipschitz-sum} we obtain
\begin{align}\label{eq:gk}
 G_k(\tau) &= \zeta(k) + \frac{(-2\pi i)^k}{(k-1)!} \sum_{m = 1}^\infty \sum_{n=1}^\infty n^{k-1} q^{mn} = \zeta(k) + \frac{(-2\pi i)^k}{(k-1)!} \sum_{n = 1}^\infty \sigma_{k-1}(n) q^{n}\,,
\end{align}
where $\sigma_k(n) = \sum_{d | n} d^k$ denote the divisor-sum. Formula \eqref{eq:gk} also makes sense for odd $k$ but does not give a modular form, since there are no non trivial modular forms of odd weight. The sum in \eqref{eq:defgklattice} vanishes for odd $k$, therefore instead of summing over the whole lattice, we restrict the summation to the positive lattice points, with positivity coming from an order on the lattice $\Z \tau + \Z$. This in turn will also enable us to give an multiple version of the Eisenstein series in an obvious way. 

Let $\Lambda_\tau = \Z \tau + \Z$ be a lattice  with $\tau \in \Ha := \left\{ x+iy \in \C \mid  y>0 \right\}$. An order $\succ$ on $\Lambda_\tau$ is defined by setting (see \cite{GKZ})
\[ \lambda_1 \succ \lambda_2 :\Leftrightarrow \lambda_1 - \lambda_2 \in P\]
for $\lambda_1,\lambda_2 \in \Lambda_\tau$ and the following set $P$, which we call the set of positive lattice points
\[ P := \left\{ l\tau + m \in \Lambda_\tau \mid {\color{blue}{l > 0}}  \vee  \left( {\color{red}{ l=0 \wedge m>0 }}\right) \right\}  = \UU \cup \RR \]
\begin{figure}[H]
\begin{center}
\input{dia_p0}\\
The set $P$ for the case $\tau =i$.
\end{center}
\end{figure}

\begin{dfn}
For $s_1 \geq 3, s_2,\dots,s_l \geq 2$ the \emph{multiple Eisenstein series} is defined by
\[ G_{s_1,\dots,s_l}(\tau) := \sum_{\substack{\lambda_1 \succ \dots \succ \lambda_l \succ 0\\ \lambda_i \in \Lambda_\tau}} \frac{1}{\lambda_1^{s_1} \dots \lambda_l^{s_l}}   \,.\]
With $k=s_1 + \dots + s_l$ we denote the weight and with $l$ its length. 
\end{dfn}It is easy to see that these are holomorphic functions in the upper half-plane and that they fulfill the stuffle product, i.e. for example
\[ G_3(\tau) \cdot G_4(\tau) = G_{4,3}(\tau) + G_{3,4}(\tau) + G_{7}(\tau)\,. \] 
By definition it is $G_{s_1,\dots,s_l}(\tau + 1) = G_{s_1,\dots,s_l}(\tau)$, i.e. there exists a Fourier expansion of $G_{s_1,\dots,s_l}$ in $q=e^{2\pi i \tau}$. To write down the Fourier expansion of multiple Eisenstein series we need to introduce the following $q$-series which will be studied in detail in Section \ref{sec:brackets}. For $s_1,\dots,s_l \geq 1$ we define 
\[ [s_1,\dots,s_l] := \sum_{\substack{u_1 > \dots > u_l> 0 \\ v_1, \dots , v_l >0}} \frac{v_1^{s_1-1} \dots v_l^{s_l-1}}{(s_1-1)!\dots(s_l-1)!}   \cdot q^{u_1 v_1 + \dots + u_l v_l}  \in \Q[\![q]\!] \,. \]  
and write $g_{s_1,\dots,s_l}(\tau) := (-2\pi i)^{s_1+\dots+s_l} [s_1,\dots,s_l]$, which is an holomorphic function in the upper half-plane by setting $q=e^{2\pi i \tau}$.

\begin{thm}(\cite{Ba}, Fourier expansion) \label{thm:Fourier}
For $s_1\geq 3, s_2,\dots,s_l \geq 2$ the $G_{s_1,\dots,s_l}(\tau)$ can be written as a $\MZ$-linear combination of the functions $g$. More precisely there are rational numbers $\lambda_{r,j} \in \Q$, for $r=(r_1,\dots,r_l)$ and $1\leq j \leq l-1$, such that (with $k=s_1+\dots+s_l$)
\begin{align*} G_{s_1,\dots,s_l}(\tau) \,= \,\zeta(s_1,\dots,s_l) \;+\; \sum_{\substack{1\leq j \leq l-1 \\ r_1 +\dots+r_l = k}} \lambda_{r,j}
\cdot  \zeta(r_1,\dots,r_j) \cdot g_{r_{j+1},\dots,r_l}(\tau) \;+\; g_{s_1,\dots,s_l}(\tau) \,. 
\end{align*}
\end{thm}
Even though the proof of this statement is the main result of \cite{Ba} we will give a shortened version of it in the following.

The condition $s_1\geq 3$ is necessary for the absolute convergence of the sum. Nevertheless we can also allow the case $s_1 = 2$ by using the Eisenstein summation as it was done in \cite{BT} Definition 2.1. This corresponds to the usual way of defining the quasi-modular form $G_2$ in length one. Since the construction of the Fourier expansion described below uses exactly this Eisenstein summation the Theorem \ref{thm:Fourier} is also valid for $s_1 \geq 2$. 

For example the triple Eisenstein series $G_{3,2,2}$ can be written as
\begin{align*}
G_{3,2,2}(\tau) =& \zeta(3,2,2) + \left( \frac{54}{5} \zeta(2,3) + \frac{51}{5} \zeta(3,2) \right) g_2(\tau) + \frac{16}{3} \zeta  (2,2)g_3(\tau) \\
&+3 \zeta(3) g_{2,2}(\tau) + 4 \zeta(2) g_{3,2}(\tau) + g_{3,2,2}(\tau)\,.
\end{align*}

To derive the Fourier expansion we introduce the following functions, that can be seen as a multiple version of the term $\sum_{n \in \Z} \frac{1}{(x+n)^k}$ appearing in the calculation of the Fourier expansion of classical Eisenstein series. 

\begin{dfn}\label{def:multitangent}
For $s_1,\dots,s_l \geq 2$ we define the multitangent function of length $l$ by
\[  \Psi_{s_1,\dots,s_l}(x) = \sum_{\substack{ n_1 > \dots > n_l \\ n_i \in \Z}} \frac{1}{(x+n_1)^{s_1} \dots (x+n_l)^{s_l}} \,.  \]
In the case $l=1$ we also refer to these as monotangent function. 
\end{dfn}
These functions were introduced and studied in detail in \cite{Bo}. One of the main results there, which is crucial for the calculation of the Fourier expansion presented here, is the following theorem which reduces the multitangent functions into monotangent functions.

\begin{thm}(\cite[Thm. 3]{Bo}, Reduction of multitangent into monotangent functions)\label{thm:reduction}
For $s_1,\dots,s_l \geq 2$ and $k=s_1+\dots+s_l$ the multitangent function can be written as a $\MZ$-linear combination of monotangent functions, more precisely there are $c_{k,h} \in \MZ_{k-h}$ such that
\[ \Psi_{s_1,\dots,s_l}(x) = \sum_{h=2}^{k} c_{k-h} \Psi_h (x) \,. \]
\end{thm}
\begin{prf}
An explicit formula for the coefficients $c_k$ is given in Theorem 3 in \cite{Bo}. The proof uses partial fraction and a non trivial relation between multiple zeta values to argue that the sum starts at $h=2$. 
For example in length two it is
\begin{align}\label{eq:psi32}
\begin{split}
&\Psi_{3,2}(x) = \sum_{m_1 > m_2} \frac{1}{(x+m_1)^3 (x+m_2)^2} \\
&= \sum_{m_1 > m_2}  \left( \frac{1}{(m_1-m_2)^2 (x+m_1)^3} +\frac{2}{(m_1-m_2)^3 (x+m_1)^2} + \frac{3}{(m_1-m_2)^4 (x+m_1)} \right) \\
&+\sum_{m_1 > m_2} \left( \frac{1}{(m_1-m_2)^3 (x+m_2)^2}  -  \frac{3}{(m_1-m_2)^4 (x+m_2)} \right) \\
&=  3 \zeta(3) \Psi_2(x) +  \zeta(2) \Psi_3(x)  \,. 
\end{split}
\end{align}
\end{prf}

The connection between the functions $g$ and the monotangent functions is given by the following 
\begin{prop}
For $s_1,\dots,s_r \geq 2$ the functions $g$ can be written as an ordered sum of monotangent functions
\begin{align*}
&g_{s_1,\ldots,s_l}(\tau) = \sum_{m_1 > \dots > m_l > 0} \Psi_{s_1}(m_1\tau) \dots \Psi_{s_r}(m_r\tau) \,.
\end{align*}
\end{prop}
\begin{prf}
This follows directly from the Lipschitz formula \eqref{lipschitz-sum} and the definition of the functions $g$. 
\end{prf}

{\bf Preparation for the Proof of Theorem \ref{thm:Fourier}:} We will now recall the construction of the Fourier expansion of multiple Eisenstein series introduced in \cite{Ba}, in order to prove Theorem \ref{thm:Fourier}. 
To calculate the Fourier expansion we rewrite the multiple Eisenstein series as 
\begin{align*}
 G_{s_1,\dots,s_l}(\tau) &= \sum_{\lambda_1 \succ \dots \succ \lambda_l \succ 0} \frac{1}{\lambda_1^{s_1} \dots \lambda_l^{s_l}} \\
 &= \sum_{(\lambda_1,\dots,\lambda_l) \in P^l} \frac{1}{(\lambda_1+\dots+\lambda_l)^{s_1} (\lambda_2 + \dots + \lambda_l )^{s_2} \dots \lambda_l^{s_l} } \,.
\end{align*}
We decompose the set of tuples of positive lattice points $P^l$  into the $2^l$ distinct subsets $A_1 \times \dots \times A_l \subset P^l$ with $A_i \in \{R,U\}$ and write 
{\small
\[ G^{A_1 \dots A_l}_{s_1,\dots,s_l}(\tau) := \sum_{(\lambda_1,\dots,\lambda_l) \in A_1 \times \dots \times A_l} \frac{1}{(\lambda_1+\dots+\lambda_l)^{s_1} (\lambda_2 + \dots + \lambda_l )^{s_2} \dots \lambda_l^{s_l} } \]
}
this gives the decomposition
\[ G_{s_1,\dots,s_l} = \sum_{A_1,\dots,A_l \in \{R,U\}} G^{A_1\dots A_l}_{s_1,\dots,s_l} \,. \]
In the following we identify the $A_1 \dots A_l$ with words in the alphabet $\{R,U\}$. In length $l=1$ we have $G_{k}(\tau) =G_k^R(\tau) + G_k^U(\tau)$ and
\begin{align*}
G_k^R(\tau) &= \sum_{\substack{m_1 = 0 \\ n_1 > 0}} \frac{1}{(0\tau + n_1)^k} = \zeta(k) \,,\\
G_k^U(\tau) &= 	\sum_{\substack{m_1 > 0 \\ n_1 \in \Z}} \frac{1}{(m_1\tau + n_1)^k} = \sum_{m_1>0} \Psi_k(m_1\tau) \,, 
\end{align*}
where $\Psi_k$ is the  monotangent function given by
\[ \Psi_k(x) = \sum_{n \in \Z} \frac{1}{(x+n)^k} \,.\]
To calculate the Fourier expansion of $G_k^U$ one uses the Lipschitz formula \eqref{lipschitz-sum}. In general the $G^{U^l}_{s_1,\dots,s_l}$ can be written as
\begin{align*}
G^{U^l}_{s_1,\dots,s_l}(\tau) &= \sum_{\substack{m_1>\dots>m_l>0 \\ n_1,\dots,n_l \in \Z}} \frac{1}{(m_1\tau+n_1)^{s_1} \dots (m_l\tau + n_l)^{s_l}}\\
&= \sum_{m_1 > \dots > m_l>0} \Psi_{s_1}(m_1\tau) \dots \Psi_{s_l}(m_l \tau) \\
&=  \frac{(-2\pi i)^{s_1+\dots+s_l}}{(s_1-1)!\dots(s_l-1)!} \sum_{\substack{m_1 > \dots > m_l>0\\d_1,\dots,d_l>0}} d_1^{s_1-1} \dots d_l^{s_l-1} q^{m_1 d_1 + \dots + m_l d_l} \\
&= g_{s_1,\dots,s_l}(\tau) \,.
\end{align*}

The other special case $G^{R^l}_{s_1,\dots,s_l}$ can also be written down explicitly:
\begin{align*}
G^{R^l}_{s_1,\dots,s_l}(\tau) =  \sum_{\substack{m_1 =\dots=m_l=0 \\ n_1 > \dots > n_l > 0}} \frac{1}{(0\tau + n_1)^{s_1} \dots (0\tau + n_l)^{s_l}} = \zeta(s_1,\dots,s_l) \,.
\end{align*}
In length $2$ we have $ G_{s_1,s_2} = G^{RR}_{s_1,s_2}+G^{UR}_{s_1,s_2} +G^{RU}_{s_1,s_2}+G^{UU}_{s_1,s_2}$ and 
\begin{align*}
 G^{UR}_{s_1,s_2} &= \sum_{\substack{m_1>0, m_2=0 \\ n_1 \in \Z, n_2 > 0}} \frac{1}{(m_1 \tau + n_1)^{s_1} (0\tau + n_2)^{s_1}} \\
 &=\sum_{m_1>0 }\Psi_{s_1}(m_1 \tau) \sum_{n_2>0} \frac{1}{n_2^{s_2}} = g_{s_1}(\tau) \zeta(s_2)  \,,
\end{align*}
\[ G^{RU}_{s_1,s_2}(\tau) = \sum_{\substack{m_1=0, m_2>0 \\ n_1 > n_2 \\ n_i \in \Z}} \frac{1}{(m_1 \tau + n_1)^{s_1} (m_1 \tau + n_2)^{s_2}} = \sum_{m>0} \Psi_{s_1,s_2}(m\tau) .\]
In the case $G^{UR}$ we saw that we could write it as $G^U$ multiplied with a zeta value. 
In general, having a word $w$ of length $l$ ending in the letter $R$, i.e. there is a word $w'$ ending in $U$ with $w = w' R^r$ and $1 \leq r \leq l$ we can write 
\begin{align*}
G^{w}_{s_1,\dots,s_l}(\tau) =  G^{w'}_{s_1,\dots,s_{l-r}}(\tau) \cdot \zeta(s_{l-r+1},\dots,s_l) \,.
\end{align*}
{\bf Example:} $G^{RUURR}_{3,4,5,6,7} = G^{RUU}_{3,4,5} \cdot \zeta( 6,7)$ 

Hence one can concentrate on the words ending in $U$ when calculating the Fourier expansion of a multiple Eisenstein series.
Let $w$ be a word ending in $U$ then there are integers $r_1,\dots,r_j \geq 0$ with $w = R^{r_1} U R^{r_2} U \dots R^{r_j} U$. With this one can write
\[ G^w_{s_1,\dots,s_l}(\tau) = \sum_{m_1 > \dots > m_j > 0} \Psi_{s_1,\dots,s_{r_1+1}}(m_1 \tau) \cdot \Psi_{s_{r_1+2},\dots}(m_2 \tau) \dots \Psi_{s_{l-r_j},\dots,s_l}(m_j \tau) \,. \] 
{\bf Example:} $w=RURRU$
\[ G_{s_1,\dots,s_l}^{RURRU} = \sum_{m_1 > m_2 >0} \Psi_{s_1,s_2}(m_1 \tau) \Psi_{s_3,s_4,s_5}(m_2 \tau) \]

\begin{figure}[h!]
\begin{center}
\begin{tikzpicture}[scale=0.7]
\draw[densely dotted,step=1,color=gray,thin] (-4.9,-4.9) grid (4.9,4.9);
\draw [->] (0,-5.5) -- (0,5.5) node (yaxis) [above] {$m$};
\draw [->] (-5.5,0) -- (5.5,0) node (xaxis) [right] {$n$};
\draw[blue] (0,0) -- node[anchor=west] {$\lambda_5$} (1,1);
\draw[red] (1,1) -- node[anchor=north] {$\lambda_4$} (3,1); 
\draw[red] (3,1)-- node[anchor=north] {$\lambda_3$} (4,1);
\draw[blue] (4,1) --  node[anchor=west] {$\lambda_2$} (-1,4);
\draw[red] (-1,4) --  node[anchor=south] {$\lambda_1$} (2,4);
\fill[black] (1,1) circle (3pt);
\fill[black] (3,1) circle (3pt);
\fill[black] (4,1) circle (3pt);
\fill[black] (-1,4) circle (3pt);
\fill[black] (2,4) circle (3pt); 
\end{tikzpicture} \linebreak
A summand of $G^{\RR \UU \RR \RR \UU}_{s_1,\dots,s_l}$.
\end{center}
\end{figure}

{\bf Proof of Theorem \ref{thm:Fourier}:}
For $s_1,\dots,s_l \geq 2$ the Fourier expansion of the multiple Eisenstein series $G_{s_1,\dots,s_l}$ can be computed in the following way
\begin{enumerate}[i)]
\item Split up the summation into $2^l$ distinct parts $G^w_{s_1,\dots,s_l}$ where $w$ are a words in $\{R,U\}$. 
\item For $w$ being a word  ending in  $R$ one can write $G^w_{s_1,\dots,s_l}$  as $G^{w'}_{s_1,\dots} \cdot \zeta(\dots,s_l)$ with a word $w'$ ending in $U$. 
\item For $w$ being a word ending in $U$ one can write $G^w_{s_1,\dots,s_l}$ as
\[ G^w_{s_1,\dots,s_l}(\tau) = \sum_{m_1 > \dots > m_l >0} \Psi_{s_1,\dots}(m_1\tau) \dots \Psi_{\dots,s_l}(m_l \tau) \,. \]
\item Using the Theorem \ref{thm:reduction} we can write the  multitangent functions in iii) as a $\MZ$-linear combination of monotangents. We therefore just have $\MZ$-linear combinations with sums of the form
\[    \sum_{m_1 > \dots > m_j >0} \Psi_{k_1}(m_1 \tau) \dots \Psi_{k_j}(m_j \tau)= g_{k_1,\dots,k_j}(\tau) = (-2 \pi i)^{k_1+\dots+k_j} [k_1,\dots,k_j] \,. \]
\end{enumerate}
\qed

An explicit formula for the Fourier expansion of the multiple Eisenstein series for arbitrary length can be found in \cite{BT} Proposition 2.4. (with a reversed order of indices). Here we just give the Fourier expansion for the length 2 and 3.  For this we define for $n_1,n_2,k > 0$ the numbers $C_{n_1,n_2}^k$ by
\[ C_{n_1,n_2}^k = (-1)^{n_2} \binom{k-1}{n_2-1} + (-1)^{k-n_1} \binom{k-1}{n_1-1} \,. \]
\begin{prop}\label{prop:fourier23}
\begin{enumerate}[i)]
\item (\cite[Formula (52)]{GKZ}, \cite{Ba}, \cite{BT}) For $s_1,s_2 \geq 2$ the Fourier expansion of the double Eisenstein series is given by
\begin{align*}
G_{s_1,s_2}(\tau) = \zeta(s_1,s_2) + \zeta(s_2) g_1(\tau) + \sum_{\substack{k_1+k_2 = s_1+s_2 \\ k_2,k_2 \geq 2}} C^{k_2}_{s_1,s_2} \zeta(k_2) g_{k_1}(\tau) + g_{s_1,s_2}(\tau) \,.
\end{align*}
\item (\cite{Ba}, \cite{BT}) For $s_1,s_2,s_3 \geq 2$ and $k=s_1+s_2+s_3$ the Fourier expansion of the triple Eisenstein series can be written as
\begin{align*}
G_{s_1,s_2,s_3}(\tau) &=  \zeta(s_1,s_2,s_3) + \zeta(s_2,s_3) g_{s_1}(\tau) + \zeta(s_3) g_{s_1,s_2}(\tau) + g_{s_1,s_2,s_3}(\tau) \\
&+\zeta(s_3) \sum_{k_1+k_2 = s_1 + s_2} C_{s_1,s_2}^{k_1} \zeta(k_1) g_{k_2}(\tau) \\
&+\sum_{ k_1+k_2 = s_1 + s_2} C_{s_1,s_2}^{k_2}  \zeta(k_2) g_{k_1,s_3}(\tau)+\sum_{ k_1+k_2 = s_2 + s_3} C_{s_2,s_3}^{k_2} \zeta(k_2) g_{s_1,k_1}(\tau) \\
&+\sum_{k_1+k_2+k_3 = k} (-1)^{s_2+s_3} \binom{k_2-1}{s_2-1}\binom{k_3-1}{s_3-1} \zeta(k_3,k_2) g_{k_1}(\tau)   \\
&+\sum_{k_1+k_2+k_3 = k}(-1)^{s_1+s_2+k_2+k_3} \binom{k_2-1}{k_3-1}\binom{k_3-1}{s_2-1} \zeta(k_3,k_2) g_{k_1}(\tau) \\
&+(-1)^{s_1+s_3}\sum_{k_1+k_2+k_3 = k} (-1)^{k_2} \binom{k_2-1}{s_1-1}\binom{k_3-1}{s_3-1} \zeta(k_3) \zeta(k_2)  g_{k_1}(\tau) \,,
\end{align*}
where in the sums we sum over all $k_i \geq 2$. 
\end{enumerate}
\qed
\end{prop} 
We finish this section with a closer look at the stuffle product of two Eisenstein series. 
Since the product of multiple Eisenstein series can be written in terms of the stuffle product it is $G_2 \cdot G_3 = G_{2,3} + G_{3,2} + G_5$. On the other hand we have
\begin{align*}
G_2 \cdot G_3 &= \left(\zeta(2) + g_2 \right) \left( \zeta(3) + g_3 \right) = \zeta(2) \zeta(3) + \zeta(3) g_2 + \zeta(2) g_3 + g_2 \cdot g_3 \,. 
\end{align*}
and by Proposition \ref{prop:fourier23} it is 
\begin{align*}
G_{2,3} &= \zeta(2,3) - 2 \zeta(3) g_2 + \zeta(2) g_3 + g_{2,3} \,, \\
G_{3,2} &= \zeta(3,2) + 3 \zeta(3) g_2 + \zeta(2) g_3 + g_{3,2} \,.
\end{align*}
In conclusion, we obtain a relation for the product of the $g$'s namely $g_2 \cdot g_3 = g_{3,2} + g_{2,3} + g_5 + 2 \zeta(2) g_3$ and dividing out $(-2\pi i)^5$ we get
\[ [2] \cdot [3] = [3,2] + [2,3] + [5] - \frac{1}{12} [3] \,. \] 
We conclude that a product of the  $q$-series $[s_1,\dots,s_l] \in \Q[\![q]\!]$ has an expression similar to the stuffle product and that conversely, a product structure on these $q$-series could be used, together with the Fourier expansion, to explain the stuffle product for multiple Eisenstein series.

One might now ask, if the multiple Eisenstein series also "fulfill" the shuffle product. As we saw above the shuffle product of $\zeta(2)$ and $\zeta(3)$ reads
\begin{align} \label{eq:shuf23}
 \zeta(2) \cdot \zeta(3) =  \zeta(2,3) + 3 \zeta(3,2) + 6 \zeta(4,1)
\end{align}
and since there is no definition of $G_{4,1}$ this question does not make sense when replacing $\zeta$ by $G$ in \eqref{eq:shuf23}.
We will see that the understanding of the product structure of the brackets, explained in the next two sections, together with the Fourier expansion of multiple Eisenstein series will help to answer this question. This will be done by introducing shuffle regularized multiple Eisenstein series $G^\sh$ in Section \ref{sec:shuffreg}. There we will see that we can replace the $\zeta$ in \eqref{eq:shuf23} by $G^\sh$ and that the $G^\sh$ are given by the original $G$, for the cases in which they are defined. 

\section{Multiple divisor-sums and their generating functions}\label{section:bracket}

The classical divisor-sums $\sigma_r(n) = \sum_{d|n} d^r$ have a long history in number theory. They are well-known examples for multiplicative functions and appear in the Fourier expansion of Eisenstein series. This section is devoted to a larger class of functions, that can be seen as a multiple version of the divisor-sums and are therefore called multiple divisor-sums. For natural numbers $r_1,\dots,r_l \geq 0$ they are defined by 
\begin{equation} \label{def:sigma} \sigma_{r_1,\dots,r_l}(n) = \sum_{\substack{u_1 v_1 + \dots + u_l v_l = n\\u_1 > \dots > u_l >0}} v_1^{r_1} \dots v_l^{r_l} \,. \end{equation}
Even though the definition of these arithmetic functions is not complicated and somehow canonical, the author could not find any results on these functions before he started studying them in his master thesis \cite{Ba}. As mentioned in the introduction, the motivation to study them was due to their appearance in the Fourier expansion of multiple Eisenstein series, but as it turned out later in \cite{BK}, they are very nice and interesting objects in their own rights. Similar to multiple zeta values they fulfill a lot of relations. For example it is 
\begin{equation}\label{eq:reldivsum1} 
 \frac{1}{2} \sigma_2(n) =  \sigma_{1,0}(n) - \frac{1}{2} \sigma_1(n) + n\sigma_0(n) \,.
\end{equation}
Having objects of this type it is natural to consider their generating functions, which we denote by 
\begin{align*}
 [s_1,\dots,s_l] &:= \frac{1}{(s_1-1)! \dots (s_l-1)!} \sum_{n>0} \sigma_{s_1-1,\dots,s_l-1}(n) q^n \\
\end{align*}
and which are, just for the sake of short notations, called brackets. The factorial factors and the "shift" of $-1$ are natural if one thinks about the Fourier expansion of Eisenstein series. With this notation the relation \eqref{eq:reldivsum1} reads as 
\begin{equation} \label{eq:relbracket1} 
[3] = [2,1] - \frac{1}{2} [2] + q \frac{d}{dq} [1]\,,
\end{equation}
which can be seen as a counterpart of the relation $\zeta(3) = \zeta(2,1)$ between multiple zeta values\footnote{Further, one can prove the relation $\zeta(3) = \zeta(2,1)$ between multiple zeta values by multiplying both sides in \eqref{eq:relbracket1} with $(1-q)^3$ and then take the limit $q\rightarrow 1$. We will discuss this in Chapter $\ref{section:qana}$.}.\newline 

In this section, we want to focus on the algebraic structure of the space spanned by all brackets, which we will denote by $\MD$. This algebraic structure was studied in \cite{BK}. We will see that the space $\MD$ has the structure of a $\Q$-algebra and that the product of two brackets can be expressed in terms of brackets in a way that looks similar to the stuffle product of multiple zeta values. The operator $\dif =  q \frac{d}{dq}$ which appears in $\eqref{eq:relbracket1}$ plays an important role in the theory of (quasi-)modular forms. We will see that the space $\MD$ is closed under this operator and that this gives a second way of expressing the product of two brackets in length one similarly to the shuffle product of multiple zeta values. This second product expression in higher length will be discussed in Chapter \ref{section:bibrackets}.

\subsection{Brackets} \label{sec:brackets}
\begin{dfn}
For any integers $s_1,\dots,s_l>0$ we define the generating function for the multiple divisor sum $\sigma_{s_1-1,\dots,s_l-1}$ by the formal power series
\begin{align*}
 [s_1,\dots,s_l] &:= \frac{1}{(s_1-1)! \dots (s_l-1)!} \sum_{n>0} \sigma_{s_1-1,\dots,s_l-1}(n) q^n \\
 &=\sum_{\substack{u_1 > \dots > u_l> 0 \\ v_1, \dots , v_l >0}} \frac{v_1^{s_1-1} \dots v_l^{s_l-1}}{(s_1-1)!\dots(s_l-1)!}   \cdot q^{u_1 v_1 + \dots + u_l v_l}  \in \Q[\![q]\!] \,.
\end{align*}
\end{dfn}
In the first section, we saw that these series, by setting $q= \exp(2\pi i \tau)$, appear in the Fourier expansion of the multiple Eisenstein series but in this section we just view them as formal power series. We refer to these generating functions of multiple divisor sums as \emph{brackets} and define the vector space $\MD$ to be the $\Q$ vector space generated by $1 \in \Q[\![q]\!]$ and all 
brackets $[s_1,\dots,s_l]$. It is important to notice that we also include the constants in the space $\MD$.

\begin{ex} We give a few examples:
\begin{align*}
 [2] &= q + 3q^2 + 4q^3 + 7q^4 + 6q^5 + 12q^6 + 8q^7 + 15q^8 + \dots \,, \\
 [4,2] &=\frac{1}{6} \left( q^3 + 3q^4 + 15q^5 + 27q^6 + 78q^7 + 135q^8 +  \dots  \right)\,,  \\
 [4,4,4] &= \frac{1}{216} \left( q^6 + 9q^7 + 45q^8 + 190q^9 + 642q^{10} + 1899q^{11}  +\dots \right)\,,  \\
 [3,1,3,1] &= \frac{1}{4}\left( q^{10} + 2q^{11} + 8q^{12} + 16q^{13} + 43q^{14} + 70q^{15}+\dots \right) \,, \\
 [1,2,3,4,5]&=\frac{1}{288} \left( q^{15} + 17q^{16} + 107q^{17} + 512q^{18} + 1985q^{19} + \dots \right)\,.
\end{align*}
\end{ex}
Notice that the first non vanishing coefficient of $q^n$ in $[s_1,\dots,s_l]$ appears at \linebreak  $n = \frac{l (l+1)}{2}$, because it belongs to the "smallest" possible partition 
\[l \cdot 1 + (l-1) \cdot 1 + \dots + 1\cdot 1 = n \,, \]
i.e. $u_j= j$ and $v_j = 1$ for $1 \leq j \leq l$.  
The number $k=s_1+\dots+s_l$ is called the \emph{weight} of  $[s_1,\dots,s_l]$ and $l$ denotes the \emph{length}.\newline 

We want to show that the brackets are closed under multiplication by proving that their product structure is an example for a quasi-shuffle product. To do this we first introduce some notations and quote some results which are needed for this. 

Recall that for $s,z \in \C$, $|z|<1$ the polylogarithm $\Li_s(z)$ of weight $s$ is given by $\Li_s(z) = \sum_{n>0} \frac{z^n}{n^s}$.
For $s \in \N$ the $\Li_{-s}(z)$ are rational functions in $z$ with a pole in $z=1$. More precisely for $|z|<1$ they can be written as 
\[ \Li_{-s}(z) = \sum_{n>0} n^s z^n =\frac{z P_{s}(z)}{(1-z)^{s+1}}  \]
where $P_s(z)$ is the $s$-th Eulerian polynomial. Such a polynomial is given by
\[ P_s(X) = \sum_{n=0}^{s-1} A_{s,n} X^n \,, \]
where the  Eulerian numbers $A_{s,n}$ are defined by
\[ A_{s,n} = \sum_{i=0}^n (-1)^i \binom{s+1}{i} (n+1-i)^s\,.\]
For our purpose we write
\[ \Lit_{1-s}(z) := \frac{ \Li_{1-s}(z) }{(s-1)!}. \]
\begin{lem} (\cite[Lemma 2.5]{BK})\label{lem:eulerpol}
For $s_1,\dots,s_l \in \N$ we have
\begin{align*} [s_1,\dots,s_l] &= \sum_{n_1 > \dots > n_l>0} \Lit_{1-s_1}\left(q^{n_1}\right) \dots \Lit_{1-s_l}\left(q^{n_l}\right)\\
&=  \frac{1}{(s_1-1)! \dots (s_l-1)!} \sum_{n_1 > \dots > n_l > 0} \prod_{j=1}^l \frac{q^{n_j} P_{s_j-1}\left( q^{n_j} \right)}{(1-q^{n_j})^{s_j}} \,. 
\end{align*} 
 \qed
\end{lem}

\begin{rem}\begin{enumerate}[i)]
\item The second expression in terms of Eulerian Polynomials will be important for the interpretation of these series as $q$-analogues of multiple zeta values in Chapter \ref{section:qana}. 
\item This representation is also used for a fast implementation of these $q$-series in Pari GP. By doing so, the authors in \cite{BK} were able to give various results on the dimensions of the (weight and length filtered) spaces of $\MD$. These results can be found in Section 5 of \cite{BK}.
\end{enumerate}
\end{rem}

The product of $[s_1]$ and $[s_2]$ can thus be written as
\begin{align*}
 [s_1] \cdot [s_2] &= \left( \sum_{n_1 > n_2 > 0} + \sum_{n_2 > n_1 > 0} \right) \Lit_{1-s_1}\left(q^{n_1}\right)\Lit_{1-s_2}\left(q^n\right)  + \sum_{n_1 = n_2 > 0} \Lit_{1-s_1}\left(q^{n_1}\right) \Lit_{1-s_2}\left(q^{n_1}\right) \\
 &= [s_1,s_2] + [s_2,s_1] + \sum_{n> 0} \Lit_{1-s_1}\left(q^n\right) \Lit_{1-s_2}\left(q^n\right) \,.
\end{align*}
In order to prove that this product is an element of $\MD$ the product $\Lit_{1-s_1}\left(q^n\right) \Lit_{1-s_2}\left(q^n\right)$ must be  a rational linear combination of $\Lit_{1-j}\left(q^n\right)$ with $1 \leq j \leq s_1+s_2$. We therefore need the following
\begin{lem}  \label{lem_multpolylog}
For $a,b \in \N$ we have
\[ \Lit_{1-a}(z) \cdot  \Lit_{1-b}(z) = \sum_{j=1}^a \lambda^j_{a,b} \Lit_{1-j}(z) + \sum_{j=1}^b \lambda^j_{b,a}  \Lit_{1-j}(z)  + \Lit_{1-(a+b)}(z) \,, \]
where the coefficient $\lambda^j_{a,b}  \in \Q$ for $1 \leq j \leq a$ is given by
\[ \lambda^j_{a,b} = (-1)^{b-1} \binom{a+b-j-1}{a-j} \frac{B_{a+b-j}}{(a+b-j)!} \,, \]
with $B_k$ being the $k$-th Bernoulli number\footnote{For convenience we recall that the Bernoulli numbers $B_k$ are defined by $\frac{X}{e^X-1} =: \sum_{k\geq 0} \frac{B_k}{k!}X^k$.}.
\label{prop:productli}
\end{lem}
\begin{prf}
We prove this by using the generating function
\[ L(X) := \sum_{k>0} \Lit_{1-k}(z) X^{k-1} = \sum_{k>0} \sum_{n>0} \frac{n^{k-1} z^n}{(k-1)!} X^{k-1} = \sum_{n>0} e^{nX} z^n = \frac{e^X z}{1-e^X z} \,. \] 
With this one can see by direct calculation that
\[ L(X) \cdot L(Y) = \frac{1}{e^{X-Y}-1} L(X) + \frac{1}{e^{Y-X}-1} L(Y) \,. \]
By the definition of the Bernoulli numbers 
\[ \frac{X}{e^X-1} = \sum_{n\geq0} \frac{B_n}{n!} X^n  \]
this can be written as
\[ L(X) \cdot L(Y) = \sum_{n>0} \frac{B_n}{n!}(X-Y)^{n-1} L(X) + \sum_{n>0} \frac{B_n}{n!}(Y-X)^{n-1} L(Y) + \frac{L(X) - L(Y)}{X-Y}   \,. \]
The statement then follows by calculating the coefficient of $X^{a-1}Y^{b-1}$ in this equation.
\end{prf}

Now we are able to interpret the product structure of brackets as an example for a quasi-shuffle product. We equip $\h^1$ with a third product, beside the stuffle product $\ast$ and the shuffle product $\sh$. This product will be denoted $\boxast$, since it can be seen as a "bracket version" of the stuffle product $\ast$. For $a,b \in \N$ and  $w,v \in \h^1$ we define recursively the product 
\begin{align*}
\begin{split}
 z_a w \boxast z_b v &= z_a(w \boxast z_b v) + z_b(z_a w \boxast v)  + z_{a+b}( w \boxast v) + \sum_{j=1}^a \lambda^j_{a,b} z_j( w \boxast v) + \sum_{j=1}^b \lambda^j_{b,a} z_j( w \boxast v) \,,
\end{split}
\end{align*}
where the coefficients $\lambda^j_{a,b} \in \Q$ are the same as in Lemma \ref{lem_multpolylog}.
We equip $\MD$ with the usual multiplication of formal $q$-series and obtain the following: 

\begin{thm}(\cite[Prop 2.10]{BK})\label{thm:hoffmannalgebra}
For the linear map $[\, . \, ]: (\h^1, \boxast) \longrightarrow (\MD, \cdot)$ defined on the generators $w=z_{s_1}\dots z_{s_l}$ by $[w] := [s_1,\dots,s_l]$ we have
\[ [ w \boxast v ]    = [w] \cdot [v]  \]
and therefore $\MD$ is a $\Q$-algebra and $[\, . \, ]$ an algebra homomorphism. \qed
\end{thm}

\begin{ex} The first products of brackets are given by
\begin{align*}
[1] \cdot [1] &= 2[1,1] + [2] - [1] \,,  \\
[1] \cdot [2] &= [1,2]+[2,1]+[3] -\frac{1}{2}[2] \,,  \\
[1] \cdot [2,1] &= [1,2,1]+2 [2,1,1] - \frac{3}{2} [2,1]+[2,2]+[3,1] \,,\\
 [2] \cdot [3] &= [3,2] + [2,3] + [5] - \frac{1}{12} [3] \,,\\
 [3] \cdot [2,1] &= [3,2,1]+[2,3,1]+[2,1,3]+[5,1]+[2,4]+\frac{1}{12}[2,2]-\frac{1}{2}[2,3]-\frac{1}{12}[3,1] \,.
 \end{align*}
\end{ex}
We end this section by some notations which are needed for the rest of this paper. 
\begin{dfn}\label{def:filtration}
On $\MD$ we have the increasing filtration $\filw_{\bullet}$ given by the weight and the  increasing filtration $\fille_{\bullet}$ given by the length. For a subset $A\subset \MD$ we write\footnote{We set $[s_1,\dots,s_l] =1$ for $l=0$. }
\begin{align*}
\filw_k(A) &:=  \big<[s_1,\dots,s_l] \in A \,\big|\, l \geq 0\,, s_1+\dots+s_l \le k \,\big>_{\Q}\,,\\
\fille_l(A) &:=  \big<[s_1,\dots,s_r] \in A \,\big|\, 0 \le r\le l \,\big>_{\Q}\,.
\end{align*}
If we consider the length and weight filtration at the same time, we use the short notation $\filwle_{k,l} := \filw_k \fille_l$.
\end{dfn}

\begin{rem} As it can be seen by Theorem \ref{thm:hoffmannalgebra}, the multiplication of two brackets respects these filtrations, i.e. 
\[ \filwle_{k_1,l_1}(\MD) \cdot  \filwle_{k_2,l_2}(\MD) \subset \filwle_{k_1+k_2,l_1+l_2}(\MD). 
\label{prop:alg}\]
\end{rem}

\subsection{Derivatives and Subalgebras} \label{sec:dersubalb}

In this section we want to give an overview of interesting subalgebras of the space $\MD$ and discuss the differential structure with respect to the differential $\dif=q \frac{d}{dq}$. One of the main results in \cite{BK} is the following 
\begin{thm} (\cite[Thm. 1.7]{BK}) \label{thm:derivative} The operator $\dif = q \frac{d}{dq}$ is a derivation on $\MD$, it maps $\filwle_{k,l}(\MD)$ to $\filwle_{k+2,l+1}(\MD)$. \qed
\end{thm}
The proof of Theorem \ref{thm:derivative} uses generating functions of the brackets. It gives explicit formulas for the derivatives $\dif [s_1,\dots,s_l]$ for all $l$ which we omit here, since they are complicated. For example we have 
\begin{align*}
\dif [2,1,1] &= -\frac{1}{6} [2,1,1] + \frac{1}{2} [2,1,2] -[2,1,2,1]+[2,1,3]+ \frac{3}{2} [2,2,1] \\
 						 &-2\, [2,2,1,1]+[2,3,1]+6 [3,1,1]-8 [3,1,1,1]+[4,1,1]. 
\end{align*} 
In the following we give a list of subalgebras and review the results on whether they are also closed under $\dif$ or not. \newline 

{\bf i) (quasi-)modular forms}: 
Next to the connection to modular forms due to their appearance in the Fourier expansion of multiple Eisenstein series, the brackets have a direct connection to quasi-modular forms for $\Sl_2(\Z)$ with rational coefficients.
In the case $l=1$ we get the divisor sums $\sigma_{k-1}(n) = \sum_{d | n} d^{k-1}$ and
\[ [k] = \frac{1}{(k-1)!} \sum_{n>0} \sigma_{k-1}(n) q^n \,. \]
These simple brackets appear in the Fourier expansion of classical Eisenstein series with rational coefficients $\widetilde{G}_k(\tau):= (-2\pi i)^{-k}G_k(\tau)$ since we also included the rational numbers in $\MD$. For example we have
\[ \widetilde{G}_2 = -\frac{1}{24} + [2] \,,\quad \widetilde{G}_4 = \frac{1}{1440} + [4] \,, \quad \widetilde{G}_6 = -\frac{1}{60480} + [6] \,.\]
Denote by $M_\Q(\Sl_2(\Z)) = \Q[G_4,G_6]$ and $\widetilde{M}_\Q(\Sl_2(\Z)) = \Q[G_2,G_4,G_6]$ the algebras of modular forms and quasi-modular forms with rational coefficients.  

It is a well-known fact that the space $\widetilde{M}_\Q(\Sl_2(\Z))$ is closed under the operator $\dif = q \frac{d}{dq}$.\newline 

{\bf ii) Admissible brackets}: We define the set of all admissible brackets $\MDA$ as the span of all brackets $[s_1,\dots,s_l]$ with $s_1 > 1$ and $1$. This space is a subalgebra of $\MD$ (\cite[Thm. 2.13]{BK}) and every bracket can be written as a polynomial in the bracket $[1]$ with coefficients in $\MDA$: 
\begin{thm}(\cite[Thm. 2.14, Prop. 3.14]{BK}) \label{thm:polyadalg}
\begin{enumerate}[i)]
\item We have  $\MD = \MDA[\,[1]\,]$. 
\item The algebra $\MD$ is a polynomial ring over $\MDA$ with indeterminate $[1]$, i.e. $\MD$ is isomorphic to $\MDA[\,T\,]$ by sending $[1]$ to $T$. 
\item The space $\MDA$ is closed under $\dif$. \qed 
\end{enumerate}
\end{thm}
The elements in $\MDA$ are the ones, where the corresponding multiple zeta values exist. It will be reviewed in more detail in Chapter \ref{section:qana}, when we consider the brackets as $q$-analogues of multiple zeta values. \newline

{\bf iii) Even brackets and brackets with no $1$'s}:
Denote by $\MDE$ the space spanned by  $1$ and all $[s_1,\dots,s_l]$ with $s_j$ even for all $0 \leq j \leq l$ and by $\MD^\sharp$ the space spanned by  $1$ and all $[s_1,\dots,s_l]$ with $s_j>1$. Both spaces $\MDE$ and $\MD^\sharp$ are subalgebras of $\MD$ (\cite[Prop. 2.15]{BK}). It is expected, that the space $\MDE$ is not closed under $\dif$, since numerical calculation suggest, that for example $\dif [4,2] \notin \MDE$. 
Whether the space $\MD^\sharp$ is closed under this operator is an open and interesting question. In \cite{BK2} it is shown, that this is actually equivalent to one part of Conjecture 1 in \cite{O} given by Okounkov.

To summarize, we have the following inclusion of $\Q$-algebras
\begin{center}
\begin{tikzcd}[column sep=3ex]
M_\Q(\Sl_2(\Z)) \arrow[r,hook] \arrow[r,bend right=50,"\dif"]& 
\widetilde{M}_\Q(\Sl_2(\Z)) \arrow[loop,"\dif",out=125,in=45,looseness = 5] \arrow[r,hook]&
\MD^{ev}  \arrow[r,bend right=50,"\dif ?",dashrightarrow]  \arrow[r,hook] &
\MD^\#  \arrow[loop,"\dif ?",out=125,in=45,looseness = 5, dashrightarrow] \arrow[r,hook]&  
\MDA \arrow[loop,"\dif",out=125,in=45,looseness = 5]  \arrow[r,hook]&
\MD  \arrow[loop,"\dif",out=125,in=45,looseness = 5] &
\end{tikzcd}
\end{center}
The dashed arrows indicate the conjectured behavior of the map $\dif$, whereas the other arrows are all known to be correct. \newline

Though in length $l=1$ we derive not just one but several expressions for $\dif[s]$ given by the following Proposition. 
\begin{prop}(\cite[Prop 3.3]{BK}) \label{prop:formularfordk}
For $s_1,s_2$ with $s_1+s_2>2$ and $s=s_1+s_2-2$ we have the following expression for $\dif[s]$:
\[ \binom{s}{s_1-1} \frac{\dif[s]}{s} = [s_1]\cdot [s_2] +\binom{s}{s_1-1} [s+1] - \sum_{a+b=s+2} \left( \binom{a-1}{s_1-1}+\binom{a-1}{s_2-1} \right) [a,b] \,.\]
\end{prop}
 If you compare this formula with the shuffle product of multiple zeta values \eqref{eq:shufflelen1} in the length one times length one case you notice that Proposition \ref{prop:formularfordk} basically states that the brackets fulfill the shuffle product up to the term $\binom{s}{s_1-1} \frac{\dif[s]}{s} - \binom{s}{s_1-1} [s+1]$. 

We end this section by using these formulas to prove the following identity 

\begin{prop}\label{cor:delta}
The unique normalized cusp form $\Delta$ in weight $12$ can be written as
\begin{align*}
 -\frac{1}{2^6\cdot 5 \cdot 691}  \Delta  &=  168 [5,7]+150 [7,5]+28 [9,3] \notag \\
&+\frac{1}{1408} [2] - \frac{83}{14400}[4] +\frac{187}{6048} [6] - \frac{7}{120} [8] - \frac{5197}{691} [12] \,.
\end{align*}
\end{prop}
\begin{prf}
With the Eisenstein series $\widetilde{G}_6$ and $\widetilde{G}_{12}$ given by
\begin{align*}
 \widetilde{G}_6 &=(-2\pi i)^{-6} \zeta(6) + [6] =  -\frac{1}{60480} + [6] \,,\\
 \widetilde{G}_{12} &= (-2\pi i)^{-12} \zeta(12) + [12] = \frac{691}{2615348736000} + [12] \,,
\end{align*}
the cusp form $\Delta$ can be written as $\Delta = -3316800 G_6^2 + 3432000 G_{12}$.
Using quasi-shuffle product of brackets one can derive
\begin{align*}
\Delta &= \frac{3455}{198} [2] - \frac{691}{6} [4] +  \frac{6910}{21} [6] + 115200 [12] - 6633600 [6, 6] \,. 
\end{align*}
and therefore
\begin{align}\label{eq2}
\begin{split}
 -\frac{1}{2^6\cdot 5 \cdot 691} \Delta &=  30 [6, 6]- \frac{1}{12672}[2] + \frac{1}{1920}[4] -  \frac{1}{672}[6] -  \frac{360}{691} [12] \,. 
\end{split}
\end{align}
Using Proposition \ref{prop:formularfordk} for $(s_1,s_2) = (4,8), (5,7), (6,6)$ we get the following three expressions for $\dif[10]$
\begin{align*}
\dif[10]=& -\frac{1}{3} [5,7]- \frac{5}{6} [6,6]- \frac{5}{3} [7,5]- \frac{35}{12} [8,4]- \frac{16}{3} [9,3]-10 [10,2]-20 [11,1]  \\
&- \frac{1}{4790016} [2]+\frac{1}{403200}[4]-\frac{1}{36288}[6]+\frac{1}{8640}[8]+10 [11]+ \frac{1}{12}[12] \,, \\
\dif[10]=& -\frac{5}{21} [6,6]- \frac{5}{7} [7,5]-2 [8,4]- \frac{14}{3} [9,3]-10 [10,2]-20 [11,1]\\
 &+ \frac{1}{4790016}[2]- \frac{1}{604800} [4]+ \frac{1}{127008} [6]+10 [11]+ \frac{1}{21} [12] \,,\\
\dif[10]=& -\frac{10}{21} [7,5]- \frac{5}{3} [8,4]- \frac{40}{9} [9,3]-10 [10,2]-20 [11,1] \\
&- \frac{1}{4790016} [2]+ \frac{1}{725760} [4]- \frac{1}{381024} [6]+10 [11]+ \frac{5}{126} [12] \,.
\end{align*}
Summing them up as $0 = -504 \dif[10] +1890 \dif[10] -1386 \dif[10]$ we get
\begin{align}\label{eq1}
\begin{split}
 0 =&168 [5,7]-30 [6,6]+150 [7,5]+28 [9,3] \\
&+ \frac{5}{6336} [2] - \frac{181}{28800} [4]+  \frac{7}{216} [6] - \frac{7}{120} [8]   -7 [12]
\end{split}
\end{align}
Combining \eqref{eq1} and \eqref{eq2}, in order to eliminate the occurrence of $[6,6]$, we obtain the desired identity.
\end{prf}

\section{Bi-brackets and a second product expression for brackets}\label{section:bibrackets}

In the previous section we have seen that the space $\MD$ of brackets has the structure of a $\Q$-algebra and that there is an explicit formula to express the product of two brackets as a linear combination of brackets similarly to the stuffle product of multiple zeta values. 
In this section we want to present a larger class of $q$-series, called bi-brackets. The quasi-shuffle product of brackets extend to this larger class and therefore the space of bi-brackets is also a $\Q$-algebra. The beautiful feature of bi-brackets is, that there is a relation, which we call partition relation, which enables one to express the product of two bi-brackets in a second different way. These two product expressions then give a large class of linear relations, similar to the double shuffle relations of multiple zeta values. A variation of the bi-brackets were also studied in \cite{Z}. Later, the bi-brackets will be used to define regularized multiple Eisenstein series in Chapter \ref{section:regmes}. 
All results in this section were studied and introduced in \cite{Ba2}.

\subsection{Bi-brackets and their generating series} \label{sec:bibrackgen}

As motivated in the introduction of this section we want to study the following $q$-series:
\begin{dfn} For $r_1,\dots,r_l \geq 0$, $s_1,\dots,s_l > 0$ and we define the following $q$-series
\begin{align*}
 \mb{s_1, \dots , s_l}{r_1,\dots,r_l} :=&\sum_{\substack{u_1 > \dots > u_l> 0 \\ v_1, \dots , v_l >0}} \frac{u_1^{r_1}}{r_1!} \dots \frac{u_l^{r_l}}{r_l!} \cdot \frac{v_1^{s_1-1} \dots v_l^{s_l-1}}{(s_1-1)!\dots(s_l-1)!}   \cdot q^{u_1 v_1 + \dots + u_l v_l} \in \Q[\![q]\!]\,
 \end{align*}
which we call \emph{bi-brackets} of weight $r_1+\dots+r_k+s_1+\dots+s_l$, upper weight $s_1+\dots+s_l$, lower weight $r_1+\dots+r_l$ and length $l$. By $\bMD$ we denote the $\Q$-vector space spanned by all bi-brackets and $1$. 
\end{dfn}

The factorial factors in the definition of bi-brackets will become natural when considering generating functions of bi-brackets and the connection to multiple zeta values. 

For $r_1=\dots=r_l=0$ the bi-brackets are just the brackets 
\[  \mb{s_1, \dots , s_l}{0,\dots,0} = [s_1,\dots,s_l]  \]  
as defined in Chapter \ref{section:bracket}. Similarly to the Definition \ref{def:filtration} of the filtration for the space $\bMD$ we write for a subset $A\in \bMD$
\begin{align*}
\filw_k(A) &:=  \big< \mb{s_1,\dots,s_l}{r_1, \dots, r_l} \in A \,\big|\, 0 \leq l \leq k \,,\, s_1+\dots+s_l \le k \,\big>_{\Q}\\
\fild_k(A) &:=  \big< \mb{s_1,\dots,s_l}{r_1, \dots, r_l} \in A \,\big|\, 0 \leq l \leq k \,,\, r_1+\dots+r_l \le k \,\big>_{\Q}\\
\fille_l(A) &:=  \big<\mb{s_1,\dots,s_t}{r_1, \dots, r_t} \in A \,\big|\, t\le l \,\big>_{\Q}\,.
\end{align*}
and again if we consider the length and weight filtration at the same time we use the short notation $\filwle_{k,l} := \filw_k \fille_l$ and similar for the other filtrations.

\begin{prop}(\cite[Prop 4.2]{Ba2}) \label{prop:derivative}
Let $\dif := q \frac{d}{dq}$, then we have
\[ \dif \mb{s_1, \dots , s_l}{r_1,\dots,r_l} = \sum_{j=1}^l \left( s_j (r_j+1) \mb{s_1\,,\dots\,,s_{j-1}\,,s_j+1\,,s_{j+1},\dots \,,s_l}{r_1 \,, \dots \,,r_{j-1}\,,r_j+1\,,r_{j+1} \,,\dots \,, r_l} \right) \, \]
and therefore $\dif \left( \operatorname{Fil}^{\operatorname{W},\operatorname{D},\operatorname{L}}_{k,d,l}(\bMD) \right) \subset  \operatorname{Fil}^{\operatorname{W},\operatorname{D},\operatorname{L}}_{k+1,d+1,l}(\bMD) $.
\end{prop}
\begin{prf}
This is an easy consequence of the definition of bi-brackets and the fact that $\dif \sum_{n>0} a_n q^n = \sum_{n>0} n a_n q^n$.
\end{prf}

Proposition \ref{prop:derivative} suggests that the bi-brackets can be somehow viewed as  partial derivatives of the brackets with total differential $\dif$. \newline

In the following we now want to discuss the algebra structure of the space $\bMD$. For this we extend the quasi-shuffle product $\boxast$ of $\h^1$ to a larger space of words. Since we have double indices we replace the alphabet $A_z = \{z_1, z_2,\dots\}$ by  $A^{\text{bi}}_{z} := \{z_{s,r} \mid s\geq 1\,, r\geq 0\}$.

We consider on $\Q A^{\text{bi}}_{z}$ the commutative and associative product
\begin{align*}
 z_{s_1,r_1} \boxcircle z_{s_2,r_2} =&  \binom{r_1+r_2}{r_1} \sum_{j=1}^{s_1} \lambda^j_{s_1,s_2} z_{j,r_1+r_2} + \binom{r_1+r_2}{r_1}\sum_{j=1}^{s_2} \lambda^j_{s_2,s_1}  z_{j,r_1+r_2}  \\\
 &+\binom{r_1+r_2}{r_1} z_{s_1+s_2,r_1+r_2} 
\end{align*} 
and on $\Q\langle A^{\text{bi}}_{z} \rangle$ the commutative and associative quasi-shuffle product
 \begin{align*}
  z_{s_1,r_1} w \boxast z_{s_2,r_2} v &= z_{s_1,r_1}(w \boxast z_{s_2,r_2} v) + z_{s_2,r_2}(z_{s_1,r_1} w \boxast v)  +  (z_{s_1,r_1} \boxcircle z_{s_2,r_2})(w \boxast v) \,,  \\
\end{align*}
where the the numbers $\lambda^j_{a,b}  \in \Q$ for $1 \leq j \leq a$ are the same as before, i.e. 
\[ \lambda^j_{a,b} = (-1)^{b-1} \binom{a+b-j-1}{a-j} \frac{B_{a+b-j}}{(a+b-j)!} \,. \]

\begin{thm}(\cite[Thm. 3.6]{Ba2}) \label{thm:algebra}
The map $\mb{.}{}: (\Q\langle A^{\text{bi}}_{z}\rangle,\boxast) \rightarrow \left( \bMD , \cdot \right)$ given by
\[  w = z_{s_1,r_1} \dots  z_{s_l,r_l} \longmapsto [w] = \mb{s_1,\dots,s_l}{r_1,\dots,r_l} \]
fulfills $[w \boxast v]=[w] \cdot [v]$ and therefore $\bMD$ is a $\Q$-algebra.
\end{thm}

\begin{dfn}\label{def:genfct}
For the generating function of the bi-brackets we write
\[ \mt{ X_1,\dots,X_l}{Y_1,\dots,Y_l}:= \sum_{\substack{s_1,\dots,s_l > 0 \\ r_1, \dots, r_l > 0}} \mb{s_1\,,\dots\,,s_l}{r_1-1\,,\dots\,,r_l-1} X_1^{s_1-1} \dots X_l^{s_l-1} \cdot Y_1^{r_1-1} \dots Y_l^{r_l-1} \,. \]
These are elements in the ring $\bMDG =\varinjlim_j \bMD[[X_1,\dots,X_j,Y_1,\dots,Y_j]]$ of all generating series of bi-brackets.
\end{dfn}
To derive relations between bi-brackets we will prove functional equations for their generating functions. The key fact for this is that there are two different ways of expressing these given by the following Theorem.
\begin{thm}(\cite[Thm. 2.3]{Ba2}) \label{propgenfct}
For $n\in \N$ set
\[ E_n(X) := e^{nX}\,\quad \text{ and } \quad L_n(X) := \frac{e^X q^n}{1-e^X q^n} \in \Q[[q,X]]\,. \]
Then for all $l\geq 1$ we have the following two different expressions for the generating functions:
\begin{align*}
\mt{X_1,\dots,X_l}{Y_1,\dots,Y_l} &= \sum_{u_1>\dots>u_l>0} \prod_{j=1}^l E_{u_j}(Y_j) L_{u_j}(X_j) \\
&=  \sum_{u_1>\dots>u_l>0} \prod_{j=1}^l E_{u_j}(X_{l+1-j}-X_{l+2-j}) L_{u_j}(Y_1+\dots+Y_{l-j+1})
\end{align*}
(with $X_{l+1} := 0$). In particular the \emph{partition relations}\footnote{The bi-brackets and their generating series also give examples of what is called a \emph{bimould} by Ecalle in \cite{E}. In his language the partition relation \eqref{eq:partition} states that the bimould of generating series of bi-brackets is swap invariant.} holds:
\begin{equation} \label{eq:partition}
 \mt{X_1,\dots,X_l}{Y_1,\dots,Y_l} \overset{P}{=} \mt{ Y_1 + \dots + Y_l ,\dots ,Y_1 + Y_2, Y_1}{X_l , X_{l-1} - X_l, \dots , X_{1}-X_{2}}\,. 
\end{equation} \qed
\end{thm}

\begin{rem} 
A nice combinatorial explanation for the partition relation \eqref{eq:partition} is the following: By a partition of a natural number $n$ with $l$ parts we denote a representation of $n$ as a sum of $l$ distinct natural numbers, i.e. $15 = 4 + 4 + 3 + 2 + 1 + 1$ is a partition of $15$ with the $4$ parts given by $4,3,2,1$. We identify such a partition with a tuple $(u,v) \in \N^l \times \N^l$ where the $u_j$'s are the $l$ distinct numbers in the partition and the $v_j$'s count their appearance in the sum. The above partition of $15$ is therefore given by the tuple $(u,v) = ((4,3,2,1),(2,1,1,2))$. By $P_l(n)$ we denote all partitions of $n$ with $l$ parts and hence we set
\[ P_l(n) := \left\{ (u,v) \in \N^l \times \N^l \, \mid \, n = u_1 v_1 + \dots + u_l v_l \, \text{ and } \, u_1 > \dots > u_l > 0  \right\}   \] 
On the set $P_l(n)$ one has an involution given by the conjugation $\rho$ of partitions which can be obtained by reflecting the corresponding Young diagram across the main diagonal. 
\begin{figure}[H]
\[ {\tiny ((4,3,2,1),(2,1,1,2)) = \yng(4,4,3,2,1,1)  \quad  \overset{\rho}{\xrightarrow{\hspace*{1cm}}} \quad \yng(6,4,3,2) = ((6,4,3,2),(1,1,1,1))} \]
  \caption{The conjugation of the partition $15 = 4 + 4 + 3 + 2 + 1 + 1$ is given by $\rho( ((4,3,2,1),(2,1,1,2)) ) =  ((6,4,3,2),(1,1,1,1))$ which can be seen by reflection the corresponding Young diagram at the main diagonal.  }
\end{figure}
On the set $P_l(n)$ the conjugation $\rho$ is explicitly given by $\rho( (u,v) ) = (u',v')$ where $u'_j = v_1 + \dots + v_{l-j+1}$ and $v'_j = u_{l-j+1}-u_{l-j+2}$ with $u_{l+1} := 0$, i.e.  
\begin{equation} \label{eq:conj}
\rho : \binom{u_1, \dots , u_l}{v_1, \dots ,v_l} {\longmapsto} \binom{v_1+\dots+v_l, \dots,v_1+v_2,v_1}{u_l,u_{l-1}-u_l,\dots,u_1-u_2} \,.
\end{equation}
By the definition of the bi-brackets its clear that with the above notation they can be written as
\begin{align*}
 \mb{s_1, \dots , s_l}{r_1,\dots,r_l} :=& \frac{1}{r_1! (s_1-1)! \dots r_l! (s_l-1)!} \sum_{n>0} \left(  \sum_{(u,v) \in P_l(n)} u_1^{r_1} v_1^{s_1-1} \dots  u_l^{r_l} v_l^{s_l-1}  \right) q^n \,. 
 \end{align*}
The coefficients are given by a sum over all elements in $P_l(n)$ and therefore it is invariant under the action of $\rho$. As an example, consider $[2,2]$ and apply $\rho$ to the sum. Then we obtain 
\begin{align} \label{eq:2200ex}
\begin{split}
[2,2] &=  \sum_{n>0} \left( \sum_{(u,v) \in P_2(n)}  v_1 \cdot v_2 \right) q^n = \sum_{n>0} \left(\sum_{ \rho((u,v))=(u',v') \in P_2(n)}  u_2' \cdot (u_1'-u_2') \right) q^n \\
&= \sum_{n>0} \left( \sum_{(u',v') \in P_2(n)}  u_2' \cdot u_1' \right) q^n - \sum_{n>0} \left( \sum_{(u',v') \in P_2(n)}  u_2'^2 \right) q^n = \mb{1,1}{1,1} - 2 \mb{1,1}{0,2} \,.
\end{split}
\end{align}
This is exactly the relation one obtains by using the partition relation. 
\end{rem}

\begin{cor}(\cite[Cor. 2.5]{Ba2}) (Partition relation in length one and two) \label{cor:part12}
For $r,r_1,r_2 \geq 0$ and $s,s_1,s_2 > 0$ we have the  following relations in length one and two
\begin{align*}
\mb{s}{r} &= \mb{r+1}{s-1} \,,\\
\mb{s_1,s_2}{r_1,r_2} &= \sum_{\substack{0 \leq j \leq r_1\\0 \leq k \leq s_2-1}} (-1)^k \binom{s_1-1+k}{k} \binom{r_2+j}{j} \mb{r_2+j+1 \,,r_1-j+1}{s_2-1-k \,, s_1-1+k} \,.
\end{align*} \qed
\end{cor}

\begin{rem}\begin{enumerate}[i)]
\item \label{def:Pinvo} If we replace in the generating series in Definition \ref{def:genfct} the bi-brackets by the corresponding bi-words in and enforce the partition relation \eqref{eq:partition} for this power series, we obtain an involution 
\[ P:  \Q\langle A^{\text{bi}}_{z}\rangle \rightarrow  \Q\langle A^{\text{bi}}_{z}\rangle\,.\] 
By Corollary \ref{cor:part12} it is for example $P(z_{s,r}) = z_{r+1,s-1}$. This will be needed to describe the second product structure in the next section. 
\item In \cite{Z} the author introduces multiple q-zeta brackets $\mathfrak{Z}\mb{s_1,\dots,s_r}{r_1,\dots,r_l}$, which can be written in terms of bi-brackets and vice versa. For these objects the partition relation has the nice form 
\[\mathfrak{Z}\mb{s_1,\dots,s_r}{r_1,\dots,r_l} = \mathfrak{Z}\mb{r_l,\dots,r_1}{s_l,\dots,s_1}\,,\]
which can be interpreted in terms of duality. This is also used in \cite{Z} to describe the second product structure for the $\mathfrak{Z}$. Similarly in \cite{ems} the authors use a duality by Zhao (\cite{Zh}) to describe a second product structure for another model of $q$-analogues. 
\end{enumerate} 
\end{rem}

\subsection{Double shuffle relations for bi-brackets} \label{sec:dsh}

The partition relation together with the quasi-shuffle product can be used to obtain a second expression for the product of two bi-brackets. Before giving the general explanation this second product expression we illustrate it in two examples. 

\begin{ex} \begin{enumerate}[i)]
\item We want to given a second product expression for the product $[2] \cdot [3]$. By the partition relation we know that $[2] = \mb{1}{1}$, $[3] = \mb{1}{2}$ and using the quasi-shuffle product we have
\begin{align*}
 \mb{1}{1} \cdot \mb{1}{2} &=\mb{1,1}{1,2} + \mb{1,1}{2,1} -3 \mb{1}{3} + 3 \mb{2}{3} \,. \\
\end{align*}
The partition relations for the length two bi-brackets on the right is given by 
\begin{align*}
\mb{1,1}{1,2}&=\mb{3,2}{0,0}+3 \mb{4,1}{0,0} = [3,2] + 3 [4,1]\,,\\
\mb{1,1}{2,1}&=\mb{2,3}{0,0}+2\mb{3,2}{0,0}+ 3 \mb{4,1}{0,0}  = [2,3] + 2 [3,2] + 3 [4,1]\,.
\end{align*} 
Combining all of this we obtain 
\begin{align*}
 \mb{2}{0} \cdot \mb{3}{0} &=  \mb{1}{1} \cdot \mb{1}{2} \\
&=\mb{1,1}{1,2} + \mb{1,1}{2,1} -3 \mb{1}{3} + 3 \mb{2}{3} \\
&=[2,3] + 3[3,2] + 6 [4,1]+ 3\mb{4}{1} - 3 [4] \,.
\end{align*}
Compare this to the shuffle product of multiple zeta values 
\[ \zeta(2)\zeta(3)= \zeta(2,3) + 3 \zeta(3,2) + 6\zeta(4,1) \,. \] 
Since $\dif [3] = 3 \mb{4}{1}$ this example exactly coincides with the formula in Proposition \ref{prop:formularfordk} for the derivative $\dif[k]$.
\item In higher length, expressing the product of two bi-brackets in a similar way as in i) becomes interesting, since then the extra terms can't be expressed with the operator $\dif$ anymore. Doing the same calculation for the product $[3] \cdot [2,1]$, i.e. using the partition relation, the quasi-shuffle product and again the partition relation we obtain
\begin{align*}
[3] \cdot [2,1] &= \mb{1}{2} \cdot \mb{1,1}{0,1} \\
&= \mb{1,1,1}{2,0,1}+\mb{1,1,1}{0,2,1}+\mb{1,1,1}{0,1,2}+3 \mb{1,2}{0,3}+\mb{2,1}{2,1}-3 \mb{1,1}{0,3}-\mb{1,1}{2,1}\\
&=  [2,1,3]+[2,2,2]+2[2,3,1]+2[3,1,2]+5[3,2,1]+9[4,1,1]\\
&+\mb{2,3}{0,1} + 2 \mb{3,2}{0,1} + 3 \mb{4,1}{1,0} - [2,3] - 2 [3,2] -6[4,1]\,. 
\end{align*}
\end{enumerate}
This product can be seen as the analogue of the shuffle product 
\[ \zeta(3) \cdot \zeta(2,1) = \zeta(2,1,3)+\zeta(2,2,2)+2\zeta(2,3,1)+2\zeta(3,1,2) +5\zeta(3,2,1)+9\zeta(4,1,1) \,.\]
Here the bi-brackets, which are not given as brackets, can not be written in terms of the operator $\dif$ in an obvious way. 
\end{ex}
This works for arbitrary lengths and yields a natural way to obtain the second product expression for bi-brackets. To be more precise, denote by $P:  \Q\langle A^{\text{bi}}_{z}\rangle \rightarrow  \Q\langle A^{\text{bi}}_{z}\rangle$ the involution defined in Remark \ref{def:Pinvo}. Using this convention the second product expression for bi-brackets can be written in $ \Q\langle A^{\text{bi}}_{z}\rangle$ for two words $u , v \in \Q\langle A^{\text{bi}}_{z}\rangle$ as $ P\left( P(u) \boxast P(v) \right)$, i.e. the two product expressions of bi-brackets which correspond to the stuffle and shuffle product of multiple zeta values are given by
\begin{equation}\label{eq:stsh}
 [u] \cdot [v] = [ u \boxast v] \,, \qquad [u] \cdot [v] = [ P\left( P(u) \boxast P(v) \right) ] \,.
\end{equation}
In contrast to multiple zeta values these two product expression are the same for some cases, as one can check for the example $[1] \cdot [1,1]$. In the smallest length case, we have the following explicit formulas for the two products expressions.

\begin{prop}(\cite[Prop. 3.3]{Ba}) \label{prop:expstsh} 
For $s_1,s_2 > 0$ and $r_1,r_2 \geq 0$ we have the following two expressions for the product of two bi-brackets of length one:
\begin{enumerate}[i)]
\item ("Stuffle product analogue for bi-brackets")
\begin{align*}
\mb{s_1}{r_1} \cdot \mb{s_2}{r_2} &=\mb{s_1,s_2}{r_1,r_2} + \mb{s_2,s_1}{r_2,r_1} + \binom{r_1+r_2}{r_1}\mb{s_1+s_2}{r_1+r_2} \\
&+\binom{r_1+r_2}{r_1}\sum_{j=1}^{s_1} \frac{(-1)^{s_2-1} B_{s_1+s_2-j}}{(s_1+s_2-j)!} \binom{s_1+s_2-j-1}{s_1-j}  \mb{j}{r_1+r_2} \\
&+\binom{r_1+r_2}{r_1}  \sum_{j=1}^{s_2} \frac{(-1)^{s_1-1} B_{s_1+s_2-j}}{(s_1+s_2-j)!} \binom{s_1+s_2-j-1}{s_2-j} \mb{j}{r_1+r_2} 
\end{align*}

\item ("Shuffle product analogue for bi-brackets") 
\begin{align*}
\mb{s_1}{r_1} \cdot \mb{s_2}{r_2} &= \sum_{\substack{1 \leq j \leq s_1 \\ 0 \leq k \leq r_2}} \binom{s_1+s_2-j-1}{s_1-j} \binom{r_1+r_2-k}{r_1} (-1)^{r_2-k} \mb{s_1+s_2-j, j}{k,r_1+r_2-k} \\
&+\sum_{\substack{1\leq j \leq s_2 \\ 0 \leq k \leq r_1}} \binom{s_1+s_2-j-1}{s_1-1} \binom{r_1+r_2-k}{r_1-k} (-1)^{r_1-k} \mb{s_1+s_2-j, j}{k,r_1+r_2-k}\\
&+\binom{s_1+s_2-2}{s_1-1} \mb{s_1+s_2-1}{r_1+r_2+1} \\
&+\binom{s_1+s_2-2}{s_1-1}\sum_{j=0}^{r_1} \frac{(-1)^{r_2} B_{r_1+r_2-j+1}}{(r_1+r_2-j+1)!} \binom{r_1+r_2-j}{r_1-j}  \mb{s_1+s_2-1}{j} \\
&+\binom{s_1+s_2-2}{s_1-1}\sum_{j=0}^{r_2} \frac{(-1)^{r_1} B_{r_1+r_2-j+1}}{(r_1+r_2-j+1)!} \binom{r_1+r_2-j}{r_2-j}  \mb{s_1+s_2-1}{j} 
\end{align*}
\end{enumerate}
\end{prop}

Having these two expressions for the product of bi-brackets we obtain a large family of linear relations between them. Computer experiments suggest that actually every bi-bracket can be written in terms of brackets and that motivates the following surprising conjecture. 
\begin{conj}\label{conj:bdmd}
The algebra $\bMD$ of bi-brackets is a subalgebra of $\MD$ and in particular we have
\[ \operatorname{Fil}^{\operatorname{W},\operatorname{D},\operatorname{L}}_{k,d,l}(\bMD) \subset  \operatorname{Fil}^{\operatorname{W},\operatorname{L}}_{k+d,l+d}(\MD) \,. \]
\end{conj}
The results towards this conjecture, beside the computer experiments which have been done up to weight $8$, are the following 
\begin{prop}(\cite[Prop. 4.4]{Ba2}) For $l=1$ the Conjecture \ref{conj:bdmd} is true. 
\end{prop}
In \cite{BK3} it will be shown, that Conjecture \ref{conj:bdmd} is also true for all length up to weight $7$. For higher weights and lengths there are no general statements. The only general statement for the length two case is given by the following Proposition. 
\begin{prop}(\cite[Prop. 5.9]{Ba2}) \label{prop:len2conj}
For all $s_1,s_2 \geq 1$ it is
\begin{align*}
\mb{s_1,s_2}{1,0}, \mb{s_1,s_2}{0,1} \in \operatorname{Fil}^{\operatorname{W},\operatorname{L}}_{s_1+s_2+1,3}(\MD)
\end{align*} \qed
\end{prop}

\subsection{The shuffle brackets}

We now want to define a $q$-series which is an element in $\bMD$ and whose products can be written in terms of the "real" shuffle product of multiple zeta values. For $e_1,\dots,e_l \geq 1$ we generalize the generating function of bi-brackets to the following 
\begin{align}\label{def:tribracket} 
\mtt{X_1, & ... & ,\, X_l}{Y_1, & ... & ,\,  Y_l}{e_1, & ... & ,\, e_l} = \sum_{u_1>\dots>u_l>0} \prod_{j=1}^l E_{u_j}(Y_j) L_{u_j}(X_j)^{e_j} \,. 
\end{align}
In particular for $e_1 = \dots = e_l=1$ these are the generating functions of the bi-brackets. To show that the coefficients of these series are in $\bMD$ for arbitrary $e_j$ we need to define the differential operator $\mathcal{D}^Y_{e_1,\dots,e_l} := D_{Y_1,e_1}  D_{Y_2,e_2} \dots  D_{Y_l,e_l}$ with
\begin{align*}
D_{Y_j,e} &=  \prod_{k=1}^{e-1} \left(\frac{1}{k}\left(\frac{\partial}{\partial Y_{l-j+1}} - \frac{\partial}{\partial Y_{l-j+2}} \right) - 1\right)\,.
\end{align*}
where we set $\frac{\partial}{\partial Y_{l+1}}=0$. 

%

\begin{lem} Let $\mathcal{A}$ be an algebra spanned by elements $a_{s_1,\dots,s_l}$ with $s_1,\dots, s_l \in \N$, let $H(X_1,\dots,X_l) = \sum_{s_j} a_{s_1,\dots,s_l} X_1^{s_1-1} \dots X_1^{s_l-1}$ be the generating functions of these elements and define for $f \in \Q[[X_1,\dots,X_l]]$
\[ f^{\sharp}(X_1,\dots,X_l) = f(X_1+\dots+X_l, X_2+\dots+X_l, \dots, X_l) \,. \] 
Then the following two statements are equivalent.
\begin{enumerate}[i)]
\item The map $(\h^1, \sh) \rightarrow \mathcal{A}$ given by $z_{s_1} \dots z_{s_j} \mapsto a_{s_1,\dots,s_l}$ is an algebra homomorphism. 
\item For all $r,s \in \N$  it is
\[ H^{\sharp}(X_1,\dots,X_r) \cdot H^{\sharp}(X_{r+1}, \dots, X_{r+s}) = H^{\sharp}(X_1,\dots,X_{r+s})_{\vert sh_r^{(r+s)}}\,, \]
where $sh_r^{(r+s)} = \sum_{\sigma\in\Sigma(r,s)}\sigma$ in the group ring $\Z[\mathfrak{S}_{r+s}]$ and the symmetric group $\mathfrak{S}_r$ acts on $\Q[[X_1,\ldots,X_r]]$ by $(f\big|\sigma)(X_1,\ldots,X_r)= f(X_{\sigma^{-1}(1)},\ldots,X_{\sigma^{-1}(r)})$\,.
\end{enumerate}
\end{lem}
\begin{prf} This can be proven by induction over $l$ together with Proposition 8 in \cite{I}.\end{prf}

\begin{thm}(\cite[Thm. 5.7]{Ba2}) \label{thm:shufflebracket} For $s_1,\dots,s_l\in \N$ define $[s_1,\dots,s_l]^\sh \in \bMD$ as the coefficients of the following generating function
\begin{align*}
&H_{\sh}(X_1,\dots,X_l)  = \sum_{s_1,\dots,s_l \geq 1} [s_1,\dots,s_l]^\sh X_1^{s_1-1} \dots X_l^{s_l-1} \\
&:=\sum_{\substack{1\leq m \leq l\\i_1 + \dots + i_m = l}}\frac{1}{i_1! \dots i_m!} \mathcal{D}^Y_{i_1,\dots,i_m} \mt{X_1,X_{i_m+1},X_{i_{m-1}+i_m+1}, \dots, X_{i_2+\dots+i_m+1}}{Y_1,\dots,Y_l} _{\big\vert Y=0} \,.
\end{align*}
Then we have the following two statements
\begin{enumerate}[i)]
\item The $[s_1,\dots,s_l]^\sh$ fulfill the shuffle product, i.e. 
\[ H^{\sharp}_{\sh}(X_1,\dots,X_r) \cdot H^{\sharp}_{\sh}(X_{r+1}, \dots, X_{r+s}) = H^{\sharp}_{\sh}(X_1,\dots,X_{r+s})_{\vert sh_r^{(r+s)}} \,. \]
\item For $s_1\geq 1,\, s_2,\dots,s_l \geq 2$ we have $[s_1,\dots,s_l]^\sh = [s_1,\dots,s_l]$.
\end{enumerate} \qed
\end{thm} 
For low lengths we obtain the following examples: 
\begin{cor}\label{cor:explicitshufflebracket}
It is $[s_1]^\sh = [s_1]$ and for $l=2,3,4$  the $[s_1,\dots,s_l]^\sh$are given by\footnote{Here $\delta_{a,b}$ denotes the Kronecker delta, i.e $\delta_{a,b}$ is $1$ for $a=b$ and $0$ otherwise.}
\begin{enumerate}[i)]
\item $\begin{aligned}[t]
[s_1,s_2]^\sh &= [s_1,s_2] + \delta_{s_2,1}\cdot \frac{1}{2} \left( \mb{s_1}{1} - [s_1] \right) \,,\\
\end{aligned}$
\item $\begin{aligned}[t]
[s_1,s_2,s_3]^\sh &= [s_1,s_2,s_3]+ \delta_{s_3,1}\cdot \frac{1}{2} \left( \mb{s_1,s_2}{0,1} - [s_1,s_2] \right)\\
&+\delta_{s_2,1}\cdot \frac{1}{2}\left(  \mb{s_1,s_3}{1,0} - \mb{s_1,s_3}{0,1}- [s_1,s_3] \right)\\
&+\delta_{s_2 \cdot s_3,1}\cdot \frac{1}{6}\left(  \mb{s_1}{2} - \frac{3}{2}\mb{s_1}{1}+ [s_1] \right) \,,\\
\end{aligned}$
\item $\begin{aligned}[t]
[s_1,s_2,s_3&,s_4]^\sh = [s_1,s_2,s_3,s_4] + \delta_{s_4,1}\cdot \frac{1}{2} \left( \mb{s_1,s_2,s_3}{0,0,1} - [s_1,s_2,s_3] \right) \\
+ \delta_{s_3,1} &\cdot\frac{1}{2} \left( \mb{s_1,s_2,s_4}{0,1,0}-\mb{s_1,s_2,s_4}{0,0,1} + [s_1,s_2,s_4] \right)\\
+ \delta_{s_2,1} &\cdot\frac{1}{2} \left( \mb{s_1,s_3,s_4}{1,0,0}-\mb{s_1,s_3,s_4}{0,1,0} + [s_1,s_3,s_4] \right)\\
+ \delta_{s_2 \cdot s_4,1} &\cdot\frac{1}{4} \left( \mb{s_1,s_3}{1,1}-2\mb{s_1,s_3}{0,2}-\mb{s_1,s_3}{1,0} + [s_1,s_3] \right)\\
+ \delta_{s_3 \cdot s_4,1} &\cdot\frac{1}{6} \left( \mb{s_1,s_2}{0,2}-\frac{3}{2}\mb{s_1,s_2}{0,1} + [s_1,s_2] \right)\\
+ \delta_{s_2 \cdot s_3,1} &\cdot\frac{1}{6} \left( \mb{s_1,s_4}{0,2}-\mb{s_1,s_4}{1,1}+\frac{3}{2}\mb{s_1,s_4}{0,1}+\mb{s_1,s_4}{2,0}-\frac{3}{2}\mb{s_1,s_4}{1,0} + [s_1,s_4] \right)\\
+ \delta_{s_2 \cdot s_3 \cdot s_4,1} &\cdot\frac{1}{24} \left( \mb{s_1}{3}-2\mb{s_1}{2}+\frac{11}{6} \mb{s_1}{1} - [s_1] \right) \,.
\end{aligned}$
\end{enumerate}
\end{cor}
\begin{prf}
This follows by calculating the coefficients of the series $G_\sh$ in Theorem \ref{thm:shufflebracket}.
\end{prf}

The shuffle brackets will be used to define shuffle regularized multiple Eisenstein series in the next section. 

\section{Regularizations of multiple Eisenstein series}\label{section:regmes}

This section is devoted to Question \ref{qu1} in the introduction, which was to find a regularization of the multiple Eisenstein series. We want to present two type of regularization: The shuffle regularized multiple Eisenstein series (\cite{BT}, \cite{Ba2}) and stuffle regularized multiple Eisenstein series (\cite{Ba2}). 

The definition of shuffle regularized multiple Eisenstein series uses a beautiful connection of the Fourier expansion of multiple Eisenstein series and the coproduct of formal iterated integrals. The other regularization, the stuffle regularized multiple Eisenstein series uses the construction of the Fourier expansion of multiple Eisenstein series together with a  result on regularization of multitangent functions by O. Bouillot (\cite{Bo}). 

We start by reviewing the definition of formal iterated integrals and the coproduct defined by Goncharov. An explicit example in length two will make the above mentioned connection of multiple Eisenstein series and this coproduct clear. After doing this, we give the definition of shuffle and stuffle regularized multiple Eisenstein series as presented in \cite{BT} and \cite{Ba2}. At the end of this section we compare these two regularizations with a help of a few examples. 

\subsection{Formal iterated integrals} \label{sec:formaliteratedintegrals}
Following Goncharov (Section 2 in \cite{G}) we consider the algebra ${\mathcal I}$ generated by the elements
\[ \mathbb{I}(a_0;a_1,\ldots,a_N;a_{N+1}), \quad a_i\in\{0,1\}, N\ge0. \]
together with the following relations
\begin{enumerate}[i)]
\item For any $a,b\in\{0,1\}$ the unit is given by $\mathbb{I}(a;b):=\mathbb{I}(a;\emptyset;b)=1$.
\item The product is given by the shuffle product $\sh$
\begin{align*}
& \mathbb{I}(a_0;a_1,\ldots,a_M;a_{M+N+1}) \mathbb{I}(a_0;a_{M+1},\ldots,a_{M+N};a_{M+N+1})\\
&=\sum_{\sigma\in sh_{M,N}} \mathbb{I}(a_0;a_{\sigma^{-1}(1)},\ldots,a_{\sigma^{-1}(M+N)};a_{M+N+1}),
\end{align*}
where $sh_{M,N}$ is the set of $\sigma \in \mathfrak{S}_{M+N}$ such that $\sigma(1)<\cdots<\sigma(M)$ and $\sigma(M+1)<\cdots<\sigma(M+N)$.
\item The path composition formula holds: for any $N\ge0$ and $a_i,x\in\{0,1\}$, one has
\[ \mathbb{I}(a_0;a_1,\ldots,a_N;a_{N+1}) = \sum_{k=0}^N \mathbb{I}(a_0;a_1,\ldots,a_k;x)  \mathbb{I}(x;a_{k+1},\ldots,a_N;a_{N+1}).\]
\item For $N\ge1$ and $a_i,a\in\{0,1\}$ it is $\mathbb{I} (a;a_1,\ldots,a_N;a)=0$.
\item The path inversion is satisfied:
\[ \mathbb{I}(a_0;a_1,\ldots,a_N;a_{N+1}) = (-1)^N \mathbb{I}(a_{N+1};a_N,\ldots,a_1;a_0)\,. \]
\end{enumerate}

\begin{dfn}(Coproduct) 
Define the coproduct $\Delta$ on ${\mathcal I}$ by
\begin{align*} 
 &\Delta \left( \mathbb{I}(a_0;a_1,\ldots,a_N;a_{N+1}) \right) :=\\
 &\sum \left( \mathbb{I}(a_0;a_{i_1},\ldots,a_{i_k};a_{N+1}) \otimes \prod_{p=0}^k \mathbb{I}(a_{i_p};a_{i_p+1},\ldots,a_{i_{p+1}-1};a_{i_{p+1}}) \right)  ,
\end{align*}
where the sum on the right runs over all $i_0=0<i_1<\cdots<i_k<i_{k+1}=N+1$ with $0\le k \le N$.
\end{dfn}
\begin{prop}(\cite[Prop. 2.2]{G})
The triple $({\mathcal I},\sh, \Delta)$ is a commutative graded Hopf algebra over $\Q$. 
\end{prop}

To calculate $\Delta \left( \mathbb{I}(a_0;a_1,\dots,a_8;a_{9}) \right)$ one sums over all possible diagrams of the following form. 
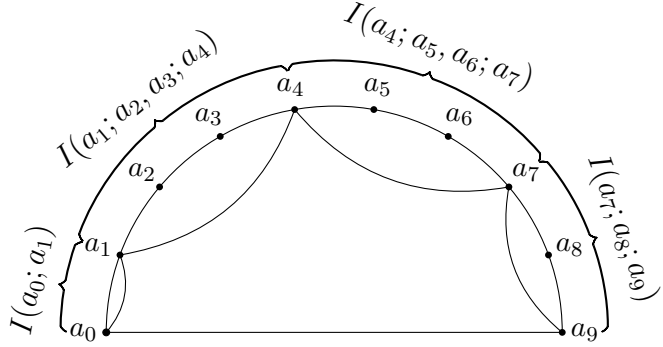
\begin{figure}[H]
\begin{center}
\input{dia0m1}
\caption{One diagram for the calculation of  $\Delta \left( \mathbb{I}(a_0;a_1,\dots,a_8;a_{9}) \right)$. It gives the term 
$I(a_0;a_1,a_4,a_7;a_9)  \otimes I(a_0;a_1)I(a_1;a_2,a_3;a_4)I(a_4;a_5,a_6;a_7)I(a_7;a_8;a_9) \,.$}
\end{center}
\end{figure}
For our purpose it will be important to consider the quotient space\footnote{If one likes to interpret the integrals as real integrals, then the passage from $\mathcal{I}$ to $\mathcal{I}^1$  regularizes these integrals such that "$-\log(0) = \int_{1 > t > 0} \frac{dt}{t} := 0$". }
\[ \mathcal{I}^1=\mathcal{I}/\mathbb{I}(1;0;0)\mathcal{I} \,.\]
Let us denote by 
\[ I(a_0;a_1,\ldots,a_N;a_{N+1})\]
an image of $\mathbb{I}(a_0;a_1,\ldots,a_N;a_{N+1})$ in $\mathcal{I}^1$. The quotient map $\mathcal{I} \rightarrow \mathcal{I}^1$ induces a Hopf algebra structure on $\mathcal{I}^1$, but for our application we just need that for any $w_1,w_2\in\mathcal{I}^1$, one has $\Delta(w_1\ \sh \ w_2)=\Delta(w_1)\ \sh \ \Delta(w_2)$. The coproduct on $\mathcal{I}^1$ is given by the same formula as before by replacing $\mathbb{I}$ with $I$.
For integers $n\ge0,s_1,\ldots,s_r\ge1$, we set 
\[I_{n}(s_1,\ldots,s_r):=I(1;\underbrace{0,0,\ldots,1}_{s_1},\ldots,\underbrace{0,0,\ldots,1}_{s_r},\underbrace{0,\ldots,0}_{n};0).\]
In particular, we write\footnote{This notion fits well with the iterated integral expression of multiple zeta values. Recall that 
\[ \zeta(2,3) = \int_{{\small 1 > t_1 > \dots > t_5 > 0}} \underbrace{\frac{dt_1}{t_1} \cdot \frac{dt_2}{1-t_2}}_{2} \cdot \underbrace{ \frac{dt_3}{t_3}\cdot  \frac{dt_4}{t_4} \cdot \frac{dt_5}{1-t_5}}_{3} \,. \]
This corresponds to $I(2,3)$ (but is of course not the same since the $I$ are formal symbols).} $I(s_1,\ldots,s_r)$ to denote $I_0(s_1,\ldots,s_r)$.
\begin{prop}(\cite[Eq. (3.5),(3.6) and Prop. 3.5]{BT})\label{prop:basisfori}
\begin{enumerate}[i)]
\item We have $I_n(\emptyset)=0$ if $n\ge1$ or $1$ if $n=0$.
\item For integers $n\ge0,s_1,\ldots,s_r\ge1$, 
\begin{equation*}
I_{n}(s_1,\ldots,s_r) = (-1)^n \sum^*\bigg(\prod_{j=1}^{r}\binom{k_j-1}{s_j-1} \bigg) I (k_1,\ldots,k_r) \,,
\end{equation*}
where the sum runs over all $k_1+\cdots+k_r=s_1+\cdots +s_r+n$ with $k_1,\ldots,k_r\ge1$.
\item The set $\{I(s_1,\ldots,s_r)\mid r\ge0,s_i\ge1\}$ forms a basis of the space $\mathcal{I}^1$.
\end{enumerate}
\end{prop}

We give an example for ii): In $\mathcal{I}^1$ it is $I(1;0;0)=0$ and therefore
\begin{align*}
0 &= I(1;0;0)  I(1;0,1;0) \\
&=I(1;0,0,1;0) + I(1;0,0,1;0) + I(1;0,1,0;0) \\
&= 2 I(3) + I_1(2) 
\end{align*}
which gives $I_1(2) = -2 I(3) = (-1)^1 \binom{2}{1} I(3)$. 

\begin{rem} \label{rmk:ih}
Statement iii) in Proposition \ref{prop:basisfori} basically states that we can identify $\mathcal{I}^1$ with $\h^1$ by sending $I(s_1,\dots,s_l)$ to $z_{s_1}\dots z_{s_l}$. In other words we can equip $\h^1$ with the coproduct $\Delta$. Instead of working with $I$ we will use this identification in the next section, when defining the shuffle regularized multiple Eisenstein series. 
\end{rem}

\begin{ex}\label{ex:coproductfourier32}
In the following we are going to calculate $\Delta(I(3,2)) = \Delta(I(1;0,0,1,0,1;0))$. Therefore we have to determine all possible markings of the diagram
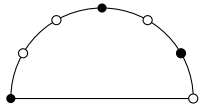
\begin{figure}[H]
\centering
\begin{tikzpicture}[scale=0.60]
	\def \ra {2}
	\def \circsize {0.1}		
	\draw (-\ra,0) -- (\ra,0)
				(\ra,0) arc (0:180:\ra);
	\fill[black] (180:\ra) circle (\circsize); 
	
	\fill[black] (30:\ra) circle (\circsize); 
	\draw (30:\ra) circle (\circsize); 	
	
	\fill[white] (60:\ra) circle (\circsize); 
	\draw (60:\ra) circle (\circsize); 	
	
	\fill[black] (90:\ra) circle (\circsize); 
	
	\fill[white] (120:\ra) circle (\circsize); 
	\draw (120:\ra) circle (\circsize); 
					
	\fill[white] (150:\ra) circle (\circsize);	
	\draw (150:\ra) circle (\circsize); 
	
	\fill[white] (0:\ra) circle (\circsize); 	
	\draw (0:\ra) circle (\circsize); 
\end{tikzpicture}
\end{figure}where the corresponding summand in the coproduct does not vanish. For simplicity we draw  $\circ$ to denote a $0$ and $\bullet$ to denote a $1$. We will consider the $4 = 2^2$ ways of marking the two $\bullet$ in the top part of the circle separately. As mentioned in the introduction, we want to compare the coproduct to the Fourier expansion of multiple Eisenstein series. Therefore, in this case we also calculate the expansion of $G_{3,2}(\tau)$ using the construction described in Section \ref{sec:mes}. Recall that we also had the $4$ different parts $G^{RR}_{3,2}$, $G^{UR}_{3,2}$, $G^{RU}_{3,2}$ and $G^{UU}_{3,2}$. We will see that the number and positions of the marked $\bullet$ correspond to the number and positions of the letter $U$ in the word $w$ of $G^w$. 
\begin{enumerate}[i)]
\item Diagrams with no marked $\bullet$: 
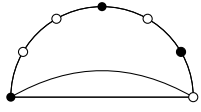
\begin{figure}[H]
\centering
\include{d32_0}
\end{figure}
Corresponding sum in the coproduct: 
\[ I(0;\emptyset;1) \otimes I(1;0,0,1,0,1;0) = 1 \otimes I(3,2) \,. \] 
The part of the Fourier expansion of $G_{3,2}$ which is associated to this, is the one with no $U$ "occurring", i.e.  $G^{RR}_{3,2}(\tau) = \zeta(3,2)$.  
\item  Diagrams with the first $\bullet$ marked: 
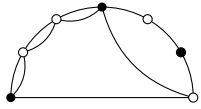
\begin{figure}[H]
\centering
\include{d32_2}
\end{figure}
Corresponding sum in the coproduct: 
\begin{align*}
 I(1;0,0,1;0) \otimes \big( I(1;0) \cdot I(0;0) \cdot I(0;1) \cdot I(1;0,1;0) \big) = I(3) \otimes I(2) \,.
\end{align*}
The associated part of the Fourier expansion of $G_{3,2}$ is $G^{UR}_{3,2}(\tau) = g_3(\tau) \cdot \zeta(2)$.  
\item Diagrams with the second $\bullet$ marked: 
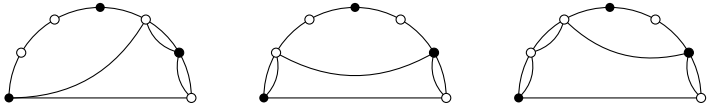
\begin{figure}[H]
\centering
\include{d32_1}
\end{figure}
Corresponding sum in the coproduct: 
\begin{align*}
 & I(1;0,1;0) \otimes \big(  I(1;0,0,1;0) \cdot I(0;1) \cdot I(1;0)\big) \\
+ & I(1;0,1;0) \otimes \big(  I(1;0) \cdot I(0;0,1,0;1) \cdot I(1;0)\big) \\
+ & I(1;0,0,1;0) \otimes \big(  I(1;0) \cdot I(0;0) \cdot I(0;1,0;1) \cdot I(1;0)\big) \\
& = I(2) \otimes I(3) - I(2) \otimes I_1(2) + I(3) \otimes I(2) \,,
\end{align*}
where we used $I(0,0,1,0;1)=-I_1(2)$ and $I(0;1,0;1) = (-1)^2 I(1;0,1;0) = I(2)$. Together with $I_1(2) = -2 I(3)$ this gives
\[ 3 I(2) \otimes I(3) + I(3) \otimes I(2) \,.\]
Also the associated part of the Fourier expansion is the most complicated one. We had $G^{RU}_{3,2}(\tau) = \sum_{m>0} \Psi_{3,2}(m\tau)$ and with $\eqref{eq:psi32}$ we derived $\Psi_{3,2}(x) = 3 \Psi_2(x) \cdot \zeta(3) + \Psi_3(x) \cdot \zeta(2)$, i.e. 
\[ G^{RU}_{3,2}(\tau) = 3 g_2(\tau) \cdot \zeta(3) + g_3(\tau) \cdot \zeta(2) \,.\] 
\item  Diagrams with both $\bullet$ marked: 
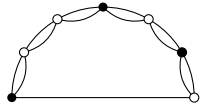
\begin{figure}[H]
\centering
\include{d32_12}
\end{figure}
Corresponding sum in the coproduct: $ I(3,2) \otimes 1$. The associated part of the Fourier expansion of $G_{3,2}$ is $G^{UU}_{3,2}(\tau) = g_{3,2}(\tau)$.  
\end{enumerate}
Summing all $4$ parts together we obtain for the coproduct
\[ \Delta( I(3,2) ) =  1 \otimes I(3,2)   +  3 I(2) \otimes I(3) + 2 I(3) \otimes I(2)+  I(3,2)  \otimes 1  \]
and for the Fourier expansion of $G_{2,3}(\tau)$: 
\begin{align*}
G_{3,2}(\tau) =\zeta(3,2) + 3  g_2(\tau) \zeta(3) + 2  g_3(\tau) \zeta(2) + g_{3,2}(\tau) \,. 
\end{align*}
\end{ex}

This shows that the left factors of the terms in the coproduct corresponds to the functions $g$ and the right factors side to the multiple zeta values. We will use this in the next section to define shuffle regularized multiple Eisenstein series. 

\subsection{Shuffle regularized multiple Eisenstein series}\label{sec:shuffreg}
In this section we present the definition of shuffle regularized multiple Eisenstein series as it was done in \cite{BT} together with the simplification developed in \cite{Ba2}. We use the observation of the section before and use the coproduct $\Delta$ of formal iterated integrals to define these series. As mentioned in Remark \ref{rmk:ih} we can equip the space $\h^1$ with the coproduct $\Delta$ instead of working with the space  $\mathcal{I}^1$.
Denote by $\MZB \subset \C[\![q]\!]$ the space of all formal power series in $q$ which can be written as a $\Q$-linear combination of products of multiple zeta values, powers of $(-2\pi i)$ and bi-brackets. In the following, we set $q=\exp(2\pi i \tau)$ with $\tau$ being an element in the upper half-plane. Since the coefficient of bi-brackets just have polynomials growth, the elements in $\MZB$ and $\bMD$ can be viewed as holomorphic functions in the upper half-plane with this identification.

In analogy to the map $Z^\sh: (\h^1, \sh) \rightarrow \MZ$ of shuffle regularized multiple zeta values (Proposition \ref{prop:mzvreg}), the map  $\mathfrak{g}^{\sh} : (\h^1, \sh) \rightarrow \Q[2\pi i][\![q]\!]$ defined on the generators $z_{t_1} \dots z_{t_l}$ by 
\[ \mathfrak{g}^{\sh}(z_{t_1} \dots z_{t_m}) = g^\sh_{t_1,\ldots,t_m}(\tau) := (-2\pi i)^{t_1+\dots+t_m}[t_1,\ldots,t_m]^\sh\,, \] 
is also an algebra homomorphism by Theorem \ref{thm:shufflebracket}. 

With this notation we can recall the definition of $G^\sh$ from \cite{Ba2} (which is a variant of the definition in \cite{BT}, where the authors did not use bi-brackets and the shuffle bracket).
\begin{dfn}
For integers $s_1,\ldots,s_l\ge1$, define the functions $G^\sh_{s_1,\ldots,s_l}(\tau)\in \MZB$, called  \emph{shuffle regularized multiple Eisenstein series}, as 
\[ G^\sh_{s_1,\ldots,s_l}(\tau) := m\left( (  \mathfrak{g}^{\sh} \otimes  Z^{\sh})\circ \Delta \big(z_{s_1} \dots z_{s_l} \big)\right)\,,\]
where $m$ denotes the multiplication given by $m: a \otimes b \mapsto a\cdot b$ and $Z^\sh$ denotes the map for shuffle regularized multiple zeta values given in Proposition \ref{prop:mzvreg}.
\end{dfn}
We can view $G^\sh$ as an algebra homomorphism $G^\sh: (\h^1, \sh) \rightarrow \MZB$ such that the following diagram commutes
\[
\xymatrix{  
(\h^1,\sh)  \ar@{->}[r]^-{\Delta} \ar@{->}[d]_{G^\sh} & (\h^1,\sh)  \otimes (\h^1,\sh)  \ar@{->}[d]^{\mathfrak{g}^{\sh} \otimes  Z^{\sh} }  \\ 
 \MZB &   \Q[2\pi i][\![q]\!] \otimes \,\MZ  \ar@{->}[l]^-{m}\\
} 
\]
\begin{thm}(\cite[Thm. 6.5 ]{Ba2}, \cite[Thm. 1.1, 1.2]{BT}) \label{thm:messh}
For all $s_1,\ldots,s_l\ge1$ the shuffle regularized multiple Eisenstein series $G^{\sh}_{s_1,\ldots,s_l}$ have the following properties:
\begin{enumerate}[i)]
\item They are holomorphic functions on the upper half-plane having a Fourier expansion with the shuffle regularized multiple zeta values as the constant term. 
\item They fulfill the shuffle product.
\item For integers $s_1,\ldots,s_l\ge2$ they equal the multiple Eisenstein series
\[ G^{\sh}_{s_1,\ldots,s_l}(\tau)=G_{s_1,\ldots,s_l}(\tau) \]
and therefore they fulfill the stuffle product in these cases. 
\end{enumerate} \qed
\end{thm}
Parts i) and ii) in this theorem follow directly by definition. The important part here is iii), which states that the connection of the Fourier expansion and the coproduct, as illustrated in Example \ref{ex:coproductfourier32}, holds in general. It also proves that the shuffle regularized multiple Eisenstein series fulfill the stuffle product in many cases. Though the exact failure of the stuffle product of these series is unknown so far. 

\subsection{Stuffle regularized multiple Eisenstein series}\label{sec:stufflereg}
Motivated by the calculation of the Fourier expansion of multiple Eisenstein series described in Section \ref{sec:mes} we consider the following construction. 

\begin{constr} \label{const}
Given a $\Q$-algebra $(A,\cdot)$ and a family of homomorphism
\[ \left\{ w \mapsto f_w(m) \right\}_{m \in \N} \]
 from  $(\h^1, \ast)$ to $(A,\cdot)$, we define for $w \in \h^1$ and $M \in \N$ 
\[ F_w(M) := \sum_{\substack{1 \leq k \leq l(w)\\w_1 \dots w_k = w\\ M> m_1 > \dots >m_k >0}}  f_{w_1}(m_1) \dots f_{w_k}(m_k) \in A \,, \]
where $l(w)$ denotes the length of the word $w$ and $w_1 \dots w_k = w$ is a decomposition of $w$ into $k$ words in $\h^1$. 
\end{constr}
\begin{prop}(\cite[Prop. 6.8]{Ba2})\label{prop:construction}
For all $M \in \N$ the assignment $w \mapsto F_w(M)$, described above, determines an algebra homomorphism from $(\h^1, \ast)$ to $(A,\cdot)$. In particular $\left\{ w \mapsto F_w(m) \right\}_{m \in \N}$ is again a family  of homomorphism as used in Construction \ref{const}. \qed
\end{prop}

For a word $w=z_{s_1} \dots z_{s_l} \in \h^1$ we also write in the following $f_{s_1,\dots,s_l}(m):=f_w(m)$ and similarly $ F_{s_1,\dots,s_l}(M):=F_w(M)$. 

\begin{ex} Let $f_w(m)$ be as in Construction \ref{const}. In small lengths the $F_w$ are given by
\[ F_{s_1}(M) = \sum_{M> m_1 >0} f_{s_1}(m_1) \,, \quad F_{s_1,s_2}(M) = \sum_{M> m_1 >0} f_{s_1,s_2}(m_1) + \sum_{M> m_1 > m_2>0} f_{s_1}(m_1)  f_{s_2}(m_2) \,\]
and one can check directly by the use of the stuffle product for the $f_w$ that 
\begin{align*}
&F_{s_1}(M) \cdot F_{s_2}(M) =  \sum_{M> m_1>0} f_{s_1}(m_1) \cdot  \sum_{M> m_2 >0} f_{s_2}(m_2) \\
&=   \sum_{M> m_1 > m_2>0} f_{s_1}(m_1)  f_{s_2}(m_2) +  \sum_{M> m_2 > m_1>0}   f_{s_2}(m_2) f_{s_1}(m_1) +  \sum_{M>  m_1> 0}   f_{s_1}(m_1) f_{s_2}(m_1)  \\
&=  \sum_{M> m_1 > m_2>0} f_{s_1}(m_1)  f_{s_2}(m_2) +  \sum_{M> m_2 > m_1>0}   f_{s_2}(m_2) f_{s_1}(m_1) \\
&+ \sum_{M> m_1>0} \left(f_{s_1,s_2}(m_1) +  f_{s_2,s_1}(m_1) + f_{s_1+s_2}(m_1) \right)  \\
&=  F_{s_1,s_2}(M)+ F_{s_2,s_1}(M) + F_{s_1+s_2}(M)  \,.
\end{align*}
\end{ex}

Let us now give an explicit example for maps $f_w$ in which we are interested. Recall (Definition \ref{def:multitangent}) that for integers $s_1,\dots,s_l \ge 2$ we defined the multitangent function by
\[ \Psi_{s_1,\ldots,s_l}(z) = \sum_{\substack{n_1>\cdots>n_l\\ n_j \in \Z}} \frac{1}{(z+n_1)^{s_1}\cdots (z+n_l)^{s_l}}.\]
In \cite{Bo}, where these functions were introduced, the author uses the notation $\mathcal{T}e^{s_1,\ldots,s_l}(z)$ which corresponds to our notation $\Psi_{s_1,\ldots,s_l}(z)$. It was shown there that the series $\Psi_{s_1,\ldots,s_l}(z)$ converges absolutely when $s_1,\ldots,s_l\ge2$. These functions fulfill (for the cases they are defined) the stuffle product. As explained in Section \ref{sec:mes} the multitangent functions appear in the calculation of the Fourier expansion of the multiple Eisenstein series $G_{s_1,\dots,s_l}$, for example in length two it is
\[ G_{s_1,s_2}(\tau) = \zeta(s_1,s_2) + \zeta(s_1) \sum_{m_1 > 0} \Psi_{s_2}(m_1 \tau) + \sum_{m_1>0} \Psi_{s_1,s_2}(m_1 \tau) +  \sum_{m_1 > m_2 > 0} \Psi_{s_1}(m_1\tau) \Psi_{s_2}(m_2\tau) \,.\] 
One nice result of \cite{Bo} is a regularization of the multitangent function to get a definition of $\Psi_{s_1,\ldots,s_l}(z)$ for all $s_1,\dots,s_l \in \N$.  We will use this result together with the above construction to recover the Fourier expansion of the multiple Eisenstein series. 
\begin{thm}(\cite{Bo})\label{thm:olivier}
For all $s_1,\dots,s_l \in \N$  there exist holomorphic functions $\Psi_{s_1,\dots,s_l}$ on $\Ha$ with the following properties
\begin{enumerate}[i)]
\item Setting $q=e^{2\pi i \tau}$ for $\tau \in \Ha$ the map  $w \mapsto \Psi_w(\tau)$ defines an algebra homomorphism from $(\h^1, \ast)$ to $(\C[\![q]\!],\cdot)$.
\item In the case $s_1,\dots,s_l \geq 2$ the  $\Psi_{s_1,\dots,s_l}$ are given by the multitangent functions in Definition \ref{def:multitangent}.
\item The monotangents functions have the $q$-expansion given by
\[ \Psi_1(\tau) = \frac{\pi}{\tan(\pi \tau)} =  (-2 \pi i)\left(\frac{1}{2} + \sum_{n>0} q^n \right),\quad  \Psi_k(\tau) =  \frac{(-2\pi i)^k}{(k-1)!} \sum_{n>0} n^{k-1} q^n \, \text{ for } k\geq 2.\] 
\item (Reduction into monotangent function) Every $\Psi_{s_1,\dots,s_l}(\tau)$ can be written as a $\MZ$-linear combination of monotangent functions. There are explicit $\epsilon^{s_1,\dots,s_l}_{i,k} \in \MZ$ s.th. 
\[\Psi_{s_1,\dots,s_l}(\tau) = \delta^{s_1,\dots,s_l} + \sum_{i=1}^l \sum_{k=1}^{s_i} \epsilon^{s_1,\dots,s_l}_{i,k} \Psi_k(\tau) \,,\]
where $ \delta^{s_1,\dots,s_l} = \frac{(\pi i)^l}{l!}$ if $s_1=\dots=s_l=1$ and $l$ even and $\delta^{s_1,\dots,s_l} = 0$ otherwise. 
For $s_1>1$ and $s_l>1$ the sum on the right starts at $k=2$, i.e. there are no $\Psi_1(\tau)$ appearing and therefore there is no constant term in the $q$-expansion.
\end{enumerate}
\end{thm}
\begin{prf}
This is just a summary of the results in Section $6$ and $7$ of \cite{Bo}. The last statement iv) is given by Theorem 6 in \cite{Bo}. 
\end{prf} 
Due to iv) in the Theorem the calculation of the Fourier expansion of multiple Eisenstein series, where ordered sums of multitangent functions appear, reduces to ordered sums of monotangent functions. The connection of these sums to the brackets, i.e. to the functions $g$, is given by the following fact which can be seen by using iii) of the above Theorem. For $n_1,\dots,n_r \geq 2$ it is
\begin{align*}
&g_{s_1,\ldots,s_r}(\tau) = \sum_{m_1 >\dots >m_l>0} \Psi_{s_1}(m_1\tau) \dots \Psi_{s_l}(m_l\tau) \,.
\end{align*}
For $w \in \h^1$ we now use the Construction \ref{const} with $A=\C[\![q]\!]$ and the family of homomorphism $\{ w \mapsto \Psi_w(n \tau) \}_{n\in \N}$ (See Theorem \ref{thm:olivier} i)\,) to define
\[ \mathfrak{g}^{\ast,M}(w) := (-2\pi i)^{|w|} \sum_{\substack{1 \leq k \leq l(w)\\w_1 \dots w_k = w}} \sum_{M>m_1 > \dots > m_k >0} \Psi_{w_1}(m_1 \tau) \dots \Psi_{w_k}(m_k \tau)  \,. \] 
From Proposition \ref{prop:construction} it follows that for all $M\in \N$ the map $\mathfrak{g}^{\ast,M}$ is an algebra homomorphism from $(\h^1,\ast)$ to $\C[\![q]\!]$.

To define stuffle regularized  multiple Eisenstein series we need the following: For an arbitrary quasi-shuffle algebra $\Q\langle A \rangle$ define the following coproduct for a word $w$ 
\[ \Delta_H(w) = \sum_{u v = w} u \otimes v\,. \] 
Then it is known due to Hoffman (\cite{H}) that the space $\left( \Q\langle A \rangle, \odot , \Delta_H \right)$ has the structure of a bialgebra. With this we try to mimic the definition of the $G^\sh$ and use the coproduct structure on the space $(\h^1, \ast, \Delta_H)$ to define for $M\geq 0$ the function $G^{\ast,M}$ and then take the limit $M\rightarrow \infty$ to obtain the stuffle regularized multiple Eisenstein series. For this we consider the following diagram
\[
\xymatrix{  
 (\h^1,\ast)  \ar@{->}^-{\Delta_H} [r]\ar@{->}[d]_{G^{\ast,M}}  &(\h^1,\ast) \otimes (\h^1,\ast) \ar@{->}[d]^{\mathfrak{g}^{\ast,M} \otimes\, Z^{\ast} }   \\ 
 \C[\![q]\!] & \C[\![q]\!] \otimes  \MZ  \ar@{->}[l]^-{m}   &\\
} \]
with the above algebra homomorphism $\mathfrak{g}^{\ast,M} : (\h^1, \ast) \rightarrow \C[\![q]\!]$ and the map  $Z^\ast$ for stuffle regularized multiple zeta values given in Proposition \ref{prop:mzvreg}. 

\begin{dfn}\label{def:gast}
For integers $s_1,\ldots,s_l\ge1$ and $M \geq 1$, we define the $q$-series $G^{\ast,M}_{s_1,\ldots,s_r}\in\C[\![q]\!]$ as the image of the word $w= z_{s_1}\dots z_{s_l} \in \h^1$ under the algebra homomorphism $(\mathfrak{g}^{\ast,M} \otimes Z^{\ast})\circ \Delta_H$:
\[ G^{\ast,M}_{s_1,\ldots,s_l}(\tau) :=m\left( (\mathfrak{g}^{\ast,M} \otimes Z^{\ast})\circ \Delta_H \big( w \big) \right) \in \C[\![q]\!]\,.\]
\end{dfn}

For $s_1,\dots,s_l \geq 2$ the limit
\begin{equation}\label{eq:stufflimit}
  G^\ast_{s_1,\ldots,s_l}(\tau) :=  \lim_{M\to\infty}   G^{\ast,M}_{s_1,\ldots,s_l}(\tau) \, 
\end{equation} 
exists and we have $G_{s_1,\ldots,s_l} = G^\ast_{s_1,\ldots,s_l}=G^\sh_{s_1,\ldots,s_l}$ (\cite[Prop. 6.13]{Ba2}).

\begin{rem} The open question is for what general $s_1,\dots,s_l$ the limit in \eqref{eq:stufflimit} exists. It is believed that this is exactly the case for $s_1\geq2$ and $s_2,\dots,s_l\geq 1$ as explained in Remark 6.14 in \cite{Ba2}. This would be the case if $\Psi_{1,\dots,1}$ are the only multitangent functions with a constant term in the decomposition of Theorem \ref{thm:olivier} iv). That this is the case is remarked, without a proof, in \cite{Bo2} in the last sentence of page 3.  
\end{rem}

\begin{thmx}(\cite{Ba2}) For all $s_1,\dots,s_l \in \N$ and $M \in \N$ the $G_{s_1,\dots,s_l}^{\ast,M} \in \C[\![q]\!]$ have the following properties:
\begin{enumerate}[i)]
\item Their product can be expressed in terms of the stuffle product.
\item In the case where the limit $G^\ast_{s_1,\dots,s_l} :=  \lim_{M\to\infty} G_{s_1,\dots,s_l}^{\ast,M}$ exists, the functions $G^\ast_{s_1,\dots,s_l}$ are elements in $\MZB$. 
\item For $s_1,\dots,s_l \geq 2$ the $G^\ast_{s_1,\dots,s_l}$ exist and equal the classical multiple Eisenstein series
\[ G_{s_1,\dots,s_l}(\tau) = G^\ast_{s_1,\dots,s_l}(\tau) \,.\] 
\end{enumerate}
\end{thmx} 

\subsection{Double shuffle relations for regularized multiple Eisenstein series}

By Theorem \ref{thm:messh} we know that the product of two shuffle regularized multiple Eisenstein series $G^\sh_{s_1,\dots,s_l}$ with $s_1,\dots,s_l \geq 1$ can be expressed by using the shuffle product formula. This means we can for example replace every $\zeta$ by $G^\sh$ in the shuffle product \eqref{eq:shuffle2} of multiple zeta values and obtain
\begin{equation}\label{eq:shmessh}
G^\sh_2 \cdot G^\sh_3 =  G^\sh_{2,3} + 3 G^\sh_{3,2} + 6 G^\sh_{4,1} \,.
\end{equation}
Due to Theorem \ref{thm:messh} iii) we know that $G^\sh_{s_1,\dots,s_l}=G_{s_1,\dots,s_l}$ whenever $s_1,\dots,s_l \geq 2$. Since the product of two multiple Eisenstein series $G_{s_1,\dots,s_l}$ can be expressed using the stuffle product formula we also have
\begin{align}\label{eq:stmessh}
\begin{split}
G^\sh_2 \cdot G^\sh_3 &= G_2 \cdot G_3 = G_{2,3} + G_{3,2} + G_{5} \\
&= G^\sh_{2,3} + G^\sh_{3,2} + G^\sh_{5}   \,.
\end{split}
\end{align}
Combining \eqref{eq:shmessh} and \eqref{eq:stmessh} we obtain the relation $G^\sh_5 = 2 G^\sh_{3,2} + 6 G^\sh_{4,1}$. In the following we will call these relations, i.e. the relations obtained by writing the product of two $G^\sh_{s_1,\dots,s_l}$ with $s_1,\dots,s_l \geq 2$ as the stuffle and shuffle product, \emph{restricted double shuffle relations}.

We know that multiple zeta values fulfill even more linear relations, in particular we can express the product of two multiple zeta values $\zeta(s_1,\dots,s_l)$ in two different ways whenever $s_1\geq 2$ and $s_2,\dots,s_l\geq 1$. A natural question therefore is, in which cases the $G^\sh$ also fulfill these additional relations. The answer to this question is that some are satisfied and some are not, as the following will show. 

In \cite[Example 6.15]{Ba2} it is shown that $G^\sh_{2,1,2} =G^\ast_{2,1,2}$ , $G^\sh_{2,1} =  G^\ast_{2,1}$, $G^\sh_{2,2,1} = G^\ast_{2,2,1}$ and $G^\sh_{4,1} = G^\ast_{4,1}$. Since the product of two $G^\ast$ can be expressed using the stuffle product we obtain
\begin{align}\label{eq:stufflew5}
\begin{split}
G^\sh_2 \cdot G^\sh_{2,1} &= G^\ast_2 \cdot G^\ast_{2,1} \\
&= G^\ast_{2,1,2} + 2 G^\ast_{2,2,1} + G^\ast_{4,1}+ G^\ast_{2,3} \\
&=  G^\sh_{2,1,2} + 2 G^\sh_{2,2,1} + G^\sh_{4,1}+ G^\sh_{2,3} \,.
\end{split}
\end{align}
Using also the shuffle product to express $G^\sh_2 \cdot G^\sh_{2,1}$ we obtain a linear relation in weight $5$ which is not covered by the restricted double shuffle relations. This linear relation was numerically observed in \cite{BT} but could not be proven there. So far it is not known exactly which products of the $G^\sh$ can be written in terms of stuffle products.   \newline 

We end this section by comparing different versions of the double shuffle relations and explain, why multiple Eisenstein series can't fulfill every double shuffle relation of multiple zeta values. For this we write for words $u,v \in \h^1$
\[ \ds(u,v) := u \shuffle v - u \ast v \in \h^1\,.\]
Recall that by $\h^0$ we denote the algebra of all admissible words, i.e. $\h^0 = 1 \cdot \Q + x \h y$. Additionally we set $\h^2=\Q\langle \{z_2,z_3,\dots \} \rangle$ to be the span of all words in $\h^1$ with no $z_1$ occurring, i.e. the words for which the  multiple Eisenstein series $G$ exists. These are also the words for which the product of two multiple Eisenstein series can be expressed as the shuffle and stuffle product by Theorem \ref{thm:messh}. Denote by $|w| \in \h^1$ the length of the word $w$ with respect to the alphabet $\{x , y\}$ and define 
\begin{align*}
\eds_k &:= \big\{ \ds(u,v) \in \h^0 \mid |u|+|v|=k, \,\,  u \in \h^0 , v \in \h^0 \cup \{ z_1 \}  \big\}  \,,\\
\fds_k &:= \big\{ \ds(u,v) \in \h^0 \mid |u|+|v|=k ,\,\, u,v \in \h^0   \big\} \,,\\
\rds_k &:=  \big\{ \ds(u,v) \in \h^0 \mid  |u|+|v|=k ,\,\, u,v \in \h^2 \big\} \,.
\end{align*}
Also set $\eds = \bigcup_{k>0} \eds_k$ and similarly $\fds$ and $\rds$. These spaces can be seen as the words in $\h^0$ corresponding to the extended\footnote{In \cite{IKZ} the authors introduced the notion of extended double shuffle relations. We use this notion here for smaller subset of these relations given there as the relations described in statement (3) on page 315.}-, finite- and the restricted double shuffle relations. We have the inclusions
\[ \rds_k \subset \fds_k \subset \eds_k \,. \] 

View $\zeta$ as a map $\h^0 \rightarrow \MZ$ by sending the word $z_{s_1} \dots z_{s_l}$ to $\zeta(s_1,\dots,s_l)$. It is known (\cite[Thm. 2]{IKZ}), that $\eds_k$ is in the kernel of the map $\zeta$ and it is expected (Statement (3) after Conjecture 1 in \cite{IKZ}) that actually $\eds_k = \ker(\zeta)$. 
Viewing $G^\sh$ in a similar way as a map $\h^0 \rightarrow \MZB$, we know that $\rds_k$ is contained in the kernel of this map (Theorem \ref{thm:messh} iv)). But due to \eqref{eq:stufflew5} we also have $\ds( z_2 , z_2 z_1) \in \ker(G^\sh)$ which is not an element of $\rds_5$. In \cite{Ba} Example 6.15 ii) it is shown that there are also elements in $\fds_k \subset \eds_k$, that are not in the kernel of $G^\sh$. We therefore expect 
\[ \rds \subsetneq \ker G^\sh \subsetneq \eds \]
and the above examples show, that it seems to be crucial to understand for which indices we have $G^\sh = G^\ast$ to answer these questions.

We now discuss applications of the extended double shuffle relations to the classical theory of (quasi-)modular forms. As we have seen in the introduction it is known due to Euler that
\begin{equation}\label{eq:euler3}
 \zeta(2)^2 = \frac{5}{2} \zeta(4)\,,\quad \zeta(4)^2 =  \frac{7}{6} \zeta(8) \,, \quad \zeta(6)^2 = \frac{715}{691} \zeta(12) \,.
\end{equation} 
In the following, we want to show how to prove these relations using extended double shuffle relations and argue why for multiple Eisenstein series the second is fulfilled but the first and the last equation of \eqref{eq:euler3} are not.
\begin{enumerate}[i)] 
\item The relation $\zeta(2)^2 = \frac{5}{2} \zeta(4)$ can be proven in the following way by using double shuffle relations. It is $z_2 \ast z_2 = 2 \ds(z_3,z_1) - \frac{1}{2} \ds(z_2,z_2) + \frac{5}{2} z_4$, since
\begin{align*}
\ds(z_3, z_1) &= z_3 z_1 + z_2 z_2 - z_4\,,\\
\ds(z_2, z_2) &= 4 z_3 z_1 - z_4 \,,\\
z_2 \ast z_2 &= 2 z_2 z_2 + z_4\,.
\end{align*}
Applying the map $\zeta$ we therefore deduce 
\[ \zeta(2)^2 = \zeta(z_2 \ast z_2) = \zeta\left(2 \ds(z_3,z_1) - \frac{1}{2} \ds(z_2,z_2) + \frac{5}{2} z_4\right) = \frac{5}{2}\zeta(4) \,.\]
This relation is not true for Eisenstein series. Though $\ds(z_2,z_2)$ is in the kernel of $G^\sh$ the element $\ds(z_3,z_1)$ is not. In fact, using the explicit formula for the Fourier expansion of $G^\sh_{3,1}$ and $G^\sh_{2,2}$ together with Proposition \ref{prop:formularfordk} for $\dif[2]$ we obtain $G^\sh(\ds(z_3,z_1)) = 6 \zeta(2) \dif G_2$, where as before $\dif = q \frac{d}{dq}$. Using this we get
 \[ G_2^2 = G^\sh(z_2 \ast z_2) = G^\sh\left(2 \ds(z_3,z_1) - \frac{1}{2} \ds(z_2,z_2) + \frac{5}{2} z_4\right) = 12 \zeta(2) \dif G_2 + \frac{5}{2}G_4 \,.\]
This is a well-known fact in the theory of quasi-modular forms (\cite{dz}). 
\item Similarly to the above example one can prove the relation $ \zeta(4)^2 =  \frac{7}{6} \zeta(8)$ by checking that
\begin{align*}
z_4 \ast z_4 = \frac{2}{3} \ds(z_4,z_4) - \frac{1}{2} \ds(z_3,z_5) + \frac{7}{6} z_8 
\end{align*}
and since $\ds(z_4,z_4),  \ds(z_3,z_5) \in \rds_8 \subset \ker G^\sh$ we also derive ${G_4}^2 =  \frac{7}{6} G_8$ by applying the map $G^\sh$ to this equation.

\item To prove the relation $\zeta(6)^2 = \frac{715}{691} \zeta(12)$ in addition to the double shuffles of the form $\ds(z_a,z_b)$ double shuffles of the form $\ds(z_a z_b, z_c)$ are needed as well. This follows indirectly from the results obtained in \cite{GKZ}. Using the computer one can check that 
\begin{align*}
z_6 \ast z_6 = 2 z_6 z_6 + z_{12} =\frac{715}{691} z_{12} + \frac{1}{2^2 \cdot 19 \cdot 113 \cdot 691} \cdot (R + E)
\end{align*}
with $R \in \rds_{12}$ and $E \in \eds_{12} \backslash \rds_{12}$ being the quite complicated elements 
\begin{align*}
R &= \,2005598\ds(z_6,z_6)-8733254\ds(z_7,z_5)+8128450\ds(z_8,z_4)+5121589\ds(z_9,z_3)\\
&+16364863\ds(z_{10},z_2)+2657760\ds(z_2 z_8, z_2)+5220600\ds(z_3 z_7, z_2)\\
&+12711531\ds(z_4 z_6 , z_2)+10460184\ds(z_5 z_5, z_2)+18601119\ds(z_6 z_4, z_2)\\
&+33877826\ds(z_7 z_3 , z_2)+39496002\ds(z_8 z_2, z_2)-13288800\ds(z_2 z_2, z_8)\\
&-5220600\ds(z_2 z_7 , z_3)-5734750\ds(z_3 z_6 , z_3)-84659\ds(z_4 z_5 , z_3)\\
&+2820467\ds(z_5 z_4, z_3)-5486485\ds(z_6 z_3 , z_3)+8462489\ds(z_7 z_2 , z_3)\\
&-6067131\ds(z_2 z_6 , z_4)-7532671\ds(z_3 z_5 , z_4)-10879336\ds(z_4 z_3 , z_5)\\
&-5151234\ds(z_4 z_4, z_4)+3440519\ds(z_5 z_3 , z_4)-1458819\ds(z_6 z_2 , z_4)\\
&+2259096\ds(z_5 z_2, z_5)-4319105\ds(z_3 z_4 , z_5)-778598\ds(z_5 z_2 , z_5)\\
&+7609581\ds(z_2 z_4, z_6)+13064898\ds(z_3 z_3 , z_6)-1281420\ds(z_3 z_2, z_7) \,, \\ \\
E &= -22681134 \ds(z_{11},z_1)+10631040\ds(z_3 z_8, z_1)+4241200\ds(z_7 z_1, z_4) \\
&+31893120\ds(z_4 z_7, z_1)+58185960\ds(z_5 z_6, z_1)+78309000\ds(z_6 z_5, z_1)\\
&+77976780\ds(z_7 z_4, z_1)+44849700\ds(z_8 z_3, z_1)-13288800\ds(z_9 z_2, z_1)\\
&-15946560\ds(z_{10} z_1, z_1)+75052824\ds(z_9 z_1 , z_2)+19477164\ds(z_8 z_1 , z_3)\\
&-12951740\ds(z_6 z_1 , z_5)-10631040\ds(z_2 z_1 , z_9) \,
\end{align*}
Here the elements $E$ and $R$ are in the kernel of $\zeta$ but $E$, in contrast to $R$, is not in the kernel of $G^\sh$. The defect here is given by the cusp form $\Delta$ in weight $12$ as one can derive 
\[ G^\sh(E) = -\frac{2147}{1200} (-2\pi i)^{12} \Delta \,. \] 
\end{enumerate}
It is still an open problem how to derive these Euler relations in general by using double shuffle relations. The last example shows that this also seems to be very complicated.  But as the examples above show, this might be of great interest to understand the connection of modular forms and multiple zeta values. This together with the question which double shuffle relations are fulfilled by multiple Eisenstein series will be considered in upcoming works by the author. 

\section{$q$-analogues of multiple zeta values}\label{section:qana}

In general, a $q$-analogue of an mathematical object is a generalization involving a new parameter $q$ that returns the original object in the limit as $q\rightarrow 1$. The easiest example of such an generalization is the $q$-analogue of a natural number $n \in \N$ given by
\[ [n]_q := \frac{1-q^n}{1-q} = 1 + q + \dots + q^{n-1} \,.\] 
Clearly this gives back the original number $n$ as $\lim_{q \to 1} [n]_q = n$. 

Several different models for $q$-analogues of multiple zeta values have been studied in recent years. A good overview of them can be found in \cite{Zh}. There are different motivations to study $q$-analogues of multiple zeta values. 

That our brackets can be seen as $q$-analogue of multiple zeta values somehow occurred by accident since their original motivation was their appearance in the Fourier expansion of multiple Eisenstein series. But as turned out, seeing them as $q$-analogues gives a direct connection to multiple zeta values.
In this section we first show how the brackets can be seen as a $q$-analogue of multiple zeta values and then discuss how one can obtain relations between multiple zeta values using the results obtained in \cite{BK}. The second section will be devoted to connecting the brackets to other $q$-analogues.

\subsection{Brackets as $q$-analogues of MZV and the map $Z_k$}
Define for $k\in \N$ the map $Z_k: \Q[\![q]\!] \rightarrow \R \cup \{ \infty \}$ by 
\[ Z_k(f) = \lim_{q \to 1} (1-q)^{k}f(q) \,.\]
Since we have seen that the brackets can be written as 
\[[s_1,\dots,s_l] = \frac{1}{(s_1-1)! \dots (s_l-1)!} \sum_{n_1 > \dots > n_l > 0} \prod_{j=1}^l \frac{q^{n_j} P_{s_j-1}\left( q^{n_j} \right)}{(1-q^{n_j})^{s_j}}\] 
and using $P_{k-1}(1) = (k-1)!$ and interchanging the summation and the limit we derive (\cite[Prop. 6.4]{BK}), that for $s_1>1$, i.e. $[s_1,\dots,s_l] \in \MDA$
\[ Z_k\left( [s_1, \dots , s_l ] \right) = \left\{
\begin{array}{cl} \zeta(s_1,\dots,s_l)\,, & k = s_1+\dots+s_l ,  \\ 0\,, &k > s_1+\dots+s_l 
\,. \end{array}   \right. 
\]
Due to $\MD = \MDA[\,[1]\,]$ (Theorem \ref{thm:polyadalg}) we can define a well-defined map\footnote{This map is similar to the evaluation map $Z^\ast: \h^1 \rightarrow \R[T]$, of stuffle regularized multiple zeta values, given in Proposition 1 in \cite{IKZ}. We used this map in the previous sections (Proposition \ref{prop:mzvreg}) with $T=0$.} on the whole space $\MD$ by
\begin{align*}
 Z^{alg}_k:  \filw_{k}(\MD) &\rightarrow \R[T]  \\
 Z^{alg}_k\left( \sum_{j=0}^{k} g_j [1]^{k-j} \right) &= \sum_{j=0}^{k} Z_j(g_j) T^{k-j} \in \R[T]\,  
\end{align*}
where $g_j \in \filw_{j}(\MDA)$.

Every relation between multiple zeta values of weight $k$ is contained in the kernel of the map $Z_k$. Therefore the kernel of $Z_k$ was studied in \cite{BK}. 
\begin{thm}(\cite[Thm. 1.13]{BK})\label{thm_Zk-kernel} For the kernel of $Z^{alg}_k \in \filw_{k}(\MD)$ we have 
\begin{enumerate}[i)]
\item If for $[s_1,\dots,s_l]$ it holds $s_1+\dots+s_l<k$, then 
$Z_k^{alg}[s_1,\dots, s_l]=0$.
\item For any $f \in \filw_{k-2}(\MD) $ we have $Z_k^{alg}\dif(f)=0$, i.e.,  $\dif \filw_{k-2}(\MD)\subseteq \ker Z_k$.
\item If $f \in \filw_{k}(\MD)$ is a cusp form for $\Sl_2(\Z)$, then $Z_k^{alg}(f)=0$.
\end{enumerate}
\end{thm} 

\begin{ex} We illustrate some applications for Theorem \ref{thm_Zk-kernel}. For this we recall identities for the derivatives and relations of brackets as they were given in \cite{BK}. All of them can be obtained by using the results explained in Chapter \ref{section:bracket}. 
\begin{align} \label{example:derivative1}
\dif [1] &= [3] + \frac{1}{2}[2] - [2,1] \,,  \\
\dif[2] &= [4] + 2 [3] - \frac{1}{6} [2] - 4 [3,1] \,,\label{example:derivative2}\\
\dif[2] &= 2[4] + [3] +\frac{1}{6} [2] -2[2,2] - 2[3,1] \,,\label{example:derivative3}\\
\dif [1,1] &=[3,1] + \frac{3}{2}[2,1] + \frac{1}{2}[1,2]+[1,3]- 2 [2,1,1]-[1,2,1] \,, \label{example:derivative4}\\
[8] &= \frac{1}{40} [4] - \frac{1}{252} [2] + 12 [4,4]\,.\label{eq:relwt8}
\end{align}
Using \thmref{thm_Zk-kernel} as immediate consequences 
and without any difficulties we recover the following well-known identities for multiple zeta values.
 \begin{enumerate}[i)]
\item If we apply $Z_3$ to \eqref{example:derivative1} we deduce
$\zeta(3)=\zeta(2,1)$.
\item If we apply  $Z_4$ to \eqref{example:derivative2} and \eqref{example:derivative3} 
we deduce $\zeta(4)= 4 \zeta(3,1) = \frac{4}{3} \zeta(2,2)$.
\item  The identity \eqref{example:derivative4} reads
 in $\MDA[\,[1]\,]$ as
\[ \dif[1,1]= \left([3]-[2,1]+\frac{1}{2}[2]\right)\cdot[1]  
+ 2 [3,1] - \frac{1}{2}[4] - \frac{1}{2}[2,1] - \frac{1}{2}[3] + \frac{1}{3}[2]\,.\,\]
Applying $Z^{alg}_4$ we deduce again the two relations
$\zeta(3)=\zeta(2,1)$ and $4\zeta(3,1)=\zeta(4)$, since
by \thmref{thm_Zk-kernel} we  have 
\[ Z^{alg}_4(\dif[1,1])=\left( \zeta(3) - \zeta(2,1) \right) T - \frac{1}{2}\zeta(4) + 2 \zeta(3,1) = 0 \,. \] 
\item If we apply $Z_8$ to \eqref{eq:relwt8} we deduce $\zeta(8)= 12 \zeta(4,4)$.
\item  As we have seen in Proposition \ref{cor:delta} the cusp form $\Delta$ can be written as
\begin{align}\label{eq:deltal2}
 -\frac{1}{2^6\cdot 5 \cdot 691}  \Delta  &=  168 [5,7]+150 [7,5]+28 [9,3] \notag \\
&+\frac{1}{1408} [2] - \frac{83}{14400}[4] +\frac{187}{6048} [6] - \frac{7}{120} [8] - \frac{5197}{691} [12] \,.
\end{align}
Letting $Z_{12}$ act on both sides of \eqref{eq:deltal2} one obtains the relation \eqref{eq:exotic}
\[ \frac{5197}{691} \zeta(12) =  168 \zeta(5,7)+150 \zeta(7,5) + 28 \zeta(9,3) \,. \] 
\end{enumerate}
\end{ex}

But as mentioned in the introduction there are also elements in the kernel of $Z_k$ that are not covered by Theorem \ref{thm_Zk-kernel}. In weight $4$ one has the following relation of multiple zeta values $\zeta(4) = \zeta(2,1,1)$, i.e. it is $[4] - [2,1,1] \in \ker Z_4$. But this element can't be written as a linear combination of cusp forms, lower weight brackets or derivatives. But using the double shuffle relations for bi-brackets described in Section \ref{sec:dsh} one can prove\footnote{That the last term $\mb{2,1}{1,0}$ in \eqref{eq:4211} is in the kernel of $Z_4$ can be proven in the following way: In Proposition 7.2 \cite{BK} it is shown, that an element $f = \sum_{n>0} a_n q^n$ with $a_n = O(n^m)$ and $m<k-1$ is in the kernel of $Z_k$. Here we have
\begin{align*}
\mb{2,1}{1,0} = \sum_{\substack{u_1 > u_2 > 0\\ v_1 , v_2 > 0}} v_1 u_1 q^{v_1 u_1 + v_2 u_2} < \sum_{\substack{u_1,u_1 0\\ v_1 , v_2 > 0}} v_1 u_1 q^{v_1 u_1 + v_2 u_2} = \dif [1] \cdot [1]\,,
\end{align*}
where the $<$ is meant to be coefficient wise. Since the coefficients of $\dif [1] \cdot [1]$ grow like $n^2 \log(n)^2$ we conclude $\mb{2,1}{1,0} \in \ker Z_4$.} that 
\begin{equation}\label{eq:4211}
[4] - [2,1,1] = \frac{1}{2}\left( \dif[1] + \dif[2] \right) - \frac{1}{3} [2] - [3] + \mb{2,1}{1,0} \,.
\end{equation}

Another way to see that many of the bi-brackets of weight $k$ are in the kernel of the map $Z_k$ is the following. Assume that $s_1 > r_1+1$ and $s_j \geq r_j + 1$ for $j=2, \dots, l$, then using again the representation with the Eulerian polynomials (See also Proposition 1 \cite{Z}) we get
\[ Z_{s_1+\dots+s_l}\left( \mb{s_1,\dots,s_l}{r_1,\dots,r_l} \right) = \frac{1}{r_1! \dots r_l!} \zeta(s_1-r_1,\dots,s_l-r_l) \] 
and in particular with this assumption it is $\mb{s_1,\dots,s_l}{r_1,\dots,r_l} \in \ker Z_{s_1+\dots+s_l+1}$.

The study of the kernel $Z_k$ is of great interest since it contains every relation of weight $k$. We expect that every element in the kernel of $Z_k$ can be described using bi-brackets of a "certain kind" and it seems to be a really interesting question to specify this "certain kind" explicitly. To determine which bi-brackets are exactly in the kernel of the map $Z_k$ and also which bi-brackets can be written in terms of brackets in $\MDA$ is an open problem. The naive guess, that exactly the bi-brackets $\mb{s_1,\dots,s_l}{r_1,\dots,r_l}$ where at least one $r_j >0$ are elements in the kernel of $Z_{s_1+\dots+s_l+r_1+\dots+r_l}$ is wrong, since for example  
\[ \lim_{q \to 1} (1-q)^3 \mb{1,1}{1,0} = \infty \,. \]

\subsection{Connection to other $q$-analogues} \label{sec:otherqana}
In \cite{Zh} the author gives an overview over several different $q$-analogues of multiple zeta values. Here we complement his work and focus on aspects related to our brackets. To compare the brackets to other $q$-analogues we first generalize the notion of a $q$-analogue of multiple zeta values as it was done in \cite{BK2}. 
This notion of a $q$-analogue does cover many but not all $q$-analogues described in \cite{Zh}.

In the following we fix a subset  $S \subset \N$, which we consider as the support for index entries, i.e. we assume $s_1,\dots,s_l \in S$. 
For each $s\in S$ we let  $Q_s(t) \in \Q[t]$ be a polynomial with $Q_s(0)=0$ and $Q_s(1) \neq 0$.
We set $Q = \left\{ Q_s(t) \right\}_{s\in S}$.
A sum of the form 
\begin{equation} \label{eq:zq}
 Z_Q(s_1,\dots,s_l) := \sum_{n_1 > \dots > n_l > 0} \prod_{j=1}^l \frac{Q_{s_j}(q^{n_j})}{(1-q^{n_j})^{s_j}} 
\end{equation}
with polynomials $Q_s$ as before, defines a $q$-analogue of a multiple zeta-value of weight $k=s_1+ \dots +s_l$ and length $l$. Observe only because of $Q_{s_1}(0)=0$ this defines an element of $\Q[\![q]\!]$. That these objects are in fact a $q$-analogue of a multiple zeta-value is justified by the following calculation. 
 \begin{align*}
 \lim\limits_{q \rightarrow 1}{(1-q)^k Z_Q(s_1,\dots,s_l)} &= \sum_{n_1 > \dots > n_l > 0} 
  \prod_{j=1}^l     \lim\limits_{q \rightarrow 1}{ 
  \left( Q_{s_j}(q^{n_j}) \frac{(1-q)^{s_j}}{(1-q^{n_j})^{s_j}} \right)}\\
 &=   Q_{s_1}(1) \dots Q_{s_l}(1) \cdot \zeta(s_1,\dots,s_l) \,.
  \end{align*}
Here we used that $ \lim\limits_{q \rightarrow 1}{ 
   (1-q)^{s}/(1-q^{n})^{s} = 1/n^s }$ and
with the same arguments as in \cite{BK} Proposition 6.4, the above interchange of the limit with the sum can be justified for all 
$(s_1,...,s_l)$ with
$s_1 > 1$. Related definitions for $q$-analogues of multiple zeta values are given in \cite{db}, \cite{yt}, \cite{Zu2} and \cite{YOZ}. 
 It is convenient to define $Z_Q(\emptyset)=1$ and then we denote the vector space spanned by all these elements by
\begin{equation} \label{eq:zqs}
 Z(Q,S) := \big <\,  Z_Q(s_1,\dots,s_l) \big| \,l \ge 0 \mbox{ and } s_1,\dots,s_l \in S \big>_\Q \,. 
\end{equation}
Note by the above convention we have, that $\Q$ is contained in this space. 

\begin{lem} (\cite[Lemma 2.1]{BK2}) \label{lem:algebra}
If for each $r,s \in S$ there exists numbers $\lambda_j(r,s) \in \Q$ 
   such that
\begin{equation}\label{eq:reduction} Q_r(t) \cdot Q_s(t) = \sum_{\substack{j \in S \\1 \le j \le r+s}} \lambda_j(r,s) (1-t)^{r+s-j} Q_j(t) \,, \end{equation}
then the vector space $Z(Q,S)$ is a $\Q$-algebra.
\end{lem} \qed

\begin{thm} (\cite[Thm. 2.4]{BK2}) \label{thm:qanainmdsharp}
Let $Z(Q,\N_{>1})$  be any family of q-analogues of multiple zeta values as in \eqref{eq:zqs}, where 
each $Q_s(t) \in Q$ is a polynomial with degree at most $s-1$, then 
 \[ Z(Q,\N_{>1}) =  \MD^\sharp \,,\]
 where $\MD^\sharp$ was the in Section \ref{sec:dersubalb} defined subalgebra of $\MD$ spanned by all brackets $[s_1,\dots,s_l]$ with $s_j \geq 2$. Therefore, all such  families of q-analogues of multiple zeta values  are $\Q$-subalgebras of $\MD$. 
\qed
\end{thm}
The following proposition allows one to write an arbitrary element in $Z(Q,\N_{>1})$ as an linear combination of $[s_1,\dots,s_l] \in \MD^\sharp $.

\begin{prop}(\cite[Prop. 2.5]{BK2}) \label{prop:qanaintermsofbrackets}  Assume $k\geq 2$.  For
$1 \leq i,j \leq k-1$ define the numbers $b^{k}_{i,j} \in \Q$ by
\[\sum_{j=1}^{k-1} \frac{b^k_{i,j}}{j!} t^j :=\binom{t+k-1-i}{k-1}
\,.\]
With this it is for $1 \leq i \leq k-1$ and $Q^E_{j}(t) = \frac{1}{(j-1)!} t P_j(t)$ 
\[ t^i = \sum_{j=2}^k b^k_{i,j-1} (1-t)^{k-j} Q^E_{j}(t) \,.\] \qed
\end{prop}

We give some examples of $q$-analogues of multiple zeta values, with some being of the above type.\newline

{\bf i)} To write the brackets in the above way we choose $Q^E_s(t) = \frac{1}{(s-1)!} t P_{s-1}(t)$, where the $P_s(t)$ are the Eulerian polynomials defined earlier by
\[ \frac{t P_{s-1}(t)}{(1-t)^{s}} = \sum_{d=1}^\infty d^{s-1} t^d  \]
for $s\geq 0$. With this we have for all $s_1,\dots,s_l \in \N$
\[ [s_1,...,s_l] := 
\sum_{n_1 > ... > n_l > 0} \prod_{j=1}^l \frac{Q_{s_j}^E(q^{n_j})}{(1-q^{n_j})^{s_j}} \,. \]  
and $\MD = Z( \{Q^E_s(t))\}_{s},\N)$.  \newline

{\bf ii)} The polynomials $Q_s^T(t) = t^{s-1}$ are considered in \cite{yt},\cite{Zu2} and sums of the form \eqref{eq:zq} with $s_1>1$ and $s_2, \dots ,s_l \geq 1$ are studied there.  Using Proposition \ref{prop:qanaintermsofbrackets} every $q$-analogue of this type can be written explicitly in terms of brackets.  \newline

{\bf iii)} Okounkov chooses the following polynomials in \cite{O}
\[  Q^O_s(t) = \begin{cases} t^{\frac{s}{2}} & s = 2,4,6,\dots  \\ 
t^{\frac{s-1}{2}} (1+t) & s=3,5,7,\dots . \end{cases} \]
and defines for $s_1,\dots,s_l \in S=\N_{>1}$
\[ \OZ(s) =  \sum_{n_1 > \dots > n_l > 0} \prod_{j=^0}^l \frac{Q^O_{s_j}(q^{n_j})}{(1-q^{n_j})^{s_j}}\,.\] 
We write for the space of the Okounkov $q$-multiple zetas
\[ \qMZV = Z(\{ Q^O_s(t) \}_s ,\N_{>1}) \,. \]  \newline 
Due to Theorem \ref{thm:qanainmdsharp} we have $\qMZV = \MD^\sharp$. In \cite{O} Okounkov conjectures, that the space $\qMZV $ is closed under the operator $\dif$. In length $1$ this is proven in Proposition 2.9 \cite{BK2}. \newline 

{\bf iv)} There are also $q$-analogues which are not of the type as in \eqref{eq:zq}. For example, the model introduced in \cite{YOZ} and further studied in \cite{JMK}. For $s_1,\dots,s_l\geq 1$ they are define by 
\[ \mathfrak{z}_q(s_1,\dots,s_l) = \sum_{n_1 > \dots > n_l > 0} \frac{q^{n_1}}{(1-q^{n_1})^{s_1} \dots (1-q^{n_l})^{s_l}} \,. \]
It is easy to see, that every $\mathfrak{z}_q(s_1,\dots,s_l)$ can be written in terms of bi-brackets. For example 
\begin{align*}
\mathfrak{z}_q(2,1) &= \sum_{n_1 > n_2 > 0} \frac{q^{n_1}}{(1-q^{n_1})^2 (1-q^{n_2})} = \sum_{n_1 > n_2 > 0} \frac{q^{n_1} (q^{n_2} + 1 - q^{n_2})}{(1-q^{n_1})^2 (1-q^{n_2})} \\
&= \sum_{n_1 > n_2 > 0} \frac{q^{n_1} q^{n_2}}{(1-q^{n_1})^2 (1-q^{n_2})} + \sum_{n_1 > n_2 > 0} \frac{q^{n_1}}{(1-q^{n_1})^2} \\
&= [2,1] + \sum_{n_1 > 0} \frac{(n_1-1)q^{n_1}}{(1-q^{n_1})^2} = [2,1] + \mb{2}{1} - [2] \,. 
\end{align*}
Similarly one can prove $\mathfrak{z}_q(2,1,1) = [2,1,1]-2[2,1]+ \mb{2,1}{1,0} + \mb{2}{2} - \frac{3}{2} \mb{2}{1} + [2]$. For higher weights this also works as illustrated in the following 
\begin{align*}
\mathfrak{z}_q(2,2) &= \sum_{n_1 > n_2 > 0} \frac{q^{n_1}}{(1-q^{n_1})^2 (1-q^{n_2})^2} = \sum_{n_1 > n_2 > 0} \frac{q^{n_1}(q^{n_2} + 1 - q^{n_2})}{(1-q^{n_1})^2 (1-q^{n_2})^2} \\ 
&= [2,2] + \mathfrak{z}_q(2,1) = [2,2]+[2,1] + \mb{2}{1} - [2] \,.
\end{align*}
Using again Proposition \ref{prop:qanaintermsofbrackets} it becomes clear for arbitrary weights $s_1,\dots,s_l \geq 2$ we can write $\mathfrak{z}_q(s_1,\dots,s_l)$ in terms of bi-brackets. \newline

Writing any $q$-analogue in terms of bi-brackets enables us to use the double shuffle structure explained in Chapter \ref{section:bibrackets} to obtain linear relations for all of these $q$-analogues. Though it might be difficult to compare our double shuffle relations to the double shuffle relations of other models. For example in the case of $\mathfrak{z}_q(s_1,\dots,s_l)$ the authors in \cite{JMK} consider $s_1,\dots,s_l \in \Z$ to describe their double shuffle relations and it is not clear if these can be written in terms of bi-brackets or not.

{\small
{\it E-mail:}\\\texttt{henrik.bachmann@math.nagoya-u.ac.jp}\\

\noindent {\sc Graduate School of Mathematics\\Nagoya University\\
Chikusa-ku, Nagoya, 464-8602\\
Japan}}

\end{document}

%% file: dia_p0.tex
\begin{tikzpicture}[scale=0.4]
\draw[dotted,step=1,color=gray,thin] (-4.9,-4.9) grid (4.9,4.9); 
\draw [->,thick] (0,-5.5) -- (0,5.5) node (yaxis) [above] {$l$};
\draw [->,thick] (-5.5,0) -- (5.5,0) node (xaxis) [right] {$m$};
\fill[black] (-5,1) circle (4pt);
\fill[black] (-4,1) circle (4pt);
\fill[black] (-3,1) circle (4pt);
\fill[black] (-2,1) circle (4pt);
\fill[black] (-1,1) circle (4pt);
\fill[black] (0,1) circle (4pt);
\fill[black] (1,1) circle (4pt);
\fill[black] (2,1) circle (4pt);
\fill[black] (3,1) circle (4pt);
\fill[black] (4,1) circle (4pt);
\fill[black] (5,1) circle (4pt);
\fill[black] (-5,2) circle (4pt);
\fill[black] (-4,2) circle (4pt);
\fill[black] (-3,2) circle (4pt);
\fill[black] (-2,2) circle (4pt);
\fill[black] (-1,2) circle (4pt);
\fill[black] (0,2) circle (4pt);
\fill[black] (1,2) circle (4pt);
\fill[black] (2,2) circle (4pt);
\fill[black] (3,2) circle (4pt);
\fill[black] (4,2) circle (4pt);
\fill[black] (5,2) circle (4pt);
\fill[black] (-5,3) circle (4pt);
\fill[black] (-4,3) circle (4pt);
\fill[black] (-3,3) circle (4pt);
\fill[black] (-2,3) circle (4pt);
\fill[black] (-1,3) circle (4pt);
\fill[black] (0,3) circle (4pt);
\fill[black] (1,3) circle (4pt);
\fill[black] (2,3) circle (4pt);
\fill[black] (3,3) circle (4pt);
\fill[black] (4,3) circle (4pt);
\fill[black] (5,3) circle (4pt);
\fill[black] (-5,4) circle (4pt);
\fill[black] (-4,4) circle (4pt);
\fill[black] (-3,4) circle (4pt);
\fill[black] (-2,4) circle (4pt);
\fill[black] (-1,4) circle (4pt);
\fill[black] (0,4) circle (4pt);
\fill[black] (1,4) circle (4pt);
\fill[black] (2,4) circle (4pt);
\fill[black] (3,4) circle (4pt);
\fill[black] (4,4) circle (4pt);
\fill[black] (5,4) circle (4pt);
\fill[black] (-5,5) circle (4pt);
\fill[black] (-4,5) circle (4pt);
\fill[black] (-3,5) circle (4pt);
\fill[black] (-2,5) circle (4pt);
\fill[black] (-1,5) circle (4pt);
\fill[black] (0,5) circle (4pt);
\fill[black] (1,5) circle (4pt);
\fill[black] (2,5) circle (4pt);
\fill[black] (3,5) circle (4pt);
\fill[black] (4,5) circle (4pt);
\fill[black] (5,5) circle (4pt);
\fill[black] (1,0) circle (4pt);
\fill[black] (2,0) circle (4pt);
\fill[black] (3,0) circle (4pt);
\fill[black] (4,0) circle (4pt);
\fill[black] (5,0) circle (4pt);
\draw[red,very thin] (1,0.5) arc (90:270:0.5);
\draw[red,very thin] (1,0.5) -- (7,0.5);
\draw[red,very thin] (1,-0.5) -- (7,-0.5);
\coordinate [label=right:\textcolor{red}{$R$}] (R) at (6.5,0);
\draw[blue,very thin] (-6,0.7) -- (7,0.7);
\coordinate [label=right:\textcolor{blue}{$U$}] (R) at (6.5,3);
\end{tikzpicture}

%% file: dia0m1.tex
\begin{tikzpicture}[scale=1.5]
	\def \ra {2}
	\def \noteradius {2.2}
	\def \circsize {0.03}	
		\def \braceradius {2.4}
	\draw (-\ra,0) -- (\ra,0)
				(\ra,0) arc (0:180:\ra);

	\fill (20:\ra) circle (\circsize); 
	\node at (20:\noteradius) {\small $a_8$}; 
				
	\fill (40:\ra) circle (\circsize); 
	\node at (40:\noteradius) {\small $a_7$}; 
	
	\fill(60:\ra) circle (\circsize); 
	\node at (60:\noteradius) {\small $a_6$}; 
	
	\fill(80:\ra) circle (\circsize); 
	\node at (80:\noteradius) {\small $a_5$}; 
	
	\fill (100:\ra) circle (\circsize); 
	\node at (100:\noteradius) {\small $a_4$}; 								

	\fill (120:\ra) circle (\circsize); 
	\node at (120:\noteradius) {\small $a_3$}; 
	
	\fill (140:\ra) circle (\circsize); 
	\node at (140:\noteradius) {\small $a_2$}; 	

	\fill (160:\ra) circle (\circsize); 
	\node at (160:\noteradius) {\small $a_1$}; 	
								
	\fill(-\ra,0) circle (\circsize ); 				
	\draw (-\ra,0) circle (\circsize ); 
	
	\fill (\ra,0) circle (\circsize ); 
	\node at (180:\noteradius) {\small $a_0$}; 	
		\node at (0:\noteradius) {\small $a_9$};

	\draw  (0:\ra) to [bend left] (40:\ra) ;
	\draw  (40:\ra) to [bend left] (100:\ra) ;	
	\draw  (100:\ra) to [bend left] (160:\ra) ;		
	\draw  (160:\ra) to [bend left] (180:\ra) ;		
	
	\braceme[thick]{\braceradius}{0}{40}{br1}{$I(a_7;a_8;a_9)$}
	\braceme[thick]{\braceradius}{40}{100}{br1}{$I(a_4;a_5,a_6;a_7)$}
	\braceme[thick]{\braceradius}{100}{160}{br1}{$I(a_1;a_2,a_3;a_4)$}
  \braceme[thick]{\braceradius}{160}{180}{br1}{$I(a_0;a_1)$}

\end{tikzpicture}

%% file: d32_0.tex
\begin{tikzpicture}[scale=0.6]
	\def \ra {2}
	\def \circsize {0.1}	
	\draw  (0:\ra) to [bend right] (180:\ra) ;			
		
	\draw (-\ra,0) -- (\ra,0)
				(\ra,0) arc (0:180:\ra);
	\draw (-\ra,0) -- (\ra,0)
				(\ra,0) arc (0:180:\ra);
	\fill[black] (180:\ra) circle (\circsize); 
	
	\fill[black] (30:\ra) circle (\circsize); 
	\draw (30:\ra) circle (\circsize); 	
	
	\fill[white] (60:\ra) circle (\circsize); 
	\draw (60:\ra) circle (\circsize); 	
	
	\fill[black] (90:\ra) circle (\circsize); 
	
	\fill[white] (120:\ra) circle (\circsize); 
	\draw (120:\ra) circle (\circsize); 
					
	\fill[white] (150:\ra) circle (\circsize);	
	\draw (150:\ra) circle (\circsize); 
	
	\fill[white] (0:\ra) circle (\circsize); 	
	\draw (0:\ra) circle (\circsize); 
\end{tikzpicture}
  		%

%% file: d32_2.tex
		%
	%
			%
%
\begin{tikzpicture}[scale=0.6]
	\def \ra {2}
	\def \circsize {0.1}

	\draw  (0:\ra) to [bend left] (90:\ra) ;	
	\draw  (90:\ra) to [bend left] (120:\ra) ;	
	\draw  (120:\ra) to [bend left] (150:\ra) ;	
	\draw  (150:\ra) to [bend left] (180:\ra) ;				
			
	\draw (-\ra,0) -- (\ra,0)
				(\ra,0) arc (0:180:\ra);
	\fill[black] (180:\ra) circle (\circsize); 
	
	\fill[black] (30:\ra) circle (\circsize); 
	\draw (30:\ra) circle (\circsize); 	
	
	\fill[white] (60:\ra) circle (\circsize); 
	\draw (60:\ra) circle (\circsize); 	
	
	\fill[black] (90:\ra) circle (\circsize); 
	
	\fill[white] (120:\ra) circle (\circsize); 
	\draw (120:\ra) circle (\circsize); 
					
	\fill[white] (150:\ra) circle (\circsize);	
	\draw (150:\ra) circle (\circsize); 
	
	\fill[white] (0:\ra) circle (\circsize); 	
	\draw (0:\ra) circle (\circsize); 

\end{tikzpicture}

%% file: d32_1.tex
\begin{tikzpicture}[scale=0.6]
	\def \ra {2}
	\def \circsize {0.1}

	\draw  (0:\ra) to [bend left] (30:\ra) ;			
	\draw  (30:\ra) to [bend left] (60:\ra) ;	
	\draw  (60:\ra) to [bend left] (180:\ra) ;				
			
	\draw (-\ra,0) -- (\ra,0)
				(\ra,0) arc (0:180:\ra);
	\fill[black] (180:\ra) circle (\circsize); 
	
	\fill[black] (30:\ra) circle (\circsize); 
	\draw (30:\ra) circle (\circsize); 	
	
	\fill[white] (60:\ra) circle (\circsize); 
	\draw (60:\ra) circle (\circsize); 	
	
	\fill[black] (90:\ra) circle (\circsize); 
	
	\fill[white] (120:\ra) circle (\circsize); 
	\draw (120:\ra) circle (\circsize); 
					
	\fill[white] (150:\ra) circle (\circsize);	
	\draw (150:\ra) circle (\circsize); 
	
	\fill[white] (0:\ra) circle (\circsize); 	
	\draw (0:\ra) circle (\circsize); 

\end{tikzpicture}
\,\,\,\,\,\,\,\,
\begin{tikzpicture}[scale=0.6]
	\def \ra {2}
	\def \circsize {0.1}

	\draw  (0:\ra) to [bend left] (30:\ra) ;			
	\draw  (30:\ra) to [bend left] (150:\ra) ;	
	\draw  (150:\ra) to [bend left] (180:\ra) ;				
			
	\draw (-\ra,0) -- (\ra,0)
				(\ra,0) arc (0:180:\ra);
	\fill[black] (180:\ra) circle (\circsize); 
	
	\fill[black] (30:\ra) circle (\circsize); 
	\draw (30:\ra) circle (\circsize); 	
	
	\fill[white] (60:\ra) circle (\circsize); 
	\draw (60:\ra) circle (\circsize); 	
	
	\fill[black] (90:\ra) circle (\circsize); 
	
	\fill[white] (120:\ra) circle (\circsize); 
	\draw (120:\ra) circle (\circsize); 
					
	\fill[white] (150:\ra) circle (\circsize);	
	\draw (150:\ra) circle (\circsize); 
	
	\fill[white] (0:\ra) circle (\circsize); 	
	\draw (0:\ra) circle (\circsize); 

\end{tikzpicture}
\,\,\,\,\,\,\,\,
\begin{tikzpicture}[scale=0.6]
	\def \ra {2}
	\def \circsize {0.1}

	\draw  (0:\ra) to [bend left] (30:\ra) ;	
	\draw  (30:\ra) to [bend left] (120:\ra) ;	
	\draw  (120:\ra) to [bend left] (150:\ra) ;	
	\draw  (150:\ra) to [bend left] (180:\ra) ;				
			
	\draw (-\ra,0) -- (\ra,0)
				(\ra,0) arc (0:180:\ra);
	\fill[black] (180:\ra) circle (\circsize); 
	
	\fill[black] (30:\ra) circle (\circsize); 
	\draw (30:\ra) circle (\circsize); 	
	
	\fill[white] (60:\ra) circle (\circsize); 
	\draw (60:\ra) circle (\circsize); 	
	
	\fill[black] (90:\ra) circle (\circsize); 
	
	\fill[white] (120:\ra) circle (\circsize); 
	\draw (120:\ra) circle (\circsize); 
					
	\fill[white] (150:\ra) circle (\circsize);	
	\draw (150:\ra) circle (\circsize); 
	
	\fill[white] (0:\ra) circle (\circsize); 	
	\draw (0:\ra) circle (\circsize); 

\end{tikzpicture}
	%
			%
%

%% file: d32_12.tex
\begin{tikzpicture}[scale=0.6]
	\def \ra {2}
	\def \circsize {0.1}	
			
	\draw  (0:\ra) to [bend left] (30:\ra) ;	
	\draw  (30:\ra) to [bend left] (60:\ra) ;	
	\draw  (60:\ra) to [bend left] (90:\ra) ;	
	\draw  (90:\ra) to [bend left] (120:\ra) ;				
	\draw  (120:\ra) to [bend left] (150:\ra) ;	
	\draw  (150:\ra) to [bend left] (180:\ra) ;					
	
	\draw (-\ra,0) -- (\ra,0)
				(\ra,0) arc (0:180:\ra);
	\fill[black] (180:\ra) circle (\circsize); 
	
	\fill[black] (30:\ra) circle (\circsize); 
	\draw (30:\ra) circle (\circsize); 	
	
	\fill[white] (60:\ra) circle (\circsize); 
	\draw (60:\ra) circle (\circsize); 	
	
	\fill[black] (90:\ra) circle (\circsize); 
	
	\fill[white] (120:\ra) circle (\circsize); 
	\draw (120:\ra) circle (\circsize); 
					
	\fill[white] (150:\ra) circle (\circsize);	
	\draw (150:\ra) circle (\circsize); 
	
	\fill[white] (0:\ra) circle (\circsize); 	
	\draw (0:\ra) circle (\circsize); 

\end{tikzpicture}